\documentclass[11pt]{article}
\usepackage{amsbsy,amsfonts,amsmath,amssymb,amsthm,euscript,enumerate}
\usepackage{graphicx,color,centernot,framed}
\definecolor{shadecolor}{gray}{0.85}
\usepackage{url}
\usepackage{multirow}
\usepackage{hyperref}
\usepackage{siunitx}
\usepackage{mathrsfs}
\hypersetup{colorlinks,citecolor=blue,linkcolor=blue,urlcolor=black}
\usepackage{bbm}

\newcommand{\indic}[1]{\mathbbm{1}{\raisebox{-2pt}{$\scriptstyle #1$}}}

\allowdisplaybreaks[1]

\theoremstyle{plain}
\numberwithin{equation}{section}
\theoremstyle{plain}
\newtheorem{thm}{Theorem}[section]

\newtheorem{lmm}[thm]{Lemma}
\newtheorem{prp}[thm]{Proposition}

\theoremstyle{definition}
\newtheorem{dfn}{Definition}
\newtheorem{rmk}{Remark}

\newcommand{\lbeq}[1]{\label{eq:#1}}

\newcommand{\nn}{\nonumber}

\newcommand{\Proof}[1]{\paragraph{\it #1}}
\newcommand{\QED}{\hspace*{\fill}\rule{7pt}{7pt}\bigskip}

\newcommand{\refeq}[1]{(\ref{eq:#1})}

\newcommand{\sss}{\scriptscriptstyle}

\newcommand{\vep}{\varepsilon}

\newcommand{\Z}{\mathbb{Z}}

\newcommand{\Zd}{\Z^d}

\newcommand{\bs}{\boldsymbol}

\makeatletter
\newcommand{\xRightarrow}[2][]{\ext@arrow 0359\Rightarrowfill@{#1}{#2}}
\makeatother

\definecolor{dgreen}{rgb}{0,0.6,0}

\title{Asymptotic expansion of the critical point for oriented percolation in high dimensions. }
\author{Noe Kawamoto\thanks{nkawamoto@ncts.ntu.edu.tw}\footnote{National Center for Theoretical Sciences, Taiwan. \url{https://orcid.org/0000-0003-2273-4497}}
}
\date{}
\begin{document}
\maketitle

\begin{abstract}
We study an asymptotic expansion of the critical point for the nearest-neighbor oriented percolation on $\Zd$ in powers of $d^{-1}$ as $d\rightarrow \infty$. The proof relies heavily on the lace expansion. 
\end{abstract}
%\tableofcontents
\section{Introduction}
In this paper, we study oriented percolation on $\Zd$. Let $\Lambda=(\Zd, \mathbb Z_{+})$. For $\bs{x}=(x,t), \bs{y}=(y,t+1) \in \Lambda$ with $x \neq y$, we define a bond $[\bs{x},\bs{y}\rangle$, which is open  with probability $pD(y-x)$, independently of all other bonds. Here, $D$ is the $1$-step transition probability of the underlying simple random walk, defined by
\begin{align}
D(x)=\frac{1}{2d}\indic{\{|x|=1\}}. 
\end{align}
We define $\theta_p$ as the probability that there exist an occupied path starting at $\bs{o}=(o,0)$ and diverging to infty. The critical point $p_c$ is then defined as 
\begin{align}
p_c:=\inf\{p: \theta_p>0 \}.
\end{align}
The purpose of this paper is to examine asymptotic expansion of $p_c$ in powers of $d^{-1}$ as $d\rightarrow \infty$. Our main result is stated below.  
\begin{thm}\label{thm:mainthm}
As $d \rightarrow \infty$, 
\begin{align}\lbeq{mainthm}
p_c=1+(2d)^{-2}+\frac{7}{2}(2d)^{-3}+\frac{129}{8}(2d)^{-4}+O(d^{-5}).
\end{align}
\end{thm}
The proof relies on the lace expansion derived by Sakai~\cite{S01}, which will be introduced in Section~\ref{LaceExpansion}. 
\begin{rmk}[Application to the contact process.]
The contact process is defined by the zero lattice-spacing of the time-discretized contact process with a discretization parameter $\varepsilon \in (0,1]$, which is the following oriented percolation on $\Zd \times \varepsilon\mathbb Z_{+}$: For $\bs{x}=(x,t), \bs{y}=(y,t+\varepsilon) \in (\Zd, \varepsilon\mathbb Z_{+})$, a bond $[\bs{x},\bs{y}\rangle$ is open independently of the other bonds with probability
\begin{align}
q_p(y-x)=
\begin{cases}
1-\vep&(x=y),\\
p\vep D(y-x)&(x\neq y). 
\end{cases}
\end{align}
The case $\varepsilon =1$ corresponds to standard oriented percolation. In the time-discretized contact process, we also have a bond from each point to itself after a time increment of $\varepsilon$. We believe that the methods introduced in this paper can be adapted to the time-discretized contact process, enabling the estimation of its critical point. 
\end{rmk}

{\bf Earlier findings on the critical point of the oriented percolation.} In \cite{B77}, Blease considers the nearest-neighbor percolation on a graph that is oriented in both time and space, which differs from our setting, and derives a non-rigorous expansion of $p_c$ for large $d$: 
\begin{align}\lbeq{Blease}
p_c=1+\frac{1}{2}d^{-2}+d^{-3}+3d^{-4}+\frac{21}{2}d^{-5}+\frac{479}{12}d^{-6}+O(d^{-8}).
\end{align}
Cox and Durret \cite{CD83} deal with the same setting considered by Blease, and obtain the upper and lower bound 
\begin{align}
1+\frac{1}{2}d^{-2}+o(d^{-2})\le p_c \le 1+d^{-2}+O(d^{-3}),
\end{align}
where the lower bound agrees with \refeq{Blease} up to terms of order $o(d^{-2})$. 
Their proof does not rely on the lace expansion: they establish a comparison between oriented percoaltion and a branching process. Due to the difference in setting, applying the methods used in this paper may yield a different value of $p_c$, as the value is  generally believed to depend on the details of the underlying graph. However, the method itself can be adapted to their setting. 

For the spread out model, in which $D$ is defined as the uniform distribution on a box
of side length $L$, a result was obtained by van der Hofstad and Sakai \cite{HS05}, namely
\begin{align}
p_c=1+\frac{1}{2}L^{-d}\sum_{n=2}^{\infty}U^{*2n}(o)+O(L^{-d-1}),
\end{align}
where $U(x)$ is the uniform distribution on a box of side length $1$ on $\mathbb R^d$ and $U^{*2n}$ denotes the $2n$-fold convolution of $U$ itself. In their analysis, the second biggest error comes from the Riemann sum approximation. They first derive an expression in terms of $D$ and then they estimate the error between $D$ and $U$. 
Therefore, terms of order higher than $L^{-d-1}$ are absorbed in the error from the
Riemann sum approximation, allowing them to focus on deriving the coefficient of the leading error term, which simplifies the proof. Although they also use the lace expansion, in contrast to the spread-out case, the nearest neighbor case requires more careful estimates, as the error terms beyond the second order cannot be neglected.  

{\bf Toward all-order asymptotic expansion of the critical point in oriented percolation.} 
In this paper, we investigate the expansion of $p_c$ up to order $d^{-4}$, providing an approximation of its value. It is natural to ask whether an expansion to all orders exists for $p_c$. In \cite{HGS05} and \cite{HS95}, van der Hofstad and Slade, and Hara and Slade, respectively, proved for bond percolation and the self-avoiding walk, there is an asymptotic expansion of the critical point in powers of $1/2d$ to all orders. Specifically, they proved that there are constants $\{\alpha_i\}_i$ such that for all $M\ge 1$, 
\begin{align}\lbeq{pcapprox}
p_c=\sum_{i=0}^{M-1}\alpha_i(2d)^{-i}+O(d^{-M}). 
\end{align}
Here, the constant in the error term may depend on $M$. In the case of oriented percolation, a corresponding result is expected to hold. However, it is expected that the full asymptotic series $\sum_{i}\alpha_i(2d)^{-i}$ does not converge for any $d$: that is, the series $\sum_{i}\alpha_iz^{-i}$ has zero radius of convergence. There exists an example in which the divergence of the full asymptotic series has been rigorously proven. For the spherical model, defined by the $N\rightarrow \infty$ limit of the $N$-vector model, it was proven in \cite{GF74} that the critical temperature admits an asymptotic expansion to all orders in powers of $1/d$, although the series is divergent. Furthermore, for the spherical model, sign changes appear in the expansion of the critical temperature \cite{GF74}. The value of the critical point generally depends on the specific details of the model, so this does not imply that the critical points in other models share the same property. However, there is no clear reasoning to suggest that these properties are unique to the spherical model. 

In \cite{S79}, Sokal proved the following: Suppose $p_c(z)$ can be expressed as in \refeq{pcapprox} with $d$ replaced by $z$, and that $p_c(z)$ is an analytic on the domain $D=\{z: \operatorname{Re} z^{-1} > 1\}$ and satisfies, for all $s \in D$ and all $M\ge 1$, 
\begin{align}\lbeq{borelbd}
\Big|p_c(z)-\sum_{i=0}^{M-1}\alpha_iz^i\Big|\le C_1^Mz^MM!, 
\end{align}
then, the Borel transform $B(z):=\sum_{i=0}^{\infty}\frac{\alpha_i}{i!}z^i$ converges and $p_c(z)$ can be represented as 
\begin{align}
p_c(z)=\frac{1}{z}\int_0^{\infty}B(\sigma)e^{-\frac{\sigma}{z}}d\sigma~~(z \in D_R). 
\end{align}
This means that the actual function of $p_c(z)$ can be recovered from its asymptotic expansion via Borel summation. Graham \cite{G10} proved the bound in \refeq{borelbd} for the self-avoiding walk, and a corresponding result is expected to hold for oriented percolation. However, the analyticity of $p_c(z)$ remains a conjecture. 

\subsection{Organization of the paper}
The remainder of this paper is organized as follows. In Section~\ref{LaceExpansion}, we introduce the lace expansion for oriented percolation, and prove the main result, Theorem~\ref{thm:mainthm}, assuming Proposition~\ref{prp:sumofeachpi}. Proposition~\ref{prp:sumofeachpi} will be proved across Section~\ref{SectionPi0} (estimating $\Pi_{p}^{\sss(0)}$), Section~\ref{Sectionpi1} (estimating $\Pi_{p}^{\sss(1)}$) and Section~\ref{Sectionpi2} (estimating $\Pi_{p}^{\sss(2)}$). Section~\ref{PreliminaryPi0} is devoted to preliminary estimates required in Section~\ref{SectionPi0}. Finally, the appendix contains proofs of two lemmas that are used in Section~\ref{3case}. 

{\bf Notations.}
\begin{itemize}
\item 
We use the Fourier transform and the inverse Fourier transform, defined for an absolutely summable finctions $f: \Lambda \rightarrow \mathbb C$ by
\begin{align}
\hat f(k,\omega)=\sum_{\bs{x} \in \Lambda}f(\bs{x})e^{i(k,\omega)\cdot \bs{x}},&& f(\bs{x})=\int_{[-\pi,\pi]^{d+1}}\hat f(k,\omega)e^{-i(k,\omega)\cdot \bs{x}}\frac{d\omega}{2\pi}\frac{d^dk}{(2\pi)^d}
\end{align}
where $\omega \in [-\pi,\pi]$ and $k \in [-\pi,\pi]^d$. 

\item The convolution of two absolutely summable functions $f,g: \Lambda \rightarrow \mathbb C$ is defined for $\bs{x} \in \Lambda$ as 
\begin{align}
(f*g)(\bs{x})=\sum_{\bs{y}: \bs{y}(t) \le \bs{x}(t)}f(\bs{y})g(\bs{x}-\bs{y}). 
\end{align}
We denote by $f^{*m}$ the $m$-fold convolution of $f$ for $m\ge 2$. 
\item We use $C$ and $C'$ to denote positive constants that may depend on $d$. Their values may change from line to line. 
\item For $\bs{x} \in \Lambda$, we denote the spatial component of $\bs{x}$ by $\bs{x}(s)$ and the temporal component of $\bs{x}$ by $\bs{x}(t)$, i.e., $\bs{x}=(\bs{x}(s),\bs{x}(t))$. For $\bs{x}$ and $\bs{y}$ such that $\bs{x}(t)=\bs{y}(t)+1$, we write $D(\bs{x}-\bs{y})$ as shorthand for $D(\bs{x}(s)-\bs{y}(s))$. 
\end{itemize}

\section{The lace expansion for oriented percolation}\label{LaceExpansion}
In this section, we introduce the lace expansion for oriented percolation and conclude the section with a proof of Theorem~\ref{thm:mainthm}. Our analysis of $p_c$ relies heavily on the lace expansion, which was introduced by Sakai~\cite{S01}. To define the lace-expansion coefficients, we first introduce some essential definitons of key connection events. 

 We say that $\bs{x}$ is connected to $\bs{y}$ and write $\bs{x} \rightarrow \bs{y}$.
For events $E_1$ and $E_2$, we write $E_1 \circ E_2$ to denote the event that $E_1$ and $E_2$ occur bond-disjointly in terms of bonds, that is, they occur using disjoint sets of bonds. We denote by $\text{piv}E_1$ the set of pivotal bonds for $E_1$. For the definition of a pivotal bond, see, for example, \cite{ms93}. 
\begin{dfn} 
\begin{itemize}
\item We define $\{\bs{o}\rightarrow \bs{x}\big\}\circ \{\bs{o}\rightarrow \bs{x}\big\}$ as $\{\bs{o}\Rightarrow \bs{x}\big\}$, and say that $\bs{o}$ and $\bs{x}$ are doubly connected. 

\item Let $S\subset \Zd$. Define $E(b,\bs{x}: S)$ as 
\begin{align}
E(b,\bs{x}: S)=\{b \rightarrow \bs{x}\}\cap \{ \bs{x} \in S\} \cap \{\forall b' \in \text{piv}\{ \bs{\bar{b}} \rightarrow \bs{x}\}\}. 
\end{align}
\item Define $\tilde C^{b}(\bs{u})$ to be the set of vertices connected to $\bs{u}$ without using $b$. 
\end{itemize}
\end{dfn}
We can now define the lace expansion coefficient. 
\begin{dfn}[Lace-expansion coefficients]Let $n \ge 1$, and $p\le p_c$. We define 
\begin{align}
\Pi_{p}^{\sss(0)}(\bs{x})&:=\mathbb P_p\big(\bs{o}\Rightarrow \bs{x}\big),\lbeq{defpi0}\\
\Pi_{p}^{\sss(n)}(\bs{x})&:=\sum_{b_1, \cdots, b_n \in \mathbb B}\mathbb P_p\Big(\{\bs{o}\Rightarrow \bs{\underline{b}_1}\} \cap \bigcap_{j=1}^{n} E\big(b_j, \bs{\underline{b}_{j+1}}: \tilde C^{b_j}(\bs{\bar{b}_{j-1}})\big)\Big).\lbeq{defpin}
\end{align}
Furthermore, we define 
\begin{align}
\Pi_{p}(\bs{x})=\sum_{n=0}^{\infty}(-1)^n\Pi_{p}^{\sss(n)}(\bs{x}).
\end{align}
\end{dfn}
The lace expansion yields an recursion equation for the two point function: 
\begin{align}
\tau_p(\bs{x}):=\mathbb P_p\big(\bs{o}\rightarrow \bs{x}\big). 
\end{align}
\begin{thm}[Sakai, \cite{S01}]The two-point function satisfies following equation: 
\begin{align}\lbeq{LE}
\tau_p(\bs{x})=\delta_{\bs{o},\bs{x}}+\Pi_{p}(\bs{x})+p\sum_{\bs{u},\bs{v}\in \Lambda}[\delta_{\bs{o},\bs{u}}+\Pi_{p}(\bs{u})]D(\bs{v}-\bs{u})\tau_p(\bs{x}-\bs{v}). 
\end{align}
Furthermore, for sufficiently large $d>4$, there exist constants $C$ and $C'$ such that for $p\le p_c$, 
\begin{align}\lbeq{Pinbd}
\big|\sum_{\bs{x}\in \Lambda}\Pi_{p}^{\sss(n)}(\bs{x}) \big|\le
\begin{cases}
Cd^{-1}&(n=0),\\
(C'd)^{-n}&(n\ge 1). 
\end{cases}
\end{align}
\end{thm}
The following bound is known as the infrared bound. 
\begin{thm}[Nguyen and Yang~\cite{NY93}]\label{lmm: infraredbd}
Let $d$ be sufficiently large, greater than $4$ and $p <p_c$. There exists a constant $C$ such that 
\begin{align}
|\hat \tau_p(k,\omega)| \le \frac{C}{|k|^2+|\omega|}. 
\end{align}
\end{thm}

\begin{lmm}\label{lmm:errorm}
Let $0\le n\le4 $ and $p\le p_c$. Then,
\begin{align}
\sum _{\substack{\bs{x}\\ \bs{x}(t)\ge 5}}\Pi_{p_c}^{\sss(n)}(\bs{x})\le O(d^{-5}).
\end{align}
\end{lmm}

\Proof{Proof of Lemma~\ref{lmm:errorm}.}
Let $p<p_c$. First we consider the case $n=0$. By the definition given in \refeq{defpi0} and applying the BK inequality, we obtain
\begin{align}\lbeq{pi0bdbk}
\sum_{\bs{x}: \bs{x}(t)\ge 5}\Pi_{p}^{\sss(0)}(\bs{x}) \le \sum_{\bs{x}: \bs{x}(t)\ge 5}\tau_p(\bs{x})^2\le p_c^{10}\sum_{\bs{x}}(D^{*5}*\tau_p)(\bs{x})^2. 
\end{align}
By the inverse Fourier transform and Lemma~\ref{lmm: infraredbd}, \refeq{pi0bdbk} can be bounded as
\begin{align}\lbeq{pi0frbd}
\int_{[-\pi,\pi]^{d+1}}\hat D(k)^{10}\hat \tau_p(k,\omega)\hat \tau_p(-k,-\omega)\frac{d\omega}{2\pi}\frac{d^dk}{(2\pi)^d}&\le \int_{[-\pi,\pi]^{d+1}}\frac{\hat D(k)^{10}}{[|k|^2+|\omega|]^2}\frac{d\omega}{2\pi}\frac{d^dk}{(2\pi)^d}\nn\\
&\le \int_{[-\pi,\pi]^{d}}\frac{\hat D(k)^{10}}{|k|^2}\frac{d^dk}{(2\pi)^d}\nn\\
&\le C\int_{[-\pi,\pi]^{d}} \frac{\hat D(k)^{10}}{1-\hat D(k)} \frac{d^dk}{(2\pi)^d}. 
\end{align}
For the last inequalitiy, we use the bound $1-\hat D(k)\le \frac{1}{2}|k|^2$ which follows from $\hat D(k)=\frac{1}{d}\sum_{j=1}^d\cos k_j$. By \cite[Appendix]{HS95}, \refeq{pi0frbd} can be bounded by $O(d^{-5})$. 

Next, we consider the case $1\le n \le 4$. We note that for $i \ge 2$, $\bs{\underline{b}_i}\neq \bs{\bar{b}_{i-1}}$. Thus,  $\Pi_{p}^{\sss(n)}(\bs{x})=0$ for $\bs{x}$ such that $\bs{x}(t)<2n-1$. In the remainder of this proof, we use the following diagrammatic representation: 
\begin{align}
\tau_p(\bs{x})=\raisebox{0pt}{\includegraphics[scale=0.6]{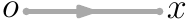}}, &&D(\bs{\bar{b}_j}-\bs{\underline{b}_j})=\raisebox{0pt}{\includegraphics[scale=0.7]{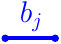}}. 
\end{align}
For example, using the above diagrammatic representation, the following can be expressed as 
\begin{align}
D(\bs{\bar{b}_j}-\bs{\underline{b}_j})\tau_p(\bs{y_{j+1}}-\bs{\bar{b}_j})\tau_p(\bs{\underline{b}_{j+1}}-\bs{y_{j+1}})\tau_p(\bs{\underline{b}_{j+1}}-\bs{y_{j}})=\raisebox{-20pt}{\includegraphics[scale=0.55]{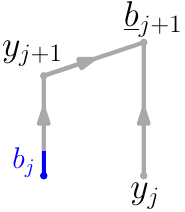}}.
\end{align}
By \cite[Lemma~5.3]{S01} and applying the BK inequality, we obtain 
\begin{align}
\Pi_{p}^{\sss(n)}(\bs{x})
&\le\sum_{b_1, \cdots, b_n \in \mathbb B}\sum_{\substack{\bs{y_1},\cdots, \bs{y_n} \\ \bs{y_i}(t)\le\bs{\underline{b}_{i}} (t)}}\Bigg\{\delta_{\bs{o}, \bs{y_1}}\bigg[\delta_{\bs{o},\bs{\underline{b}_1}}+[1-\delta_{\bs{o},\bs{\underline{b}_1}}]\raisebox{-20pt}{\includegraphics[scale=0.55]{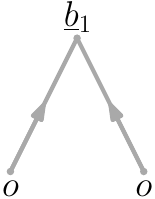}}\bigg]+[1-\delta_{\bs{o},\bs{\underline{b}_1}}]\raisebox{-20pt}{\includegraphics[scale=0.55]{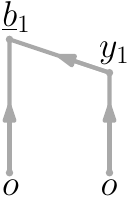}}\bigg\}\nn\\
&\times \prod_{j=1}^{n-1}\Bigg\{\raisebox{-20pt}{\includegraphics[scale=0.55]{Pindecj1}}+\raisebox{-20pt}{\includegraphics[scale=0.55]{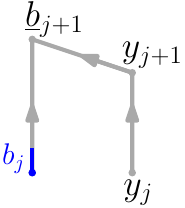}}\Bigg\}\raisebox{-20pt}{\includegraphics[scale=0.55]{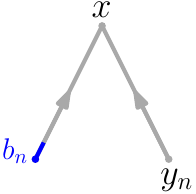}}. 
\end{align}
By considering the difference of the time components between $\bs{\bar{b}_i}$ and $\bs{\underline{b}_{i+1}}$ for $0\le i\le n$, where $\bs{\bar{b}_0}=\bs{o}$ and $\bs{\underline{b}_{n+1}}=\bs{x}$, the summation can be decomposed as
\begin{align}
\sum_{\bs{x}: \bs{x}(t)\ge 5}\sum_{b_1, \cdots, b_n \in \mathbb B}=\sum_{m=5\vee2n-1 }^{\infty}\sum_{\bs{x}: \bs{x}(t)=m}\sum_{\substack{t_1, \cdots, t_{n+1} \\ \sum_{i=1}^{n+1}t_i=m-n \\ t_i\ge 1~(i\ge 2)}}\sum_{\substack{b_1, \cdots, b_n \in \mathbb B \\ \bs{\underline{b}_i}=\sum_{l=1}^{i}t_l+i-1}}. 
\end{align}
Therefore, we obtain 
\begin{align}\lbeq{pinupperbd}
&\sum_{\bs{x}: \bs{x}(t)\ge 5}\Pi_{p}^{\sss(n)}(\bs{x})\nn\\
&\le \sum_{m=5\vee2n-1 }^{\infty}\sum_{\substack{t_1, \cdots, t_{n+1} \\ \sum_{i=1}^{n+1}t_i=m-n \\ t_i\ge 1~(i\ge 2)}}[\delta_{0,t_1}+T_{1,p}(t_1)+T_{2,p}(t_1)]\prod_{j=1}^{n-1}2T_{2,p}(t_{j+1}+1)\cdot T_{1,p}(t_{n+1}+1)
\end{align}
where, for $l\ge1$ and $m=1,2$, 
\begin{align}
T_{m,p}(l)= \sup_{\bs{y},\bs{z}}\sum_{\substack{\bs{x} \\ \bs{x}(t)\ge \bs{y}(t), \bs{z}(t) }}p^{2l}(D^{*l}*\tau_p)(\bs{x}-\bs{y})(D^{*l}*\tau_p^{*m})(\bs{x}-\bs{z}). 
\end{align}
As in \refeq{pi0frbd}, by applying \cite[Appendix]{HS95}, we obtain $T_{1,p}(l)\le O(d^{-l})$ and $T_{2,p}(l)\le O(d^{-l})$. Therefore, \refeq{pinupperbd} can be bounded as 
\begin{align}
\sum_{\bs{x}: \bs{x}(t)\ge 5}\Pi_{p}^{\sss(n)}(\bs{x})\le \sum_{m=5\vee2n-1 }^{\infty}\sum_{\substack{t_1, \cdots, t_{n+1} \\ \sum_{i=1}^{n+1}t_i=m-n \\ t_i\ge 1~(i\ge 2)}}(Cd)^{-m}\le \sum_{m=5\vee2n-1 }^{\infty}m^4(C'd)^{-m}\le O(d^{-5}). 
\end{align}
For the case where $p=p_c$, since $\Pi_{p}^{\sss(n)}(\bs{x})$ is increasing in $p$, we take the limit as $p \rightarrow p_c$. Then, we complete the proof. 
\QED

Before concluding this section, we prove Theorem~\ref{thm:mainthm}, assuming the following lemma, whose proof will occupy the remaining sections. The proof of \refeq{abc}, \refeq{def} and \refeq{ghi} will be given in Sections~\ref{SectionPi0}, \ref{Sectionpi1} and \ref{Sectionpi2}, respectively. 
\begin{prp}\label{prp:sumofeachpi}
Let $s:=\frac{1}{2d}$. 
\begin{align}
\sum_{\substack{\bs{x}\\ \bs{x}(t)=2,3,4}}\Pi_{p_c}^{\sss(0)}(\bs{x})&=p_c^4s^2+\frac{3p_c^4}{2}[3p_c^2-1]s^3-\frac{111}{8}s^4+O(d^{-5}),\lbeq{abc}\\
 \sum_{\substack{\bs{x}\\ \bs{x}(t)=2,3,4}}\Pi_{p_c}^{\sss(1)}(\bs{x})&=2p_c^4s^2+\frac{p_c^4}{2}[19p_c^2-6]s^3+29s^4+O(d^{-5}),\lbeq{def} \\
\sum_{\substack{\bs{x}\\ \bs{x}(t)=3,4}}\Pi_{p_c}^{\sss(2)}(\bs{x})&=4s^4+O(d^{-5}). \lbeq{ghi}
\end{align}
\end{prp}
\Proof{Proof of Theorem~\ref{thm:mainthm}, assuming Proposition~\ref{prp:sumofeachpi}.}
By summing both sides of \refeq{LE} over $\bs{x}$, and solving the equation for $\chi_p$, we obtain
\begin{align}
\chi_p=\frac{1+\sum_{{\bs{x}}}\Pi_p({\bs{x}})}{1-p-p\sum_{{\bs{x}}}\Pi_p({\bs{x}})}.
\end{align}
The critical point $p_c$ is defined as the velue of $p$ at which $\chi_p$ diverges. By \refeq{Pinbd}, the numerator of $\chi_p$ remains finite as $p \rightarrow p_c$. This leads to an identity of $p_c$: 
\begin{align}\lbeq{pcidentity}
p_c=1-p_c\sum_{\bs{x}}\Pi_{p_c}(\bs{x}).
\end{align}
We note that as a consequence of \refeq{Pinbd}, 
\begin{align}\lbeq{easybdpc}
|p_c-1|\le Cd^{-1}
\end{align}
which will be used frequently throughout the paper. 

By \refeq{Pinbd}, \refeq{pcidentity} and Lemma~\ref{lmm:errorm}, we obtain
\begin{align}\lbeq{recursionE}
p_c&=1-p_c\left\{\sum_{\substack{\bs{x}\\ \bs{x}(t)=2,3,4}}\Pi_{p_c}^{\sss(0)}(\bs{x})-\sum_{\substack{\bs{x}\\ \bs{x}(t)=2,3,4}}\Pi_{p_c}^{\sss(1)}(\bs{x})+\sum_{\substack{\bs{x}\\ \bs{x}(t)=3,4}}\Pi_{p_c}^{\sss(2)}(\bs{x})\right\}+O(d^{-5}).
\end{align}
By \refeq{recursionE} and Proposition~\ref{prp:sumofeachpi}, we obtain 
\begin{align}\lbeq{kk}
p_c=1+p_c^5s^2+\frac{p_c^5}{2}[10p_c^2-3]s^3+\frac{89}{8}s^4+O(d^{-5}). 
\end{align}
Thus, substituting \refeq{kk} into $p_c$ in the right-hand side of itself yields \refeq{mainthm}. 
\QED

\section{Preliminary estimates for the lace expansion coefficients.}\label{PreliminaryPi0}
This section is devoted to the preparation for estimating $\Pi_{p_c}^{\sss(0)}(\bs{x})$. Throughout the remainder of the paper, we frequently use the following notations. Let $\{\bs{w_i}\}_i$ be a sequence of vertices in $\Lambda$. Using diagrams, we then define
\begin{align}
\mathcal P_p^{\bs{u},\bs{v}}&:=\mathbb P_p\big(\{\bs{o}\rightarrow \bs{u}\}\circ\{\bs{o}\rightarrow \bs{v}\}\big)=\raisebox{-13pt}{\includegraphics[scale=0.85]{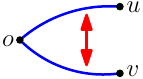}},\lbeq{defPpuv}\\
\mathcal P_p^{\bs{u},\bs{v}}(\{\bs{w_i}\}_i)&:=\mathbb P_p\Big(\big\{\{\bs{o}\rightarrow \bs{u}\}\circ\{\bs{o}\rightarrow \bs{v}\}\big\}\cap \bigcap_{i}\{\bs{o}\rightarrow \bs{w_i}\}\Big)=\raisebox{-30pt}{\includegraphics[scale=0.85]{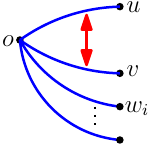}}.\lbeq{defPpuvw} 
\end{align}
In the diagram, blue lines represents paths and the red arrow indicates the mutual avoidance between the two paths. A black disk represent a vertex.  For example, the diagram in \refeq{defPpuv} indicates that there exists at least one pair of disjoint paths from $\bs{o}$ to $\bs{u}$ and $\bs{o}$ to $\bs{v}$. 

We order the support of $D$ in an arbitrary but fixed manner. We say that $x$ is lower than $y$ in that order, which is denoted $x\prec y$, and we also use $x \succ y$ to say that $x$ is higher than $y$. If we write $\bs{x}\prec \bs{y}$ or $\bs{x}\succ \bs{y}$, we mean that $\bs{x}(s)\prec \bs{y}(s)$ or $\bs{x}(s)\succ \bs{y}$, respectively.  
\begin{lmm}\label{doubleconnection}
For $\bs{x}$ and $\bs{y}$ such that $\bs{x}(t)=\bs{y}(t)$, there exist function $\mathcal H_{p}(\bs{x},\bs{y})$ such that for $p\le p_c$, 
\begin{align}
\mathcal P_p^{\bs{x},\bs{y}}=\sum_{\substack{\bs{u},\bs{v}\\ \bs{u}(t)=\bs{v}(t)=\bs{x}(t)-1 \\ \bs{u}\neq \bs{v}}}p^2D(\bs{x}-\bs{u})D(\bs{y}-\bs{v})\mathcal P_p^{\bs{u},\bs{v}}-\mathcal H_{p}(\bs{x},\bs{y}),\lbeq{lmm:oxoy}\\
\mathcal P_p^{\bs{x},\bs{x}}=\sum_{\substack{\bs{u}, \bs{v}\\ \bs{u}(t)=\bs{v}(t)=\bs{x}(t)-1 \\ \bs{x}-\bs{u}\prec \bs{x}-\bs{v} }}p^2D(\bs{x}-\bs{u})D(\bs{x}-\bs{v})\mathcal P_p^{\bs{u},\bs{v}}-\mathcal H_{p}(\bs{x},\bs{x}).\lbeq{lmm:oxox}
\end{align}
Moreover, when $\bs{x}(t)\ge 4$, 
\begin{align}\lbeq{Hxxbdx4}
\mathcal H_{p_c}(\bs{x},\bs{x})\le  \frac{{p_c}^3}{2}(D*\tau_{p_c})(\bs{x})^3+p_c^4\sup_{\bs{x}}(D*\tau_{p_c})(\bs{x})(D*\tau_{p_c})(\bs{x})(D^{*2}*\tau_{p_c}^{*2})(\bs{x}).
\end{align}
\end{lmm}
Concerning the estimates of $\mathcal H_{p_c}(\bs{x},\bs{y})$ in the case where $\bs{x}(t)=\bs{y}(t)\le 3$, we establish Lemma~\ref{lmm: Hxyxt2} and Lemma~\ref{lmm: Hxyxt3}. 

\begin{rmk}[{\bf Sketch of the estimate of $\Pi_{p}^{\sss(0)}$}]~Lemma~\ref{doubleconnection} plays a key role in estimating $\Pi_{p_c}^{\sss(0)}$. Before presenting its proof, we briefly outline its application in this context. 

By using the diagramatic representation, $\Pi_{p_c}^{\sss(0)}(\bs{x})$ can be representad as 
\begin{align}
\Pi_{p}^{\sss(0)}(\bs{x})=\mathcal P_{p}^{\bs{x},\bs{x}}=\raisebox{-11pt}{\includegraphics[scale=0.85]{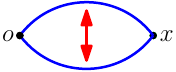}}.
\end{align}
The event that $\bs{o}$ and $\bs{x}$ are doubly connected can be expressed as a union over pairs of points $\bs{u}, \bs{v}$ attached to $\bs{x}$, such that $\bs{x}-\bs{u}\prec \bs{x}-\bs{v}$, the event $\{\bs{o}\rightarrow \bs{u}\}\circ\{\bs{o}\rightarrow \bs{v}\}$ occurs, and the bonds $[\bs{u},\bs{x}\rangle$ and $[\bs{v},\bs{x}\rangle$ are open. Therefore, we obtain 
\begin{align}\lbeq{pi0diag}
\raisebox{-11pt}{\includegraphics[scale=0.85]{Ppxx}}~\sim\sum_{\substack{\bs{u}, \bs{v}\\ \bs{x}-\bs{u}\prec \bs{x}-\bs{v} }}\raisebox{-15pt}{\includegraphics[scale=0.85]{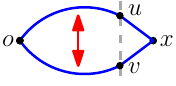}} 
\end{align}
where the right-hand side corresponds to the first term of \refeq{lmm:oxox}, contributing to the leading term of $\Pi_{p}^{\sss(0)}(\bs{x})$. The notation $\sim$ indicates that we neglect  the error, denoted by $\mathcal H_{p}(\bs{x},\bs{x})$ in \refeq{lmm:oxoy}, arising from the fact that the union event is not over disjoint events. 

For evaluating the right-hand side of \refeq{pi0diag}, we need an estimate of $\mathcal P_p^{\bs{u},\bs{v}}$ for $\bs{u}, \bs{v}$ satisfying $\bs{u}(t)=\bs{v}(t)=\bs{x}(t)-1$. The key idea is to estimate $\mathcal P_p^{\bs{x},\bs{y}}$ in order of increasing time components of $\bs{x}$ and $\bs{y}$, starting from the smallest: we first have $\mathcal P_p^{\bs{x},\bs{y}}=p^2D(\bs{x})D(\bs{y})$ for $\bs{x}\neq \bs{y}$ with $\bs{x}(t)=\bs{y}(t)=1$ (see \refeq{ouovDD}), which will be used in estimating $\mathcal P_p^{\bs{x},\bs{x}}$ and  $\mathcal P_p^{\bs{x},\bs{y}}$ for $\bs{x}, \bs{y}$ with $\bs{x}(t)=\bs{y}(t)=2$. Then, using the estimate of $\mathcal P_p^{\bs{x},\bs{y}}$ for $\bs{x}, \bs{y}$ with $\bs{x}(t)=\bs{y}(t)=2$, we estimate $\mathcal P_p^{\bs{x},\bs{y}}$ for $\bs{x}, \bs{y}$ with $\bs{x}(t)=\bs{y}(t)=3$, which in turn is used to estimate $\mathcal P_p^{\bs{x},\bs{x}}$ for $\bs{x}$ with $\bs{x}(t)=4$. 
\end{rmk}

To obtain alternative expression for a double connection, we introduce an order for a set of paths of which starting point and end point are the same in similar manner to that used in \cite{HS05}. 
\begin{itemize}
\item For a pair of paths, $\omega=(b_1,\cdots, b_N)$ and $\omega'=(b_1',\cdots, b_N')$ with $\underline{b}_1=\underline{b}'_1$ and $\bar{b}_N=\bar{b}'_N$, we use $\omega \prec \omega'$ to denote that $\omega$ is lower than $\omega'$, if for the last time $n \in \{1, \cdots, N\}$ such that $b_i=b'_i$ for all $i>n$, $\bar b_n=\bar b'_n$ and $\bar{b}_n -\underline{b}_n\prec \bar{b}'_n-\underline{b}'_n$. We also say that $\omega$ is higher than $\omega'$, which is denoted by $\omega \succ \omega'$.

\item We define an event $E_{\prec}(\omega)$ that $\omega$ is the lowest occupied path among all occupied paths from $\underline{b}_1$ to $\bar{b}_{|\omega|}$ and that there is another occupied path with the length $|\omega|$ which is bond-disjoint from $\omega$ but shares the same starting point. Additionally, given a pair of paths $\omega$ and $\omega'$, we define $E_{\succ}(\omega':\omega)$ to be an event that $\omega'$ is the highest occupied path among all occupied paths that start from $\underline{b}_1$ but are bond-disjoint from $\omega$. 
\end{itemize}

\Proof{Proof of Lemma~\ref{doubleconnection}.}
First, using the above notations, we express the probabilities in \refeq{lmm:oxoy}--\refeq{lmm:oxox} as follows.
\begin{align}\lbeq{oxoy}
\mathbb P_p\big(\{\bs{o}\rightarrow \bs{x}\}\circ\{\bs{o}\rightarrow \bs{y}\}\big)&=\sum_{\substack{\omega_1:\bs{o \rightarrow x}\\\omega_2:\bs{ o\rightarrow y}\\ \omega_1 \cap \omega_2=\varnothing}}\mathbb P_p\big(\{\omega_1,\omega_2~\text{occupied}\} \cap E_{\prec}(\omega_1)\cap E_{\succ}(\omega_2:\omega_1)\big).
\end{align}
Making a condition about the last points that $\omega_1$ and $\omega_2$ visit before they get to their respective end points, we decompose the sum over $\omega_1$ and $\omega_2$ as, if $\bs{x}\neq\bs{y}$,
\begin{align}
\sum_{\substack{\omega_1:\bs{o} \rightarrow \bs{x}\\\omega_2:\bs{o}\rightarrow \bs{y}\\ \omega_1 \cap \omega_2=\varnothing}}=\sum_{\substack{\bs{u},\bs{v}\\ \bs{u}(t)=\bs{x}(t)-1\\ \bs{v}(t)=\bs{y}(t)-1\\ \bs{u} \neq \bs{v}}}\sum_{\substack{\omega_1':\bs{o} \rightarrow \bs{u}\\ \omega'_2:\bs{o}\rightarrow \bs{v}\\ \omega'_1\cap \omega'_2=\varnothing}},
\end{align}
and if $\bs{x}=\bs{y}$, 
\begin{align}\lbeq{sumoxox}
\sum_{\substack{\omega_1,\omega_2:\bs{o}\rightarrow \bs{x}\\ \omega_1 \cap \omega_2=\varnothing}}=\sum_{\substack{\bs{u},\bs{v}\\  \bs{u}(t)=\bs{x}(t)-1\\ \bs{v}(t)=\bs{x}(t)-1\\ \bs{x}-\bs{u}\prec \bs{x}-\bs{v}}}\sum_{\substack{\omega'_1:\bs{o}\rightarrow \bs{u}\\\omega'_2:\bs{o}\rightarrow \bs{v}\\ \omega'_1 \cap \omega'_2=\varnothing}}
\end{align}
Notice that $\omega_1=\omega'_1\cup [\bs{u},\bs{x}\rangle$ and $\omega_2=\omega'_2\cup [\bs{v},\bs{y}\rangle$. Given $\bs{w}, \bs{z}$ such that the bond $[\bs{z},\bs{w}\rangle$ can be defined, $|w-v|=1$, we define $V_{\prec}^{\bs{o}\rightarrow \bs{w}}(\bs{w}-\bs{z})$ to be an event that $\bs{z}$ realizes that $\bs{w}-\bs{z}$ is the lowest among all neighbors of $\bs{w}$ which at least one occupied path from $\bs{o}$ to $\bs{w}$ visit. Furthermore, given $\omega$, we use $V_{\succ}^{\bs{o}\rightarrow \bs{w}}(\bs{w}-\bs{z}: \omega)$ to denote an event that $\bs{z}$ realizes that $\bs{w}-\bs{z}$ is the highest among all neighbors of $\bs{w}$ which at least one occupied path from $\bs{o}$ to $\bs{w}$ that is bond-disjoint from $\omega$ visit. Then, we obtain
\begin{align}
E_{\prec}(\omega'_1\cup[\bs{u},\bs{x}\rangle)&=E_{\prec}(\omega'_1)\cap V_{\prec}^{\bs{o}\rightarrow \bs{x}}(\bs{x}-\bs{u}),\lbeq{Esucc0}\\
E_{\succ}(\omega'_2\cup[\bs{v},\bs{y}\rangle: \omega'_1\cup [\bs{u},\bs{x}\rangle)&=E_{\succ}(\omega'_2: \omega'_1)\cap V_{\succ}^{\bs{o}\rightarrow \bs{y}}(\bs{y}-\bs{v}: \omega'_1).\lbeq{Esucc}
\end{align}
Therefore, defining the events
\begin{align}
Y_{\sss \omega'_1,\omega'_2}&:=\{\omega'_1,\omega'_2~\text{occupied}\}\cap E_{\prec}(\omega'_1)\cap E_{\succ}(\omega'_2:\omega'_1),\lbeq{Ydef}\\
X_{\sss \omega'_1, \omega'_2}^{\sss [\bs{u},\bs{x}\rangle, [\bs{v},\bs{y}\rangle}&:=\{[\bs{u},\bs{x}\rangle, [\bs{v},\bs{y}\rangle~\text{open}\}\cap Y_{\sss \omega'_1,\omega'_2}
\end{align}
the event of \refeq{oxoy} can be expressed as 
\begin{align}\lbeq{XVV}
X_{\sss \omega'_1, \omega'_2}^{\sss [\bs{u},\bs{x}\rangle, [\bs{v},\bs{y}\rangle} \cap V_{\prec}^{\bs{o}\rightarrow \bs{x}}(\bs{u}-\bs{x})\cap V_{\succ}^{\bs{o}\rightarrow \bs{y}}(\bs{v}-\bs{y}: \omega'_1).
\end{align}
From now on, for the sake of simplicity, we write $V_{\prec}^{\bs{o}\rightarrow \bs{x}}(\bs{x}-\bs{u})=V_{\sss[\bs{u}, \bs{x}\rangle}$, $V_{\succ}^{\bs{o}\rightarrow \bs{y}}(\bs{y}-\bs{v}: \omega'_1)=V'_{\sss[\bs{v}, \bs{y}\rangle}$ and $V_{\sss [\bs{u},\bs{x}\rangle,[\bs{v},\bs{y}\rangle}=V_{\sss[\bs{u}, \bs{x}\rangle} \cap V'_{\sss[\bs{v}, \bs{y}\rangle}$. 

If we disregard $V_{\sss [\bs{u},\bs{x}\rangle,[\bs{v},\bs{y}\rangle}$ in \refeq{XVV}, we obtain the main contribution of the probability
\begin{align}\lbeq{oxoymain1}
\sum_{\substack{\omega_1':\bs{o} \rightarrow \bs{u}\\ \omega'_2:\bs{o}\rightarrow \bs{v}\\ \omega'_1\cap \omega'_2=\varnothing}}\mathbb P_p\big(X_{\sss \omega'_1, \omega'_2}^{\sss [\bs{u},\bs{x}\rangle, [\bs{v},\bs{y}\rangle}\big)=p^2D(\bs{x}-\bs{u})D(\bs{y}-\bs{v})\sum_{\substack{\omega_1':\bs{o} \rightarrow \bs{u}\\ \omega'_2:\bs{o}\rightarrow \bs{v}\\ \omega'_1\cap \omega'_2=\varnothing}}\mathbb P_p\big(Y_{\sss \omega'_1,\omega'_2}\big).
\end{align}
Then, the summation in \refeq{oxoymain1} simplifies to $\mathbb P_p\big(\{\bs{o}\rightarrow \bs{u}\}\circ\{\bs{o}\rightarrow \bs{v}\}\big)$, which leads to the main term of \refeq{lmm:oxoy} and \refeq{lmm:oxox}.

Let
\begin{align}
\mathcal H_p(\bs{x},\bs{x})&:=\sum_{\substack{\bs{u},\bs{v}\\ \bs{u}(t)=\bs{v}(t)=\bs{x}(t)-1\\ \bs{x}-\bs{u}\prec \bs{x}-\bs{v}}}\sum_{\substack{\omega_1':\bs{o} \rightarrow \bs{u}\\ \omega'_2:\bs{o}\rightarrow \bs{v}\\ \omega'_1\cap \omega'_2=\varnothing}}\mathbb P_p\big(X_{\sss \omega'_1, \omega'_2}^{\sss [\bs{u},\bs{x}\rangle, [\bs{v},\bs{x}\rangle}\cap V_{\sss [\bs{u},\bs{x}\rangle,[\bs{v},\bs{x}\rangle}^c\big),\\
\mathcal H_p(\bs{x},\bs{y})&:=\sum_{\substack{\bs{u},\bs{v}\\ \bs{u}(t)=\bs{x}(t)-1\\ \bs{v}(t)=\bs{y}(t)-1\\ \bs{u} \neq \bs{v}}}\sum_{\substack{\omega_1':\bs{o} \rightarrow \bs{u}\\ \omega'_2:\bs{o}\rightarrow \bs{v}\\ \omega'_1\cap \omega'_2=\varnothing}}\mathbb P_p\big(X_{\sss \omega'_1, \omega'_2}^{\sss [\bs{u},\bs{x}\rangle, [\bs{v},\bs{y}\rangle}\cap V_{\sss [\bs{u},\bs{x}\rangle,[\bs{v},\bs{y}\rangle}^c\big).
\end{align} 
In the rest of the proof, we aim to show \refeq{Hxxbdx4}. 

{\bf{Estimations of $\mathcal H_{p_c}(\bs{x},\bs{x})$ for the case $\bs{x}(t)\ge 4$.}} ~From now on, we omit the condition on the value of time component in the summation over vertices, assuming that its value is $\bs{x}(t)-1=\bs{y}(t)-1$, unless explicitly required in specific cases. We will bound the following:
\begin{align}\lbeq{erroroxoy0}
&\sum_{\substack{\bs{u},\bs{v}\\ \bs{x}-\bs{u}\prec \bs{x}-\bs{v}}}\sum_{\substack{\omega_1':\bs{o} \rightarrow \bs{u}\\ \omega'_2:\bs{o}\rightarrow \bs{v}\\ \omega'_1\cap \omega'_2=\varnothing}}\mathbb P_p\big(X_{\sss \omega'_1, \omega'_2}^{\sss [\bs{u},\bs{x}\rangle, [\bs{v},\bs{x}\rangle}\cap V_{\sss [\bs{u},\bs{x}\rangle,[\bs{v},\bs{y}\rangle}^{\rm c}\big)\nn\\
&\le \sum_{\substack{\bs{u},\bs{v}\\ \bs{x}-\bs{u}\prec \bs{x}-\bs{v}}}\sum_{\substack{\omega_1':\bs{o} \rightarrow \bs{u}\\ \omega'_2:\bs{o}\rightarrow \bs{v}\\ \omega'_1\cap \omega'_2=\varnothing}}\left\{\mathbb P_p\big(X_{\sss \omega'_1, \omega'_2}^{\sss [\bs{u},\bs{x}\rangle, [\bs{v},\bs{x}\rangle}\cap V_{\sss[\bs{u}, \bs{x}\rangle}^{\rm c}\big)+\mathbb P_p\big(X_{\sss \omega'_1, \omega'_2}^{\sss [\bs{u},\bs{x}\rangle, [\bs{v},\bs{x}\rangle}\cap V_{\sss[\bs{v}, \bs{y}\rangle}^{'\rm c}\big)\right\}.
\end{align}
We only consider the upper bound of the first term of \refeq{erroroxoy0}. That is, we will show 
\begin{align}\lbeq{erroroxoy}
&\sum_{\substack{\bs{u},\bs{v}\\ \bs{x}-\bs{u}\prec \bs{x}-\bs{v}}}\sum_{\substack{\omega_1':\bs{o} \rightarrow \bs{u}\\ \omega'_2:\bs{o}\rightarrow \bs{v}\\ \omega'_1\cap \omega'_2=\varnothing}}\mathbb P_p\big(X_{\sss \omega'_1, \omega'_2}^{\sss [\bs{u},\bs{x}\rangle, [\bs{v},\bs{x}\rangle}\cap V_{\sss[\bs{u}, \bs{x}\rangle}^{\rm c}\big)\nn\\
&\le \frac{p^3}{4}(D*\tau_p)(\bs{x})^3+\frac{p^4}{2}\sup_{\bs{x}}(D*\tau_p)(\bs{x})(D*\tau_p)(\bs{x})(D^{*2}*\tau_p^{*2})(\bs{x})
\end{align}
In a similar way, the second term of \refeq{erroroxoy0} can be bounded above by the right-hand side of \refeq{erroroxoy}.

Since $V_{\sss[\bs{u}, \bs{x}\rangle}^{\rm c}$ is the event that there exists $\bs{w}$ of neighbor of $\bs{x}$ with $\bs{x}-\bs{w}\prec \bs{x}-\bs{u}$ such that there exists at least one occupied paths  from $\bs{o}$ to $[\bs{w},\bs{x}\rangle$, we obtain 
\begin{align}\lbeq{DDDowouov}
\sum_{\substack{\omega_1':\bs{o} \rightarrow \bs{u}\\ \omega'_2:\bs{o}\rightarrow \bs{v}\\ \omega'_1\cap \omega'_2=\varnothing}}\mathbb P_p\big(X_{\sss \omega'_1, \omega'_2}^{\sss [\bs{u},\bs{x}\rangle, [\bs{v},\bs{x}\rangle}\cap V_{\sss[\bs{u}, \bs{x}\rangle}^{\rm c}\big)
&\le \sum_{\substack{\bs{w}\\ \bs{x}-\bs{w}\prec \bs{x}-\bs{u}}}\sum_{\substack{\omega_1':\bs{o} \rightarrow \bs{u}\\ \omega'_2:\bs{o}\rightarrow \bs{v}\\ \omega'_1\cap \omega'_2=\varnothing}}\mathbb P_p\big(X_{\sss \omega'_1, \omega'_2}^{\sss [\bs{u},\bs{x}\rangle, [\bs{v},\bs{x}\rangle}\cap \{\bs{o}\rightarrow [\bs{w},\bs{x}\rangle\}\big)\nn\\
&=p^3D(\bs{x}-\bs{u})D(\bs{x}-\bs{v})\nn\\
&\hskip5mm\times \sum_{\substack{\bs{w}\\ \bs{x}-\bs{w}\prec \bs{x}-\bs{u}}}D(\bs{x}-\bs{w})\sum_{\substack{\omega_1':\bs{o} \rightarrow \bs{u}\\ \omega'_2:\bs{o}\rightarrow \bs{v}\\ \omega'_1\cap \omega'_2=\varnothing}}\mathbb P_p\big(Y_{\sss \omega'_1,\omega'_2}\cap \{\bs{o}\rightarrow \bs{w}\}\big).
\end{align}
Thus, to complete the proof, we consider an upper bound for the expression: 
\begin{align}\lbeq{oucircovow}
\sum_{\substack{\omega_1':\bs{o} \rightarrow \bs{u}\\ \omega'_2:\bs{o}\rightarrow \bs{v}\\ \omega'_1\cap \omega'_2=\varnothing}}\mathbb P_p\big(Y_{\sss \omega'_1,\omega'_2}\cap \{\bs{o}\rightarrow \bs{w}\}\big)=\mathbb P_p\big(\left\{\{\bs{o}\rightarrow \bs{u}\}\circ\{\bs{o}\rightarrow \bs{v}\}\right\}\cap \{\bs{o}\rightarrow \bs{w}\}\big).
\end{align}
Here, since $\bs{x}-\bs{u} \prec \bs{x}-\bs{v}$ and $\bs{w}$ is required to satisfy $\bs{x}-\bs{w}\prec\bs{x}-\bs{u}$, it follows that $\bs{w}$ cannot be $\bs{v}$ when $\bs{x}=\bs{y}$. 

Applying the BK inequality allows us to bound the probability in \refeq{oucircovow} by a sum involving $\tau_p(\bs{u})\tau_p(\bs{v})\tau_p(\bs{w})$ and an additional term
\begin{align}\lbeq{owouovossw}
\mathbb P_p\big(\left\{\{\bs{o}\rightarrow \bs{u}\}\circ\{\bs{o}\rightarrow \bs{v}\}\right\}\cap \{\bs{o}\rightarrow \bs{w}\} \backslash\{\bs{o}\rightarrow \bs{u}\}\circ\{\bs{o}\rightarrow \bs{v}\}\circ\{\bs{o}\rightarrow \bs{w}\}\big).
\end{align}
The contribution from $\tau_p(\bs{w})\tau_p(\bs{u})\tau_p(\bs{v})$ to \refeq{erroroxoy} is bounded above by
\begin{align}\lbeq{wuvcont}
p^3\sum_{\substack{\bs{u},\bs{v}\\ \bs{x}-\bs{u}\prec \bs{x}-\bs{v}}}D(\bs{x}-\bs{u})D(\bs{x}-\bs{v})\tau_p(\bs{u})\tau_p(\bs{v})\sum_{\substack{\bs{w}\\ \bs{x}-\bs{w}\prec \bs{x}-\bs{u}}}D(\bs{x}-\bs{w})\tau_p(\bs{w})
&\le \frac{p^3}{4}(D*\tau_p)(\bs{x})^3.
\end{align}

Now we consider the estimete of \refeq{owouovossw}. We note that, with $\bs{s}$ taken as the last intersection point of the path from $\bs{o}$ to $\bs{w}$ with the path of $\bs{o}\rightarrow \bs{u}$ or $\bs{o}\rightarrow \bs{v}$, the event in \refeq{owouovossw} is contained in the following event
\begin{align}
\bigcup_{\bs{t}=\bs{u},\bs{v}}\bigcup_{\bs{s}\neq \bs{o},\bs{t}}\left\{\{\bs{o}\rightarrow \bs{s} \rightarrow \bs{t}\}\circ \{\bs{o}\rightarrow \bs{\dot{t}}\}\circ \{\bs{s}\rightarrow \bs{w}\}\right\}\cap \{\bs{o}\rightarrow \bs{s}\},
\end{align}
where $\bs{\dot{t}}$ is used to represent point other than $\bs{t}$, i.e., $\bs{\dot{u}}=\bs{v}$. Hence, disregarding the event $\bs{o}\rightarrow \bs{s}$ and applying the BK inequality again along with the Markov property, \refeq{owouovossw} can be bounded above by
\begin{align}
\sum_{\bs{t}=\bs{u},\bs{v}}\sum_{\bs{s}\neq\bs{o},\bs{t}}\tau_p(\bs{s})\tau_p(\bs{t}-\bs{s})\tau_p(\bs{\dot{t}})\tau_p(\bs{w}-\bs{s}). 
\end{align}
Therefore, the contribution from \refeq{owouovossw} to \refeq{erroroxoy} is bounded above by
\begin{align}\lbeq{owouovos}
&2p^3\sum_{\substack{\bs{u},\bs{v}\\ \bs{x}-\bs{u}\prec \bs{x}-\bs{v}}}D(\bs{x}-\bs{u})D(\bs{x}-\bs{v})\sum_{\substack{\bs{w}\\ \bs{x}-\bs{w}\prec \bs{x}-\bs{u}}}D(\bs{x}-\bs{w})\sum_{\bs{s}\neq\bs{o},\bs{u}}\tau_p(\bs{s})\tau_p(\bs{u}-\bs{s})\tau_p(\bs{v})\tau_p(\bs{w}-\bs{s})\nn\\
&\le \frac{p^4}{2}(D*\tau_p)(\bs{x})\sum_{\bs{s}\neq\bs{o}}(D*\tau_p)(\bs{s})(D*\tau_p)(\bs{x}-\bs{s})^2\nn\\
&\le \frac{p^4}{2}\sup_{\bs{x}}(D*\tau_p)(\bs{x})(D*\tau_p)(\bs{x})(D^{*2}*\tau_p^{*2})(\bs{x})
\end{align}
which together with\refeq{wuvcont}, provides \refeq{erroroxoy}. \QED

\subsection{Estimations of $\mathcal H_{p_c}(\bs{x},\bs{x})$ and $\mathcal H_{p_c}(\bs{x},\bs{y})$ for the case $\bs{x}(t)=2,3$.}\label{time2and3}
This section provides preliminaries for estimating $\mathcal H_{p_c}(\bs{x},\bs{x})$ and $\mathcal H_{p_c}(\bs{x},\bs{y})$ in the cases where $\bs{x}(t)=2$ and $\bs{x}(t)=3$. 

Since $V_{\sss [\bs{u},\bs{x}\rangle,[\bs{v},\bs{y}\rangle}:=V_{\sss[\bs{u}, \bs{x}\rangle}\cap V_{\sss[\bs{v}, \bs{y}\rangle}^{'}$, the probability under consideration is exactly decomposed as
\begin{align}\lbeq{XVc}
\mathbb P_p\big(X\cap V_{\sss [\bs{u},\bs{x}\rangle,[\bs{v},\bs{y}\rangle}^{c}\big)=\mathbb P_p\big(X\cap V_{\sss [\bs{u}, \bs{x}\rangle}^{c}\big)+\mathbb P_p\big(X\cap V_{\sss[\bs{v}, \bs{y}\rangle}^{'c}\big)-\mathbb P_p\big(X\cap V_{\sss[\bs{u}, \bs{x}\rangle}^{c}\cap V_{\sss[\bs{v}, \bs{y}\rangle}^{'c}\big).
\end{align}
Here, $X$ is the abbreviation for $X_{\sss \omega'_1,\omega'_2}^{\sss [\bs{u},\bs{x}\rangle, [\bs{v},\bs{y}\rangle}$ and we sometimes use this abbreviation throughout the rest of the section. 

Notice that $V_{\sss[\bs{u}, \bs{x}\rangle}^{c}$ and $V_{\sss[\bs{v}, \bs{y}\rangle}^{'c}$ can be rewritten as 
\begin{align}
V_{\sss[\bs{u}, \bs{x}\rangle}^{c}&=\bigcup_{\substack{\bs{w}\\ \bs{w}(t)=\bs{x}(t)-1\\ \bs{x}-\bs{w}\prec\bs{x}-\bs{u}}}\{\bs{o}\rightarrow [\bs{w}, \bs{x}\rangle\}\cap V_{\sss[\bs{w}, \bs{x}\rangle},\lbeq{Voxc}\\
V_{\sss[\bs{v}, \bs{y}\rangle}^{'c}&=\bigcup_{\substack{\bs{w}\\ \bs{w}(t)=\bs{y}(t)-1 \\ \bs{y}-\bs{w}\succ\bs{y}-\bs{v}}}\{\bs{o}\rightarrow[\bs{w}, \bs{y}\rangle\}\cap V_{\sss[\bs{w}, \bs{y}\rangle}^{'},\lbeq{Voyc}
\end{align}
and the events on the right-hand side of \refeq{Voxc}--\refeq{Voyc} are disjoint events. Therefore, we have 
\begin{align}
\mathbb P_p\big(X\cap V_{\sss[\bs{u}, \bs{x}\rangle}^{c}\big)&=\sum_{\substack{\bs{w}\\  \bs{x}-\bs{w}\prec\bs{x}-\bs{u}}}\mathbb P_p\big(X\cap \{\bs{o}\rightarrow [\bs{w}, \bs{x}\rangle\} \cap V_{\sss[\bs{w}, \bs{x}\rangle}\big), \lbeq{Ppxvxc0}\\
\mathbb P_p\big(X\cap V_{\sss[\bs{v}, \bs{y}\rangle}^{'c}\big)&=\sum_{\substack{\bs{w}\\ \bs{y}-\bs{w}\succ\bs{y}-\bs{v}}}\mathbb P_p\big(X\cap \{\bs{o}\rightarrow [\bs{w}, \bs{y}\rangle\} \cap V'_{\sss[\bs{w}, \bs{y}\rangle}\big)\lbeq{Ppxvxc1}
\end{align}
and 
\begin{align}\lbeq{Ppxvxc2}
&\mathbb P_p\big(X\cap V_{\sss[\bs{u}, \bs{x}\rangle}^{c}\cap V_{\sss[\bs{v}, \bs{y}\rangle}^{'c}\big)\nn\\
&=\sum_{\substack{\bs{w},\bs{z} \\ \bs{x}-\bs{w}\prec\bs{x}-\bs{u} \\ \bs{y}-\bs{z}\succ\bs{y}-\bs{v}}}\mathbb P_p\big(X \cap\{\bs{o}\rightarrow [\bs{w}, \bs{x}\rangle\}\cap\{\bs{o}\rightarrow [\bs{z}, \bs{y}\rangle\} \cap V_{\sss[\bs{w}, \bs{x}\rangle}\cap V'_{\sss[\bs{z}, \bs{y}\rangle} \big).
\end{align}
To estimate $\mathcal H_p(\bs{x},\bs{x})$ and $\mathcal H_p(\bs{x},\bs{y})$ in the cases where $\bs{x}(t)=2$ and $\bs{x}(t)=3$, we expand \refeq{Ppxvxc0}--\refeq{Ppxvxc2} by iteratively using the relation of the events given in \refeq{Voxc}--\refeq{Voyc}. For example, the probability of the summand in \refeq{Ppxvxc0} can be rewritten as
\begin{align}\lbeq{Ppxvxc}
\mathbb P_p\big(X\cap  \{\bs{o}\rightarrow [\bs{w}, \bs{x}\rangle\}\cap V_{\sss[\bs{w}, \bs{x}\rangle}\big)&=\mathbb P_p\big(X\cap  \{\bs{o}\rightarrow [\bs{w}, \bs{x}\rangle\}\big)\nn\\
&\hskip5mm-\mathbb P_p\big(X\cap  \{\bs{o}\rightarrow [\bs{w}, \bs{x}\rangle\}\cap V_{\sss[\bs{w}, \bs{x}\rangle}^{c}\big)
\end{align}
and we apply the relation given in \refeq{Voxc}, replacing $[\bs{u}, \bs{x}\rangle$ with $[\bs{w}, \bs{x}\rangle$, to the second term of \refeq{Ppxvxc}. Similarly, we will follow the same approach for \refeq{Ppxvxc1}--\refeq{Ppxvxc2}. The number of times we need to iterate this procedure to obtain the explicit expression up to the order of $d^{-4}$ differs between the cases $\bs{x}(t)=2$ and $\bs{x}(t)=3$. For the case $\bs{x}(t)=2$, we stop the iteration after 2 steps, while for the case $\bs{x}(t)=3$, we stop it after 1 step. We use $\mathcal E^{(n)}_p$, $\mathcal E'^{(n)}_p$ and $\mathcal E''^{(n)}_p$, which depend on $[\bs{u},\bs{x}\rangle$ and $[\bs{v},\bs{y}\rangle$, to represent the contributions to $\mathcal H(\bs{x},\bs{x})$ or $\mathcal H(\bs{x},\bs{y})$ from the terms that arise by neglecting the events $V$, $V'$ and $V\cap V'$, respectively after applying the relation in \refeq{Voxc}--\refeq{Voyc} $n$ times. Furthermore, $E^{(n+1)}_p$, $E'^{(n+1)}_p$ and $E''^{(n+1)}_p$ represent the remaining contributions coming from $V^c$, $V'^c$ and $V^c\cap V'^c$ respectively after extracting $\mathcal E^{(n)}_p$, $\mathcal E'^{(n)}_p$ and $\mathcal E''^{(n)}_p$. Specifically, by applying \refeq{Voxc} to $E^{(n+1)}_p$, we obtain the relation $E^{(n+1)}_p=\mathcal E^{(n+1)}_p-E^{(n+2)}_p$. A similar relation holds for $E'^{(n+1)}_p$ and $E''^{(n+1)}_p$. 
For instance, since the first term of \refeq{Ppxvxc} results from the first application of \refeq{Voxc} to the right-hand side of \refeq{Ppxvxc0}, $\mathcal E^{(1)}_p=\mathcal E^{(1)}_p(\substack{\sss [\bs{u},\bs{x}\rangle \\ \sss [\bs{v},\bs{y}\rangle})$ and $E^{(2)}_p=E^{(2)}_p(\substack{\sss [\bs{u},\bs{x}\rangle \\ \sss [\bs{v},\bs{y}\rangle})$ are defined as
\begin{align}
\mathcal E^{(1)}_p(\substack{\sss [\bs{u},\bs{x}\rangle \\ \sss [\bs{v},\bs{y}\rangle})&=\sum_{\substack{\omega_1':\bs{o} \rightarrow \bs{u}\\ \omega'_2:\bs{o}\rightarrow \bs{v}\\ \omega'_1\cap \omega'_2=\varnothing}}\sum_{\substack{\bs{w}\\ \bs{x}-\bs{w}\prec\bs{x}-\bs{u}}}\mathbb P_p\big(X_{\sss \omega'_1,\omega'_2}^{\sss [\bs{u},\bs{x}\rangle, [\bs{v},\bs{y}\rangle}\cap  \{\bs{o}\rightarrow [\bs{w}, \bs{x}\rangle\}\big),\\
E^{(2)}_p(\substack{\sss [\bs{u},\bs{x}\rangle \\ \sss [\bs{v},\bs{y}\rangle})&=\sum_{\substack{\omega_1':\bs{o} \rightarrow \bs{u}\\ \omega'_2:\bs{o}\rightarrow \bs{v}\\ \omega'_1\cap \omega'_2=\varnothing}}\sum_{\substack{\bs{w}\\ \bs{x}-\bs{w}\prec\bs{x}-\bs{u}}}\mathbb P_p\big(X_{\sss \omega'_1,\omega'_2}^{\sss [\bs{u},\bs{x}\rangle, [\bs{v},\bs{y}\rangle}\cap  \{\bs{o}\rightarrow [\bs{w}, \bs{x}\rangle\}\cap V_{\sss[\bs{w}, \bs{x}\rangle}^{c}\big). \lbeq{exE2}
\end{align} 
Furthermore, by applying \refeq{Voxc} to $V_{\sss[\bs{w}, \bs{x}\rangle}^{c}$ in \refeq{exE2}, we obtain $E^{(2)}_p=\mathcal E^{(2)}_p-E^{(3)}_p$. 
It is easy to verify that by applying the relation in \refeq{Voxc}--\refeq{Voyc} to $E^{(n)}_p$ for $n\ge2$, we obtain $E^{(n)}_p\le \mathcal E^{(n)}_p$. Similarly, it can be verified that $E'^{(n)}_p\le \mathcal E'^{(n)}_p$ and $E''^{(n)}_p\le \mathcal E''^{(n)}_p$. Then, applying the procedure twice for the case $\bs{x}(t)=2$ and once for the case $\bs{x}(t)=3$ and using \refeq{XVc}, we obtain the followings: In the case where $\bs{x}(t)=2$, for both cases $\bs{x}=\bs{y}$ and $\bs{x}\neq\bs{y}$, 
\begin{align}\lbeq{Hxxsummandt2}
&\sum_{\substack{\omega_1':\bs{o} \rightarrow \bs{u}\\ \omega'_2:\bs{o}\rightarrow \bs{v}\\ \omega'_1\cap \omega'_2=\varnothing}}\mathbb P_p\big(X_{\sss \omega'_1,\omega'_2}^{\sss [\bs{u},\bs{x}\rangle, [\bs{v},\bs{y}\rangle}\cap V_{\sss [\bs{u},\bs{x}\rangle,[\bs{v},\bs{y}\rangle}^{c}\big)\nn\\
&=\left[\mathcal E^{(1)}_p(\substack{\sss [\bs{u},\bs{x}\rangle \\ \sss [\bs{v},\bs{y}\rangle})+\mathcal E'^{(1)}_p(\substack{\sss [\bs{u},\bs{x}\rangle \\ \sss [\bs{v},\bs{y}\rangle})\right]-\left[\mathcal E^{(2)}_p(\substack{\sss [\bs{u},\bs{x}\rangle \\ \sss [\bs{v},\bs{y}\rangle})+\mathcal E'^{(2)}_p(\substack{\sss [\bs{u},\bs{x}\rangle \\ \sss [\bs{v},\bs{y}\rangle})\right]-\mathcal E''^{(1)}_p(\substack{\sss [\bs{u},\bs{x}\rangle \\ \sss [\bs{v},\bs{y}\rangle})\nn\\
&+\left[E^{(3)}_p(\substack{\sss [\bs{u},\bs{x}\rangle \\ \sss [\bs{v},\bs{y}\rangle})+E'^{(3)}_p(\substack{\sss [\bs{u},\bs{x}\rangle \\ \sss [\bs{v},\bs{y}\rangle})\right]+E''^{(2)}_p(\substack{\sss [\bs{u},\bs{x}\rangle \\ \sss [\bs{v},\bs{y}\rangle}).
\end{align}
In the case where $\bs{x}(t)=3$, for both cases $\bs{x}=\bs{y}$ and $\bs{x}\neq\bs{y}$, 
\begin{align}\lbeq{Hxxsummandt3}
&\sum_{\substack{\omega_1':\bs{o} \rightarrow \bs{u}\\ \omega'_2:\bs{o}\rightarrow \bs{v}\\ \omega'_1\cap \omega'_2=\varnothing}}\mathbb P_p\big(X_{\sss \omega'_1,\omega'_2}^{\sss [\bs{u},\bs{x}\rangle, [\bs{v},\bs{y}\rangle}\cap V_{\sss [\bs{u},\bs{x}\rangle,[\bs{v},\bs{y}\rangle}^{c}\big)\nn\\
&=\left[\mathcal E^{(1)}_p(\substack{\sss [\bs{u},\bs{x}\rangle \\ \sss [\bs{v},\bs{y}\rangle})+\mathcal E'^{(1)}_p(\substack{\sss [\bs{u},\bs{x}\rangle \\ \sss [\bs{v},\bs{y}\rangle})\right]-\left[E^{(2)}_p(\substack{\sss [\bs{u},\bs{x}\rangle \\ \sss [\bs{v},\bs{y}\rangle})+E'^{(2)}_p(\substack{\sss [\bs{u},\bs{x}\rangle \\ \sss [\bs{v},\bs{y}\rangle})\right]-E''^{(1)}_p(\substack{\sss [\bs{u},\bs{x}\rangle \\ \sss [\bs{v},\bs{y}\rangle}). 
\end{align}
To obtain the explicit expressions of $\mathcal H_{p_c}(\bs{x},\bs{x})$ and $\mathcal H_{p_c}(\bs{x},\bs{y})$ up to the terms that contribute to the expression of \refeq{mainthm} up to order $d^{-4}$, we will consider the explicit expressions of $\mathcal E^{(n)}_p$, $\mathcal E'^{(n)}_p$ and $\mathcal E''^{(n)}_p$ for $n=1$ or $2$ and derive upper bounds of $E^{(n)}_p$, $E'^{(n)}_p$ and $E''^{(n)}_p$ for $n=1,2$ or $3$. In both cases, where $\bs{x}(t)=2$ and $\bs{x}(t)=3$, we need to estimate $\big[\mathcal E^{(1)}_p+\mathcal E'^{(1)}_p\big]$, $\big[\mathcal E^{(2)}_p+\mathcal E'^{(2)}_p\big]$ and $\mathcal E''^{(1)}_p$. However, for the case where $\bs{x}(t)=3$, we estimate $\big[\mathcal E^{(2)}_p+\mathcal E'^{(2)}_p\big]$ and $\mathcal E''^{(1)}_p$ as upper bounds for $\big[E^{(2)}_p+E'^{(2)}_p\big]$ and $E''^{(1)}_p$ to demonstrate that they contribute to the error terms. As for $\big[E^{(3)}_p+E'^{(3)}_p\big]$ and $E''^{(2)}_p$, we need to estimate them as error terms only for the case $\bs{x}(t)=2$. 

In the remainder of Section~\ref{time2and3}, we rewrite $\big[\mathcal E^{(i)}_p+\mathcal E'^{(i)}_p\big]$ for $i=1,2$, $E''^{(1)}_p$, $\big[E^{(3)}_p+E'^{(3)}_p\big]$ and $E''^{(2)}_p$ in a form that facilitates their estimation. 

\subsubsection{Estimation of $\mathcal E^{(1)}_p$ and $\mathcal E'^{(1)}_p$}
First, we recall the definitions of $\mathcal E^{(1)}_p(\substack{\sss [\bs{u},\bs{x}\rangle \\ \sss [\bs{v},\bs{y}\rangle})$ and $\mathcal E'^{(1)}_p(\substack{\sss [\bs{u},\bs{x}\rangle \\ \sss [\bs{v},\bs{y}\rangle})$: 
\begin{align}
\mathcal E^{(1)}_p(\substack{\sss [\bs{u},\bs{x}\rangle \\ \sss [\bs{v},\bs{y}\rangle})&=\sum_{\substack{\omega_1':\bs{o} \rightarrow \bs{u}\\ \omega'_2:\bs{o}\rightarrow \bs{v}\\ \omega'_1\cap \omega'_2=\varnothing}}\sum_{\substack{\bs{w}\\ \bs{x}-\bs{w}\prec\bs{x}-\bs{u}}}\mathbb P_p\big(X_{\sss \omega'_1, \omega'_2}^{\sss [\bs{u},\bs{x}\rangle, [\bs{v},\bs{y}\rangle} \cap  \{\bs{o}\rightarrow [\bs{w}, \bs{x}\rangle\}\big)\nn\\
&=p^3D(\bs{x}-\bs{u})D(\bs{y}-\bs{v})\sum_{\substack{\bs{w}\\ \bs{x}-\bs{w}\prec\bs{x}-\bs{u}}}D(\bs{x}-\bs{w})\mathcal P_p^{\bs{u},\bs{v}}(\bs{w}), \lbeq{omegaXVcux}\\
\mathcal E'^{(1)}_p(\substack{\sss [\bs{u},\bs{x}\rangle \\ \sss [\bs{v},\bs{y}\rangle})&=\sum_{\substack{\omega_1':\bs{o} \rightarrow \bs{u}\\ \omega'_2:\bs{o}\rightarrow \bs{v}\\ \omega'_1\cap \omega'_2=\varnothing}}\sum_{\substack{\bs{w}\\ \bs{y}-\bs{w}\succ\bs{y}-\bs{v}}}\mathbb P_p\big(X_{\sss \omega'_1, \omega'_2}^{\sss [\bs{u},\bs{x}\rangle, [\bs{v},\bs{y}\rangle} \cap \{\bs{o}\rightarrow [\bs{w}, \bs{y}\rangle\}\big)\nn\\
&=p^3D(\bs{x}-\bs{u})D(\bs{y}-\bs{v})\sum_{\substack{\bs{w}\\ \bs{y}-\bs{w}\succ\bs{y}-\bs{v}}}D(\bs{y}-\bs{w})\mathcal P_p^{\bs{u},\bs{v}}(\bs{w}).\lbeq{omegaXVV}
\end{align}
Here, in the second equalities of both \refeq{omegaXVcux} and \refeq{omegaXVV}, we use the same argument used to obtain \refeq{erroroxoy}. 
In this section, we simplify

\begin{align}
\sum_{\substack{\bs{u},\bs{v}\\ \bs{x}-\bs{u}\prec \bs{x}-\bs{v}}}\Big[\mathcal E^{(1)}_{p_c}(\substack{\sss [\bs{u},\bs{x}\rangle \\ \sss [\bs{v},\bs{x}\rangle})+\mathcal E'^{(1)}_{p_c}(\substack{\sss [\bs{u},\bs{x}\rangle \\ \sss [\bs{v},\bs{x}\rangle})\Big]&&\text{and}&&\sum_{\substack{\bs{u},\bs{v}\\ \bs{u} \neq \bs{v}}}\Big[\mathcal E^{(1)}_{p_c}(\substack{\sss [\bs{u},\bs{x}\rangle \\ \sss [\bs{v},\bs{y}\rangle})+\mathcal E'^{(1)}_{p_c}(\substack{\sss [\bs{u},\bs{x}\rangle \\ \sss [\bs{v},\bs{y}\rangle})\Big].
\end{align}

We now provide comments on the sum over $\bs{w}$ in \refeq{omegaXVcux}--\refeq{omegaXVV}. For \refeq{omegaXVcux}, since $\bs{x}-\bs{u} \prec \bs{x}-\bs{v}$ when $\bs{x}=\bs{y}$ and $\bs{w}$ is required to satisfy $\bs{x}-\bs{w}\prec\bs{x}-\bs{u}$, it follows that $\bs{w}$ cannot be $\bs{v}$ when $\bs{x}=\bs{y}$. In contrast, for the case $\bs{x}\neq\bs{y}$, the only restriction is $\bs{u}\neq\bs{v}$, which allows $\bs{w}$ to be $\bs{v}$. A similar observation applies to the sum over $\bs{w}$ in \refeq{omegaXVV}. Sepecifically, if $\bs{x}=\bs{y}$, $\bs{z}$ cannot be $\bs{u}$, whereas if $\bs{x}\neq\bs{y}$, $\bs{z}$ can be $\bs{u}$. Thus, observing that when $\bs{w}=\bs{v}$ or $\bs{w}=\bs{u}$
\begin{align}
\mathcal P_p^{\bs{u},\bs{v}}(\bs{w})=\mathcal P_p^{\bs{u},\bs{v}},
\end{align}
we obtain
\begin{align}
\sum_{\substack{\bs{u},\bs{v}\\ \bs{x}-\bs{u}\prec \bs{x}-\bs{v}}}\Big[\mathcal E^{(1)}_{p_c}(\substack{\sss [\bs{u},\bs{x}\rangle \\ \sss [\bs{v},\bs{x}\rangle})+\mathcal E'^{(1)}_{p_c}(\substack{\sss [\bs{u},\bs{x}\rangle \\ \sss [\bs{v},\bs{x}\rangle})\Big]&=\frac{1}{4}F_{2,p_c}(\bs{x},\bs{x}),\lbeq{F2xx}\\
\sum_{\substack{\bs{u},\bs{v}\\ \bs{u} \neq \bs{v}}}\Big[\mathcal E^{(1)}_{p_c}(\substack{\sss [\bs{u},\bs{x}\rangle \\ \sss [\bs{v},\bs{y}\rangle})+\mathcal E'^{(1)}_{p_c}(\substack{\sss [\bs{u},\bs{x}\rangle \\ \sss [\bs{v},\bs{y}\rangle})\Big]&=\frac{1}{2}F_{1,p_c}(\bs{x},\bs{y})+\frac{1}{2}F_{2,p_c}(\bs{x},\bs{y})\lbeq{F1xyF2xy}
\end{align}
where
\begin{align}
F_{1,p}(\bs{x},\bs{y})&=p^3\sum_{\substack{\bs{u},\bs{v}\\ \bullet}}D(\bs{x}-\bs{u})D(\bs{y}-\bs{v})\big[D(\bs{x}-\bs{v})+D(\bs{y}-\bs{u})\big]\mathcal P_p^{\bs{u},\bs{v}},\lbeq{F1def}\\
F_{2,p}(\bs{x},\bs{y})
&=p^3\sum_{\substack{\bs{u},\bs{v},\bs{w}\\ \bullet}}D(\bs{x}-\bs{u})D(\bs{y}-\bs{v})\big[D(\bs{x}-\bs{w})+D(\bs{y}-\bs{w})\big]\mathcal P_p^{\bs{u},\bs{v}}(\bs{w}).\lbeq{F2def}
\end{align}

Here, the symbol $\bullet$ indicates that the summation is taken over mutually distinct points. For example, 
\begin{align}
\sum_{\substack{\bs{u},\bs{v}\\ \bs{u} \neq \bs{v}}}=\sum_{\substack{\bs{u},\bs{v}\\ \bullet}},&& \sum_{\substack{\bs{u},\bs{v},\bs{w}\\ \bs{u} \neq \bs{v}\\\bs{w}\neq \bs{u}, \bs{v}}}=\sum_{\substack{\bs{u},\bs{v},\bs{w}\\ \bullet}}. 
\end{align}

By using the diagrammatic representation defined in \refeq{defPpuv}--\refeq{defPpuvw}, we rewrite \refeq{F1def} and \refeq{F2def} as
\begin{align}
F_{1,p}(\bs{x},\bs{y})&=\raisebox{-15pt}{\includegraphics[scale=0.85]{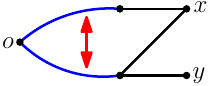}}+\raisebox{-15pt}{\includegraphics[scale=0.85]{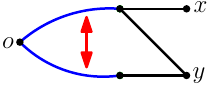}},\\
F_{2,p}(\bs{x},\bs{y})&=\raisebox{-15pt}{\includegraphics[scale=0.85]{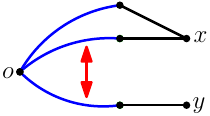}}+\raisebox{-15pt}{\includegraphics[scale=0.85]{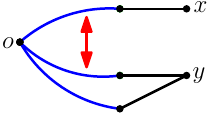}}. 
\end{align}
Here, an unlabelled black vertex is summed over all possible vertices, ensuring that all black vertices are distinct. A black line connecting two disks represents $pD$.

\subsubsection{Estimation of $\mathcal E^{(2)}_p$ and $\mathcal E'^{(2)}_p$}
Recall the definitions of $\mathcal E^{(2)}_p(\substack{\sss [\bs{u},\bs{x}\rangle \\ \sss [\bs{v},\bs{y}\rangle})$ and $\mathcal E'^{(2)}_p(\substack{\sss [\bs{u},\bs{x}\rangle \\ \sss [\bs{v},\bs{y}\rangle})$: 
\begin{align}
\mathcal E^{(2)}_p(\substack{\sss [\bs{u},\bs{x}\rangle \\ \sss [\bs{v},\bs{y}\rangle})&=\sum_{\substack{\omega_1':\bs{o} \rightarrow \bs{u}\\ \omega'_2:\bs{o}\rightarrow \bs{v}\\ \omega'_1\cap \omega'_2=\varnothing}}\sum_{\substack{\bs{w}, \bs{z}\\ \bs{x}-\bs{w}\prec\bs{x}-\bs{u} \\ \bs{x}-\bs{z}\prec\bs{x}-\bs{w}}}\mathbb P_p\big(X\cap \{\bs{o}\rightarrow [\bs{w}, \bs{x}\rangle\}\cap \{\bs{o}\rightarrow [\bs{z}, \bs{x}\rangle\}\big)\nn\\
&=p^4D(\bs{x}-\bs{u})D(\bs{y}-\bs{v})\sum_{\substack{\bs{w}, \bs{z}\\ \bs{x}-\bs{w}\prec\bs{x}-\bs{u} \\ \bs{x}-\bs{z}\prec\bs{x}-\bs{w}}}D(\bs{x}-\bs{w})D(\bs{x}-\bs{z})\mathcal P_p^{\bs{u},\bs{v}}(\bs{w},\bs{z}).\\
\mathcal E'^{(2)}_p(\substack{\sss [\bs{u},\bs{x}\rangle \\ \sss [\bs{v},\bs{y}\rangle})&=\sum_{\substack{\omega_1':\bs{o} \rightarrow \bs{u}\\ \omega'_2:\bs{o}\rightarrow \bs{v}\\ \omega'_1\cap \omega'_2=\varnothing}}\sum_{\substack{\bs{w}, \bs{z}\\ \bs{y}-\bs{w}\succ\bs{y}-\bs{v} \\ \bs{y}-\bs{z}\succ\bs{y}-\bs{w}}}\mathbb P_p\big(X\cap \{\bs{o}\rightarrow [\bs{w}, \bs{y}\rangle\}\cap \{\bs{o}\rightarrow [\bs{z}, \bs{y}\rangle\}\big)\nn\\
&=p^4D(\bs{x}-\bs{u})D(\bs{y}-\bs{v})\sum_{\substack{\bs{w}, \bs{z}\\ \bs{y}-\bs{w}\succ\bs{y}-\bs{u} \\ \bs{y}-\bs{z}\succ\bs{y}-\bs{w}}}D(\bs{y}-\bs{w})D(\bs{y}-\bs{z})\mathcal P_p^{\bs{u},\bs{v}}(\bs{w},\bs{z}).
\end{align}
Then, by considering the case where $\bs{w}$ or $\bs{z}$ equals to $\bs{v}$ for $\mathcal E^{(2)}_p(\substack{\sss [\bs{u},\bs{x}\rangle \\ \sss [\bs{v},\bs{y}\rangle})$, and the case where $\bs{w}$ or $\bs{z}$ equals to $\bs{u}$ for $\mathcal E'^{(2)}_p(\substack{\sss [\bs{u},\bs{x}\rangle \\ \sss [\bs{v},\bs{y}\rangle})$ when $\bs{x}\neq \bs{y}$, we obtain
\begin{align}
\sum_{\substack{\bs{u},\bs{v}\\ \bs{x}-\bs{u}\prec \bs{x}-\bs{v}}}\Big[\mathcal E^{(2)}_{p_c}(\substack{\sss [\bs{u},\bs{x}\rangle \\ \sss [\bs{v},\bs{x}\rangle})+\mathcal E'^{(2)}_{p_c}(\substack{\sss [\bs{u},\bs{x}\rangle \\ \sss [\bs{v},\bs{x}\rangle})\Big]&=\frac{1}{8}F_{4,p}(\bs{x},\bs{x}),\lbeq{F4xx}\\
\sum_{\substack{\bs{u},\bs{v}\\ \bs{u} \neq \bs{v}}}\Big[\mathcal E^{(2)}_{p_c}(\substack{\sss [\bs{u},\bs{x}\rangle \\ \sss [\bs{v},\bs{y}\rangle})+\mathcal E'^{(2)}_{p_c}(\substack{\sss [\bs{u},\bs{x}\rangle \\ \sss [\bs{v},\bs{y}\rangle})\Big]&=\frac{1}{2}F_{3,p}(\bs{x},\bs{y})+\frac{1}{4}F_{4,p}(\bs{x},\bs{y})\lbeq{F3xyF4xy}
\end{align}
where
\begin{align}
F_{3,p}(\bs{x},\bs{y})&=\raisebox{-15pt}{\includegraphics[scale=0.85]{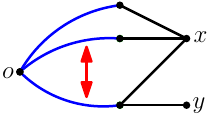}}+\raisebox{-15pt}{\includegraphics[scale=0.85]{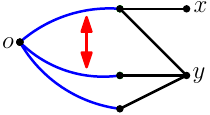}},\nn\\
F_{4,p}(\bs{x},\bs{y})&=\raisebox{-15pt}{\includegraphics[scale=0.85]{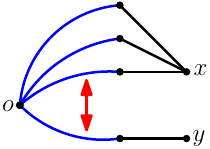}}+\raisebox{-15pt}{\includegraphics[scale=0.85]{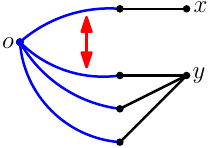}}. \lbeq{F4def}
\end{align}

\subsubsection{Estimation of $\mathcal E''^{(1)}_p$}
Recall the definitions of $\mathcal E''^{(1)}_p(\substack{\sss [\bs{u},\bs{x}\rangle \\ \sss [\bs{v},\bs{y}\rangle})$: 
\begin{align}
\mathcal E''^{(1)}_p(\substack{\sss [\bs{u},\bs{x}\rangle \\ \sss [\bs{v},\bs{y}\rangle})&=\sum_{\substack{\omega_1':\bs{o} \rightarrow \bs{u}\\ \omega'_2:\bs{o}\rightarrow \bs{v}\\ \omega'_1\cap \omega'_2=\varnothing}}\sum_{\substack{\bs{w},\bs{z}\\ \bs{x}-\bs{w}\prec \bs{x}-\bs{u} \\ \bs{y}-\bs{z} \succ \bs{y}-\bs{v}}}\mathbb P_p\big( X \cap \{\bs{o} \rightarrow [\bs{w},\bs{x}\rangle\} \cap \{\bs{o}\rightarrow [ \bs{z}, \bs{y} \rangle \}\big)\nn\\
&=p^4D(\bs{x}-\bs{u})D(\bs{y}-\bs{v})\sum_{\substack{\bs{w},\bs{z}\\ \bs{x}-\bs{w}\prec \bs{x}-\bs{u} \\ \bs{y}-\bs{z} \succ \bs{y}-\bs{v}}}D(\bs{x}-\bs{w})D(\bs{y}-\bs{z})\mathcal P_p^{\bs{u},\bs{v}}(\bs{w},\bs{z}).
\end{align}
Then, by considering the cases where $\bs{w}=\bs{z}$, $\bs{w}=\bs{v}$, and $\bs{z}=\bs{u}$, we obtain
\begin{align}
\sum_{\substack{\bs{u},\bs{v}\\ \bs{x}-\bs{u}\prec \bs{x}-\bs{v}}}\mathcal E''^{(1)}_{p_c}(\substack{\sss [\bs{u},\bs{x}\rangle \\ \sss [\bs{v},\bs{x}\rangle})&=\frac{1}{8}F_{9,p}(\bs{x},\bs{x}),\lbeq{F8xxF9xx}\\
\sum_{\substack{\bs{u},\bs{v}\\ \bs{u} \neq \bs{v}}}\mathcal E''^{(1)}_{p_c}(\substack{\sss [\bs{u},\bs{x}\rangle \\ \sss [\bs{v},\bs{y}\rangle})&=\frac{1}{4}\big[F_{7,p}(\bs{x},\bs{y})+F_{8,p}(\bs{x},\bs{y})+F_{9,p}(\bs{x},\bs{y})+F_{10,p}(\bs{x},\bs{y})\big]\lbeq{F78910xy}
\end{align}
where
\begin{align}
&F_{7,p}(\bs{x},\bs{y})=\raisebox{-15pt}{\includegraphics[scale=0.85]{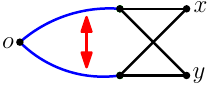}},&&F_{8,p}(\bs{x},\bs{y})=\raisebox{-15pt}{\includegraphics[scale=0.85]{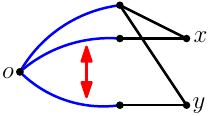}}\\
&F_{9,p}(\bs{x},\bs{y})=\raisebox{-25pt}{\includegraphics[scale=0.85]{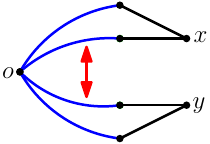}},&&F_{10,p}(\bs{x},\bs{y})=\raisebox{-25pt}{\includegraphics[scale=0.85]{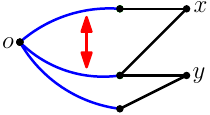}}+\raisebox{-25pt}{\includegraphics[scale=0.85]{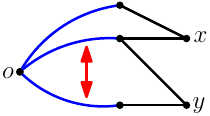}}.
\end{align}

\subsubsection{Estimation of $\mathcal E^{(3)}_p$, $\mathcal E'^{(3)}_p$ and $E''^{(2)}_p$ }
The terms $\mathcal E^{(3)}_p$, $\mathcal E'^{(3)}_p$ and $E''^{(2)}_p$, which contribute to the error, need to be estimated only in the case where $\bs{x}(t)=2$. First, we note that in the case where $\bs{x}(t)=\bs{y}(t)=2$, the time components of all neighbor of $\bs{x}$ or $\bs{y}$ are all equal to $1$. Therefore, when $\bs{u}, \bs{v}, \bs{w}, \bs{z}$ and $\bs{r}$ are pairwise distinct, we have
\begin{align}
\mathcal P_p^{\bs{u},\bs{v}}&=p^2D(\bs{u})D(\bs{v}),\lbeq{ouovDD}\\
\mathcal P_p^{\bs{u},\bs{v}}(\bs{w})&=p^3D(\bs{u})D(\bs{v})D(\bs{w}), \lbeq{ouovowDDD}\\
\mathcal P_p^{\bs{u},\bs{v}}(\bs{w},\bs{z})&=p^4D(\bs{u})D(\bs{v})D(\bs{w})D(\bs{z})\lbeq{ouovowozDDD},\\
\mathcal P_p^{\bs{u},\bs{v}}(\bs{w},\bs{z},\bs{r})&=p^5D(\bs{u})D(\bs{v})D(\bs{w})D(\bs{z})D(\bs{r}),\lbeq{DDDD}
\end{align}
where we use that $\{\bs{o}\rightarrow \bs{u}\}\circ\{\bs{o}\rightarrow \bs{v}\}=\{[\bs{o},\bs{u}\rangle, [\bs{o}, \bs{v}\rangle~\text{open}\}$ for $\bs{u}, \bs{v}$ satisfying $\bs{u}(t)=\bs{v}(t)=1$ and $\bs{u}\neq\bs{v}$. By recalling the definitions of $\mathcal E^{(3)}_p$, $\mathcal E'^{(3)}_p$: 
\begin{align}
\mathcal E^{(3)}_p(\substack{\sss [\bs{u},\bs{x}\rangle \\ \sss [\bs{v},\bs{y}\rangle})
&=p^5D(\bs{x}-\bs{u})D(\bs{y}-\bs{v})\sum_{\substack{\bs{w}, \bs{z},\bs{r}\\ \bs{x}-\bs{w}\prec\bs{x}-\bs{u} \\ \bs{x}-\bs{z}\prec\bs{x}-\bs{w}\\ \bs{x}-\bs{r}\prec\bs{x}-\bs{z}}}D(\bs{x}-\bs{w})D(\bs{x}-\bs{z})D(\bs{x}-\bs{r})\mathcal P_p^{\bs{u},\bs{v}}(\bs{w},\bs{z},\bs{r}),\\
\mathcal E'^{(3)}_p(\substack{\sss [\bs{u},\bs{x}\rangle \\ \sss [\bs{v},\bs{y}\rangle})
&=p^5D(\bs{x}-\bs{u})D(\bs{y}-\bs{v})\sum_{\substack{\bs{w}, \bs{z},\bs{r}\\ \bs{y}-\bs{w}\succ\bs{y}-\bs{v} \\ \bs{y}-\bs{z}\succ\bs{y}-\bs{w}\\ \bs{y}-\bs{r}\succ\bs{y}-\bs{z}}}D(\bs{y}-\bs{w})D(\bs{y}-\bs{z})D(\bs{y}-\bs{r})\mathcal P_p^{\bs{u},\bs{v}}(\bs{w},\bs{z},\bs{r}),
\end{align}
we arrive at 
\begin{align}
\sum_{\substack{\bs{u},\bs{v}\\ \bs{x}-\bs{u}\prec \bs{x}-\bs{v}}}\Big[\mathcal E^{(3)}_{p_c}(\substack{\sss [\bs{u},\bs{x}\rangle \\ \sss [\bs{v},\bs{x}\rangle})+\mathcal E'^{(3)}_{p_c}(\substack{\sss [\bs{u},\bs{x}\rangle \\ \sss [\bs{v},\bs{x}\rangle})\Big]&=\frac{1}{16}F_{6,p}(\bs{x},\bs{x}),\lbeq{F6xx}\\
\sum_{\substack{\bs{u},\bs{v}\\ \bs{u} \neq \bs{v}}}\Big[\mathcal E^{(3)}_{p_c}(\substack{\sss [\bs{u},\bs{x}\rangle \\ \sss [\bs{v},\bs{y}\rangle})+\mathcal E'^{(3)}_{p_c}(\substack{\sss [\bs{u},\bs{x}\rangle \\ \sss [\bs{v},\bs{y}\rangle})\Big]&=\frac{3}{8}F_{5,p}(\bs{x},\bs{y})+\frac{1}{8}F_{6,p}(\bs{x},\bs{y})
\end{align}
where
\begin{align}
F_{5,p}(\bs{x},\bs{y})&=\raisebox{-20pt}{\includegraphics[scale=0.75]{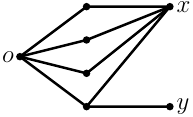}}+\raisebox{-20pt}{\includegraphics[scale=0.75]{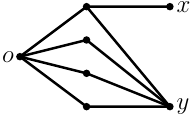}}\le p_c^9sD^{*2}(\bs{x})D^{*2}(\bs{y})\big[D^{*2}(\bs{x})^2+D^{*2}(\bs{y})^2\big], \lbeq{F5}\\
F_{6,p}(\bs{x},\bs{y})&=\raisebox{-25pt}{\includegraphics[scale=0.7]{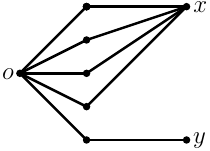}}+\raisebox{-25pt}{\includegraphics[scale=0.7]{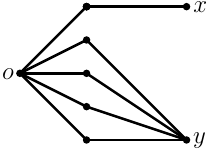}}\le p_c^{10}D^{*2}(\bs{x})D^{*2}(\bs{y})\big[D^{*2}(\bs{x})^3+D^{*2}(\bs{y})^3\big].\lbeq{F6}
\end{align}

Furthermore, since $E''^{(2)}_p$ is defined as
\begin{align}\lbeq{E''2bd}
E''^{(2)}_p(\substack{\sss [\bs{u},\bs{x}\rangle \\ \sss [\bs{v},\bs{y}\rangle})&=\sum_{\substack{\omega_1':\bs{o} \rightarrow \bs{u}\\ \omega'_2:\bs{o}\rightarrow \bs{v}\\ \omega'_1\cap \omega'_2=\varnothing}}\sum_{\substack{\bs{w},\bs{z}\\ \bs{x}-\bs{w}\prec \bs{x}-\bs{u} \\ \bs{y}-\bs{z} \succ \bs{y}-\bs{v}}}\mathbb P_p\big( X \cap \{\bs{o} \rightarrow [\bs{w},\bs{x}\rangle\} \cap \{\bs{o}\rightarrow [ \bs{z}, \bs{y} \rangle \} \cap (V_{\sss[\bs{w},\bs{x}\rangle} \cap V'_{\sss[ \bs{z}, \bs{y} \rangle})^c\big),
\end{align}
by applying \refeq{Voxc}--\refeq{Voyc}, we obtain
\begin{align}\lbeq{E''}
\sum_{\substack{\bs{u},\bs{v}\\ \bs{x}-\bs{u}\prec \bs{x}-\bs{v}}}E''^{(2)}_p(\substack{\sss [\bs{u},\bs{x}\rangle \\ \sss [\bs{v},\bs{x}\rangle})\le\frac{1}{8}F_{6,p}(\bs{x},\bs{x})\le \frac{p_c^{10}}{8}D^{*2}(\bs{x})^5
\end{align}
and
\begin{align}\lbeq{E''2bd}
\sum_{\substack{\bs{u},\bs{v}\\ \bs{u}\neq\bs{v}}}E''^{(2)}_p(\substack{\sss [\bs{u},\bs{x}\rangle \\ \sss [\bs{v},\bs{y}\rangle})
&\le p^5D(\bs{x}-\bs{u})D(\bs{y}-\bs{v})\sum_{\substack{\bs{w},\bs{z}\\ \bs{x}-\bs{w}\prec \bs{x}-\bs{u} \\ \bs{y}-\bs{z} \succ \bs{y}-\bs{v}}}D(\bs{x}-\bs{w})D(\bs{y}-\bs{z})\nn\\
&\times \Big[\sum_{\substack{\bs{r}\\ \bs{x}-\bs{r}\prec \bs{x}-\bs{w} }}D(\bs{x}-\bs{r})\mathcal P_p^{\bs{u},\bs{v}}(\bs{w},\bs{z},\bs{r})+\sum_{\substack{\bs{r}\\\bs{y}-\bs{r}\succ \bs{y}-\bs{z} }}D(\bs{y}-\bs{r})\mathcal P_p^{\bs{u},\bs{v}}(\bs{w},\bs{z},\bs{r})\Big]\nn\\
&\le C\Bigg[\raisebox{-18pt}{\includegraphics[scale=0.73]{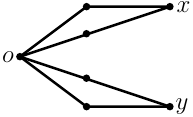}}+\raisebox{-18pt}{\includegraphics[scale=0.73]{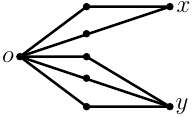}}+\raisebox{-18pt}{\includegraphics[scale=0.73]{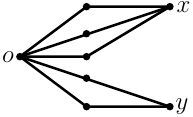}}+\raisebox{-18pt}{\includegraphics[scale=0.73]{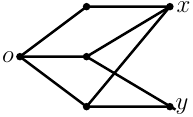}}\nn\\
&\hskip10mm+\raisebox{-18pt}{\includegraphics[scale=0.73]{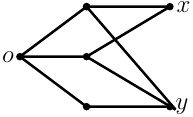}}+\raisebox{-18pt}{\includegraphics[scale=0.73]{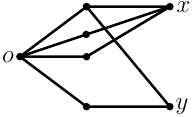}}+\raisebox{-18pt}{\includegraphics[scale=0.73]{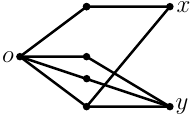}}~\Bigg]\nn\\
&\le Cs^2D^{*2}(\bs{x})D^{*2}(\bs{y}). 
\end{align}

\section{Analysis of $\Pi_{p_c}^{(0)}(\bs{x})$}\label{SectionPi0}
In this section, we investigate $\Pi_{p_c}^{(0)}(\bs{x})$. Since Lemma~\ref{doubleconnection} gives 
\begin{align}\lbeq{pi0x3}
\Pi_{p_c}^{\sss(0)}(\bs{x})=\mathcal P_{p_c}^{\bs{x},\bs{x}}=p_c^2\sum_{\substack{\bs{u},\bs{v}\\ \bs{x}-\bs{u}\prec \bs{x}-\bs{v}}}D(\bs{x}-\bs{u})D(\bs{x}-\bs{v})\mathcal P_{p_c}^{\bs{u},\bs{v}}-\mathcal H_{p_c}(\bs{x},\bs{x}),
\end{align}
we estimate the first term and the remainder term $\mathcal H_{p_c}(\bs{x},\bs{x})$ separately. We consider the case where $\bs{x}(t)=2$, $\bs{x}(t)=3$ and $\bs{x}(t)=4$, in Section~\ref{t2case}, \ref{3case} and \ref{case4},  respectively. 

\subsection{The case where the time component equals 2}\label{t2case}
In this section, we estimate $\Pi_{p_c}^{(0)}(\bs{x})$ for $\bs{x}(t)=2$ by proving the following lemma. 
\begin{prp}\label{prp: Pi0xxt2}For $\bs{x}$ such that $\bs{x}(t)=2$,  
\begin{align}\lbeq{prpPi0xxt2}
\Pi_{p_c}^{\sss(0)}(\bs{x})&=
\begin{cases}
\frac{p_c^4}{2}s^2-\frac{p_c^4}{2}[1+p_c^2]s^3+\frac{15}{8}s^4+O(d^{-5})&[\bs{x}(s)=o],\\
p_c^4s^4&[\bs{x}(s)=e_1+e_2].
\end{cases}
\end{align}
Then, we have
\begin{align}\lbeq{sumofpixt2}
\sum_{\bs{x}:\bs{x}(t)=2}\Pi_{p_c}^{\sss(0)}(\bs{x})=p_c^4s^2-\frac{p_c^4}{2}[3+p_c^2]s^3+\frac{15}{8}s^4+O(d^{-5}).
\end{align}
Moreover, for $\bs{x}$ and $\bs{y}$ such that $\bs{x}(t)=\bs{y}(t)=2$, 
\begin{align}\lbeq{Pxyx2pic}
\mathcal P_p^{\bs{x},\bs{y}}
&=\big[1+O(d^{-2})\big]\raisebox{-13pt}{\includegraphics[scale=0.75]{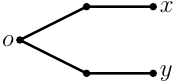}}-\frac{1}{2}\bigg[\raisebox{-15pt}{\includegraphics[scale=0.8]{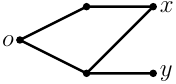}}+\raisebox{-15pt}{\includegraphics[scale=0.8]{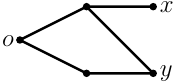}}~\bigg]\nn\\
&-\frac{1}{2}\bigg[\raisebox{-13pt}{\includegraphics[scale=0.75]{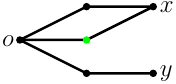}}+\raisebox{-13pt}{\includegraphics[scale=0.75]{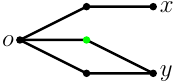}}~\bigg].
\end{align}
\end{prp}
In the above, an unlabelled green vertex is also summed over all possible vertices but it is allowed to coincide with other vertices that share the same time component. 

Before proving Proposition~\ref{prp: Pi0xxt2}, we first establish the following lemma, which provides an estimate of $\mathcal H_{p_c}(\bs{x},\bs{x})$ when $\bs{x}(t)=2$. 

\begin{lmm}\label{lmm: Hxyxt2}
For $\bs{x}$ and $\bs{y}$ such that $\bs{x}(t)=\bs{y}(t)=2$, 
\begin{align}
\mathcal H_{p_c}(\bs{x},\bs{x})&=
\begin{cases}
\frac{p_c^6}{2}s^3-\frac{15}{8}s^4+O(d^{-5})&[\bs{x}(s)=o],\\
0&[\bs{x}(s)=e_1+e_2],
\end{cases}
\lbeq{mathcalHxx2bound}\\
\mathcal H_{p_c}(\bs{x},\bs{y})
&=\frac{1}{2}\bigg[\raisebox{-13pt}{\includegraphics[scale=0.75]{Hxy21}}+\raisebox{-13pt}{\includegraphics[scale=0.75]{Hxy22}}~\bigg]+\frac{1}{2}\bigg[\raisebox{-13pt}{\includegraphics[scale=0.75]{Hxy23}}+\raisebox{-13pt}{\includegraphics[scale=0.75]{Hxy24}}~\bigg]\nn\\
&+O(d^{-2})D^{*2}(\bs{x})D^{*2}(\bs{y}).\lbeq{mathcalHxy2bound}\nn\\
\end{align}
\end{lmm}

\Proof{\it Proof of Lemma~\ref{lmm: Hxyxt2}.}First we prove \refeq{mathcalHxx2bound}. Let 
\begin{align}
\text{\bf Er}(\bs{x}(s),2)\le 
\begin{cases}
O(d^{-5})&[\bs{x}(s)=o],\\
O(d^{-7})&[\bs{x}(s)=e_1+e_2].
\end{cases}
\end{align}

By \refeq{Hxxsummandt2}, \refeq{F2xx}, \refeq{F4xx}, \refeq{F8xxF9xx}, we have 
\begin{align}\lbeq{Hpcxxt2}
\mathcal H_{p_c}(\bs{x},\bs{x})&=\frac{1}{4}F_{2,p_c}(\bs{x},\bs{x})-\frac{1}{8}F_{4,{p_c}}(\bs{x},\bs{x})-\frac{1}{8}F_{9,{p_c}}(\bs{x},\bs{x})\nn\\
&+\sum_{\substack{\bs{u},\bs{v}\\ \bs{u}-\bs{x}\prec \bs{v}-\bs{x}}}\left[E^{(3)}_{p_c}(\substack{\sss [\bs{u},\bs{x}\rangle \\ \sss [\bs{v},\bs{x}\rangle})+E'^{(3)}_{p_c}(\substack{\sss [\bs{u},\bs{x}\rangle \\ \sss [\bs{v},\bs{x}\rangle})\right]+\sum_{\substack{\bs{u},\bs{v}\\ \bs{u}-\bs{x}\prec \bs{v}-\bs{x}}}E''^{(2)}_{p_c}(\substack{\sss [\bs{u},\bs{x}\rangle \\ \sss [\bs{v},\bs{x}\rangle}). 
\end{align}
By $E^{(n)}_{p_c}\le \mathcal E^{(n)}_{p_c}$, $E'^{(n)}_{p_c}\le \mathcal E'^{(n)}_{p_c}$ and \refeq{F6xx}, we obtain
\begin{align}\lbeq{dwdmsldnjskndkjsdjjs}
\sum_{\substack{\bs{u},\bs{v}\\ \bs{u}-\bs{x}\prec \bs{v}-\bs{x}}}\left[E^{(3)}_{p_c}(\substack{\sss [\bs{u},\bs{x}\rangle \\ \sss [\bs{v},\bs{x}\rangle})+E'^{(3)}_{p_c}(\substack{\sss [\bs{u},\bs{x}\rangle \\ \sss [\bs{v},\bs{x}\rangle})\right]\le \sum_{\substack{\bs{u},\bs{v}\\ \bs{u}-\bs{x}\prec \bs{v}-\bs{x}}}\left[\mathcal E^{(3)}_{p_c}(\substack{\sss [\bs{u},\bs{x}\rangle \\ \sss [\bs{v},\bs{x}\rangle})+\mathcal E'^{(3)}_{p_c}(\substack{\sss [\bs{u},\bs{x}\rangle \\ \sss [\bs{v},\bs{x}\rangle})\right]=\frac{1}{16}F_{6,{p_c}}(\bs{x},\bs{x}).
\end{align}
We note that $F_{i,p_c}(\bs{x},\bs{x})=0$ for $\bs{x}(s)=e_1+e_2$ and $i=2,4,6,9$. 
Therefore, by \refeq{E''}, we obtain $\mathcal H_{p_c}(\bs{x},\bs{x})=0$ for $\bs{x}(s)=e_1+e_2$. 

The remaining task is to estimate the followings for $\bs{x}=(o,2)$: 
\begin{align}
F_{2,p}(\bs{x},\bs{x})&=2p^6D^{*2}(\bs{x})^3-6p^6D^{*2}(\bs{x})\sum_{\bs{u}}D(\bs{x}-\bs{u})^2D(\bs{u})^2+4p^6\sum_{\bs{u}}D(\bs{x}-\bs{u})^3D(\bs{u})^3,\lbeq{F2x2}\\
F_{4,p}(\bs{x},\bs{x})&=2p^8D^{*2}(\bs{x})^4+\text{\bf Er}(\bs{x}(s),2)\lbeq{F4x2}\\
F_{9,p}(\bs{x},\bs{x})&=p^8D^{*2}(\bs{x})^4+\text{\bf Er}(\bs{x}(s),2).\lbeq{F9est}
\end{align}
By \refeq{E''} and \refeq{dwdmsldnjskndkjsdjjs} the last two terms in the right-hand side of \refeq{Hpcxxt2} are $\text{\bf Er}(\bs{x}(s),2)$. We note that for $\bs{x}=(\bs{x}(s),2)$ and for $m,n\ge 1$, 
\begin{align}\lbeq{DmDn}
\sum_{\substack{\bs{u} \\ \bs{u}(t)=1}}D(\bs{x}-\bs{u})^mD(\bs{u})^n
=
\begin{cases}
s^{m+n-1}&[\bs{x}(s)=o],\\
2s^{m+n}&[\bs{x}(s)=e_1+e_2],\\
s^{m+n}&[\bs{x}(s)=2e_1],
\end{cases}
\end{align}
where the special case of $m=n=1$ corresponds to $D^{*2}(\bs{x})$. Therefore, using $|p_c-1|\le O(d^{-1})$, we obtain \refeq{mathcalHxx2bound} for $\bs{x}(s)=o$.

As for \refeq{mathcalHxy2bound}, we establish that there exists a positive constant $C$ such that 
\begin{align}\lbeq{Hxyt2arror}
\Big|\mathcal H_{p_c}(\bs{x},\bs{y})-\sum_{\substack{\bs{u},\bs{v}\\ \bs{u} \neq \bs{v}}}\left[\mathcal E^{(1)}_{p_c}(\substack{\sss [\bs{u},\bs{x}\rangle \\ \sss [\bs{v},\bs{y}\rangle})+\mathcal E'^{(1)}_{p_c}(\substack{\sss [\bs{u},\bs{x}\rangle \\ \sss [\bs{v},\bs{y}\rangle})\right]\Big|&\le \sum_{\substack{\bs{u},\bs{v}\\ \bs{u} \neq \bs{v}}}\left[\mathcal E^{(2)}_{p_c}(\substack{\sss [\bs{u},\bs{x}\rangle \\ \sss [\bs{v},\bs{y}\rangle})+\mathcal E'^{(2)}_{p_c}(\substack{\sss [\bs{u},\bs{x}\rangle \\ \sss [\bs{v},\bs{y}\rangle})\right]+\sum_{\substack{\bs{u},\bs{v}\\ \bs{u} \neq \bs{v}}}\mathcal E''^{(1)}_{p_c}(\substack{\sss [\bs{u},\bs{x}\rangle \\ \sss [\bs{v},\bs{y}\rangle})\nn\\
&+\sum_{\substack{\bs{u},\bs{v}\\ \bs{u} \neq \bs{v}}}\left[\mathcal E^{(3)}_{p_c}(\substack{\sss [\bs{u},\bs{x}\rangle \\ \sss [\bs{v},\bs{y}\rangle})+\mathcal E'^{(3)}_{p_c}(\substack{\sss [\bs{u},\bs{x}\rangle \\ \sss [\bs{v},\bs{y}\rangle})\right]+\sum_{\substack{\bs{u},\bs{v}\\ \bs{u} \neq \bs{v}}}E''^{(2)}_{p_c}(\substack{\sss [\bs{u},\bs{x}\rangle \\ \sss [\bs{v},\bs{y}\rangle})\nn\\
&\le Cs^2\raisebox{-12pt}{\includegraphics[scale=0.75]{Pxy2main}}.
\end{align}
For the first inequality in \refeq{Hxyt2arror}, we use \refeq{Hxxsummandt2} and $E^{(3)}_{p_c}+E'^{(3)}_{p_c}\le \mathcal E^{(3)}_{p_c}+\mathcal E'^{(3)}_{p_c}$. For the second inequality, we use \refeq{F3xyF4xy}, \refeq{F78910xy} and  \refeq{ouovDD}--\refeq{DDDD} to estimate the first two terms, and use \refeq{F5}--\refeq{F6} and \refeq{E''2bd} for the remaining terms.  

Since $F_{1,p_c}(\bs{x},\bs{y})$ can be represented diagramatically as
\begin{align}
F_{1,p_c}(\bs{x},\bs{y})=\frac{1}{2}\bigg[\raisebox{-13pt}{\includegraphics[scale=0.75]{Hxy21}}+\raisebox{-13pt}{\includegraphics[scale=0.75]{Hxy22}}~\bigg],
\end{align}
by \refeq{F1xyF2xy}, the remaining task is to show evaluate 
\begin{align}\lbeq{F2minusmain}
\bigg|F_{2,p_c}(\bs{x},\bs{y})-\bigg[\raisebox{-13pt}{\includegraphics[scale=0.75]{Hxy23}}+\raisebox{-13pt}{\includegraphics[scale=0.75]{Hxy24}}\bigg]\bigg|\le Cs^2D^{*2}(\bs{x})D^{*2}(\bs{y}). 
\end{align}
By \refeq{F2def}, the left-hand side of \refeq{F2minusmain} can be bounded diagrammatically above by
\begin{align}
\raisebox{-15pt}{\includegraphics[scale=0.85]{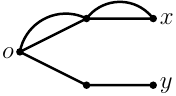}}+\raisebox{-15pt}{\includegraphics[scale=0.85]{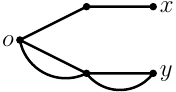}}\le 2p_c^2s^2D^{*2}(\bs{x})D^{*2}(\bs{y}).
\end{align}
This completes the proof of \refeq{mathcalHxy2bound}, and thus proof of Lemma~\ref{lmm: Hxyxt2}. 
\QED

\Proof{\it Proof of Proposition~\ref{prp: Pi0xxt2}.}
We have
\begin{align}\lbeq{sumofpi0x2}
\sum_{\bs{x}:\bs{x}(t)=2}\Pi_{p_c}^{\sss(0)}(\bs{x})=\Pi_{p_c}^{\sss(0)}(o,2)+2d(d-1)\Pi_{p_c}^{\sss(0)}(e_1+e_2,2),
\end{align}
which, togtether with \refeq{prpPi0xxt2}, yields \refeq{sumofpixt2}. 

By \refeq{ouovDD} and \refeq{DmDn}, the first term of \refeq{pi0x3} can be rewritten as 
\begin{align}\lbeq{Pi02xxmain}
p_c^4\sum_{\substack{\bs{u},\bs{v}\\ \bs{x}-\bs{u}\prec \bs{x}-\bs{v}}}D(\bs{x}-\bs{u})D(\bs{x}-\bs{v})D(\bs{u})D(\bs{v})&=\frac{1}{2}\raisebox{-13pt}{\includegraphics[scale=0.85]{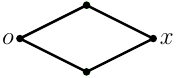}}\nn\\
&=\frac{p_c^4}{2}\Big[D^{*2}(\bs{x})^2-\sum_{\substack{\bs{u}\\ \bs{u}(t)=1}}D(\bs{x}-\bs{u})^2D(\bs{u})^2\Big]\nn\\
&=
\begin{cases}
\frac{p_c^4s^2}{2}[1-s]&[\bs{x}(s)=o]\\
p_c^4s^4&[\bs{x}(s)=e_1+e_2].
\end{cases}
\end{align}
Combining \refeq{mathcalHxx2bound}, \refeq{pi0x3} and \refeq{Pi02xxmain} yields \refeq{prpPi0xxt2}.

Finally, combining \refeq{lmm:oxoy}, \refeq{ouovDD} and \refeq{mathcalHxy2bound} yields \refeq{Pxyx2pic}. 
\QED

\subsection{The case where the time component equals 3}\label{3case}
This section provides estimates $\Pi_p^{(0)}(\bs{x})$ for $\bs{x}(t)=3$, establishing the following proposition. 
\begin{prp}\label{prp: Pi0xxt3}For $\bs{x}$ such that $\bs{x}(t)=3$, 
\begin{align}\lbeq{prpPi0xxt3}
\Pi_{p_c}^{\sss(0)}(\bs{x})&=
\begin{cases}
\frac{7p_c^6}{2}s^4-16s^5+O(d^{-6})&[\bs{x}(s)=e_1],\\
s^6+O(d^{-7})&[\bs{x}(s)=2e_1+e_2],\\
9p_c^6s^6-3s^7+O(d^{-8})&[\bs{x}(s)=e_1+e_2+e_3].
\end{cases}
\end{align}
Then, we have
\begin{align}\lbeq{sumPi0x3}
\sum_{\bs{x}:\bs{x}(t)=3}\Pi_{p_c}^{\sss(0)}(\bs{x})=5p_c^6s^3-25s^4+O(d^{-5}). 
\end{align}
Moreover, for $\bs{x}$ and $\bs{y}$ such that $\bs{x}(t)=\bs{y}(t)=3$, 
\begin{align}\lbeq{LLL}
\Big|\mathcal P_p^{\bs{x},\bs{y}}-\raisebox{-11pt}{\includegraphics[scale=0.75]{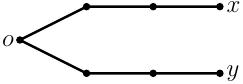}}\Big|
&\le CsD^{*3}(\bs{x})D^{*3}(\bs{y}). 
\end{align}
\end{prp}

As a step toward proving Proposition~\ref{prp: Pi0xxt3}, we first establish the following lemma. 
\begin{lmm}\label{lmm: Hxyxt3}
For $\bs{x}$ and $\bs{y}$ such that $\bs{x}(t)=\bs{y}(t)=3$, 
\begin{align}
\mathcal H_{p_c}(\bs{x},\bs{x})&=
\begin{cases}
2s^5+O(d^{-6})&[\bs{x}(s)=e_1]\\
O(d^{-8})&[\bs{x}(s)=2e_1+e_2]\\
O(d^{-8})&[\bs{x}(s)=e_1+e_2+e_3],
\end{cases}\lbeq{mathcalHxx3bound}\\
\mathcal H_{p_c}(\bs{x},\bs{y})&\le CsD^{*3}(\bs{x})D^{*3}(\bs{y}).\lbeq{mathcalHxy3bound}
\end{align}
\end{lmm}
To prove Lemma~\ref{lmm: Hxyxt3}, we require the following two lemmas, Lemma~\ref{lmm:ouovowestimate} and \ref{lmm:osxoy}, which together yield estimates of the propbabilities $\mathcal P_p^{\bs{u},\bs{v}}(\bs{w})$ and $\mathcal P_p^{\bs{u},\bs{v}}(\bs{w},\bs{z})$. 
When the time comopnents of $\bs{u}, \bs{v}, \bs{w}$ and $\bs{z}$ are equal to $2$, $\mathcal P_p^{\bs{u},\bs{v}}(\bs{w})$ and $\mathcal P_p^{\bs{u},\bs{v}}(\bs{w},\bs{z})$ cannot be simplified, in contrast to the case where their time components are 1. We will prove Lemma~\ref{lmm:ouovowestimate} and \ref{lmm:osxoy} in Appendix~\ref{Appendix}. 
\begin{lmm}\label{lmm:ouovowestimate}
Let $\bs{u},\bs{v}$ be points such that $\bs{u}(t)=\bs{v}(t)=2$, and  let $\bs{w}$ and $\bs{z}$ be distincts points with $\bs{w}(t)=\bs{z}(t)=2$, where $\bs{w}\neq \bs{u}, \bs{v}$ and $\bs{z}\neq \bs{u}, \bs{v}$. Then, for $p\le p_c$, 
\begin{align}\lbeq{ouovowestimate}
\Big|\mathcal P_p^{\bs{u},\bs{v}}(\bs{w})-\zeta_p(\bs{u},\bs{v},\bs{w})\Big|\le 
p^6D^{*2}(\bs{u})D^{*2}(\bs{v})D^{*2}(\bs{w})
\end{align}
where
\begin{align}
\mathcal P_{p,[\bs{o},\bs{s}\rangle}^{\bs{x},\bs{y}}&:=\mathbb P_p\big(\{[\bs{o},\bs{s}\rangle \rightarrow \bs{x} \}\circ \{\bs{o}\rightarrow \bs{y}\}\big), \\
\zeta_p(\bs{u},\bs{v},\bs{w})&:=\frac{2-\delta_{\bs{u},\bs{v}}}{2}p\sum_{\substack{\bs{s} \\ \bs{s}(t)=1}}D(\bs{w}-\bs{s})\big[\mathcal P_{p,[\bs{o},\bs{s}\rangle}^{\bs{u},\bs{v}}+\mathcal P_{p,[\bs{o},\bs{s}\rangle}^{\bs{v},\bs{u}}\big].\lbeq{wsosuovouosv}
\end{align}
Furthermore, in the case $\bs{u}=\bs{v}$, 
\begin{align}\lbeq{ouovowestimate1}
\Big|\mathcal P_p^{\bs{u},\bs{u}}(\bs{w})-\Big[\zeta_p(\bs{u},\bs{u},\bs{w})+\frac{1}{2}\raisebox{-12pt}{\includegraphics[scale=0.75]{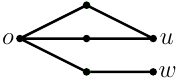}}~\Big]\Big|\le 
p^8sD^{*2}(\bs{u})^2D^{*2}(\bs{w})
\end{align}
Additionally, 
\begin{align}\lbeq{Puvwz3}
\mathcal P_p^{\bs{u},\bs{v}}(\bs{w},\bs{z})\le 
\begin{cases}
Cs^2D^{*2}(\bs{u})D^{*2}(\bs{v})&[\bs{u} \neq\bs{v}],\\
Cs^2D^{*2}(\bs{u})D^{*2}(\bs{w})&[\bs{u} =\bs{v}].
\end{cases}
\end{align}
\end{lmm}
\begin{lmm}\label{lmm:osxoy}
Let $\bs{x}$ and $\bs{y}$ be points such that $\bs{x}(t)=\bs{y}(t)=2$. For the both case where $\bs{x}\neq \bs{y}$ and $\bs{x}=\bs{y}$, 
\begin{align}\lbeq{osxoy2}
\Big|\mathcal P_{p,[\bs{o},\bs{s}\rangle}^{\bs{x},\bs{y}}-\raisebox{-12pt}{\includegraphics[scale=0.75]{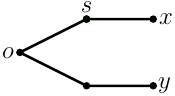}}\Big|\le \frac{1}{2}\raisebox{-12pt}{\includegraphics[scale=0.75]{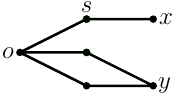}}
\end{align}
\end{lmm}
\Proof{\it Proof of Lemma~\ref{lmm: Hxyxt3}, assuming Lemma~\ref{lmm:ouovowestimate} and \ref{lmm:osxoy}.}
We first consider \refeq{mathcalHxx3bound}. By \refeq{Hxxsummandt3} and \refeq{F2xx}, we have 
\begin{align}\lbeq{Hpcxxx3}
\mathcal H_{p_c}(\bs{x},\bs{x})&=\frac{1}{4}F_{2,p_c}(\bs{x},\bs{x})-\sum_{\substack{\bs{u},\bs{v}\\ \bs{x}-\bs{u}\prec \bs{x}-\bs{v}}}\left[E^{(2)}_{p_c}(\substack{\sss [\bs{u},\bs{x}\rangle \\ \sss [\bs{v},\bs{x}\rangle})+E'^{(2)}_{p_c}(\substack{\sss [\bs{u},\bs{x}\rangle \\ \sss [\bs{v},\bs{x}\rangle})\right]-\sum_{\substack{\bs{u},\bs{v}\\ \bs{x}-\bs{u}\prec \bs{x}-\bs{v}}}E''^{(1)}_{p_c}(\substack{\sss [\bs{u},\bs{x}\rangle \\ \sss [\bs{v},\bs{x}\rangle}).
\end{align}
By \refeq{F4xx}, \refeq{F4def} and \refeq{F8xxF9xx}, the second term and the third term can be bounded resprectively as
\begin{align}
\sum_{\substack{\bs{u},\bs{v}\\ \bs{x}-\bs{u}\prec \bs{x}-\bs{v}}}\left[E^{(2)}_{p_c}(\substack{\sss [\bs{u},\bs{x}\rangle \\ \sss [\bs{v},\bs{x}\rangle})+E'^{(2)}_{p_c}(\substack{\sss [\bs{u},\bs{x}\rangle \\ \sss [\bs{v},\bs{x}\rangle})\right]&\le \sum_{\substack{\bs{u},\bs{v}\\ \bs{x}-\bs{u}\prec \bs{x}-\bs{v}}}\left[\mathcal E^{(2)}_{p_c}(\substack{\sss [\bs{u},\bs{x}\rangle \\ \sss [\bs{v},\bs{x}\rangle})+\mathcal E'^{(2)}_{p_c}(\substack{\sss [\bs{u},\bs{x}\rangle \\ \sss [\bs{v},\bs{x}\rangle})\right]=\frac{1}{8}F_{4,p_c}(\bs{x},\bs{x}),\lbeq{F4xxx3}\\
\sum_{\substack{\bs{u},\bs{v}\\ \bs{x}-\bs{u}\prec \bs{x}-\bs{v}}}E''^{(1)}_{p_c}(\substack{\sss [\bs{u},\bs{x}\rangle \\ \sss [\bs{v},\bs{x}\rangle})&\le \sum_{\substack{\bs{u},\bs{v}\\ \bs{x}-\bs{u}\prec \bs{x}-\bs{v}}}\mathcal E''^{(1)}_{p_c}(\substack{\sss [\bs{u},\bs{x}\rangle \\ \sss [\bs{v},\bs{x}\rangle})=\frac{1}{16}F_{4,p_c}(\bs{x},\bs{x}),\lbeq{F89xxx3}
\end{align}
where we use $F_{9,p_c}(\bs{x},\bs{x})=\frac{1}{2}F_{4,p_c}(\bs{x},\bs{x})$. Let 
\begin{align}
\text{\bf Er}(\bs{x}(s),3)\le
\begin{cases}
O(d^{-6})&[\bs{x}(s)=e_1]\\
O(d^{-7})&[\bs{x}(s)=2e_1+e_2]\\
O(d^{-8})&[\bs{x}(s)=e_1+e_2+e_3].
\end{cases}
\end{align}
By \refeq{Puvwz3}, we obtain 
\begin{align}\lbeq{F4x3iserror}
F_{4,p_c}(\bs{x},\bs{x})\le Cs^2D^{*3}(\bs{x})^2\le \text{\bf Er}(\bs{x}(s),3)
\end{align}
where, in the last equality, we use 
\begin{align}\lbeq{D3x3}
D^{*3}(\bs{x})
=\begin{cases}
3s^2\left[1-s\right]&[\bs{x}(s)=e_1]\\
3s^3&[\bs{x}(s)=2e_1+e_2]\\
6s^3&[\bs{x}(s)=e_1+e_2+e_3].
\end{cases}
\end{align}

Thus, the remaining task is to estimate $F_{2,p_c}(\bs{x},\bs{x})$. By \refeq{F2def} and \refeq{ouovowestimate}
we obtain
\begin{align}\lbeq{F2xxminusmain}
F_{2,p_c}(\bs{x},\bs{x})-2p_c^3\sum_{\substack{\bs{u},\bs{v},\bs{w}\\ \bullet}}D(\bs{x}-\bs{u})D(\bs{x}-\bs{v})D(\bs{x}-\bs{w})\zeta_{p_c}(\bs{u},\bs{v},\bs{w})&\le 2p_c^9D^{*3}(\bs{x})^3\nn\\
&\le \text{\bf Er}(\bs{x}(s),3). 
\end{align}

Therefore, it remains to derive an explicit expression for the second term on the left-hand side of \refeq{F2xxminusmain}. 
Therefore, recalling the definition of $\zeta_p(\bs{u},\bs{v},\bs{w})$ given in \refeq{wsosuovouosv}, it suffices to estimate 
\begin{align}\lbeq{DDDzeta}
&2p_c^3\sum_{\substack{\bs{u},\bs{v},\bs{w}\\ \bullet }}D(\bs{x}-\bs{u})D(\bs{x}-\bs{v})D(\bs{x}-\bs{w})\zeta_{p_c}(\bs{u},\bs{v},\bs{w})\nn\\
&=4p_c^4\sum_{\substack{\bs{u},\bs{v},\bs{w}\\ \bullet }}\sum_{\substack{\bs{s}\\\bs{s}(t)=1}}D(\bs{x}-\bs{u})D(\bs{x}-\bs{v})D(\bs{x}-\bs{w})D(\bs{w}-\bs{s})\mathcal P_{p,[\bs{o},\bs{s}\rangle}^{\bs{u},\bs{v}}
\end{align}
Applying \refeq{osxoy2} provides
\begin{align}
&\Big|4p_c^4\sum_{\substack{\bs{u},\bs{v},\bs{w}\\ \bullet }}\sum_{\substack{\bs{s}\\\bs{s}(t)=1}}D(\bs{x}-\bs{u})D(\bs{x}-\bs{v})D(\bs{x}-\bs{w})D(\bs{w}-\bs{s})\mathcal P_{p,[\bs{o},\bs{s}\rangle}^{\bs{u},\bs{v}}-4\raisebox{-11pt}{\includegraphics[scale=0.75]{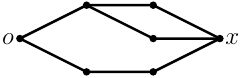}}\Big|\nn\\
&\le Cs^2D^{*3}(\bs{x})^2\nn\\
&\le \text{\bf Er}(\bs{x}(s),3). 
\end{align}

Therefore, to obtain an explicit expression for the leading term of \refeq{DDDzeta}, it suffices to estimate
\begin{align}\lbeq{mainF2xxx3}
4\raisebox{-11pt}{\includegraphics[scale=0.75]{F2main}}
&=4\raisebox{-11pt}{\includegraphics[scale=0.75]{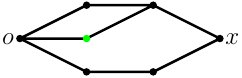}}+\text{\bf Er}(\bs{x}(s),3),
\end{align}
which will be estimated as in \refeq{M2error} (in the proof of Proposition~\ref{prp: Pi0xxt3}). 
Thus, we obtain
\begin{align}\lbeq{F2xxx3mainerror}
F_{2,p_c}(\bs{x},\bs{x})=
\begin{cases}
8s^5+O(d^{-6})&[\bs{x}(s)=e_1],\\
O(d^{-8})&[\bs{x}(s)=2e_1+e_2],\\
O(d^{-8})&[\bs{x}(s)=e_1+e_2+e_3].
\end{cases}
\end{align}

Combining \refeq{Hpcxxx3}, \refeq{F4x3iserror} and \refeq{F2xxx3mainerror} yields \refeq{mathcalHxx3bound}.

For \refeq{mathcalHxy3bound}, since we have $E^{(2)}_p\le \mathcal E^{(2)}_p$, $E'^{(2)}_p\le \mathcal E'^{(2)}_p$ and $E''^{(1)}_p\le \mathcal E''^{(1)}_p$, by using \refeq{F1xyF2xy}, \refeq{F3xyF4xy} and \refeq{F78910xy}, we obtain
\begin{align}
\mathcal H_p(\bs{x},\bs{y})&=\sum_{\substack{\bs{u},\bs{v}\\ \bs{u} \neq \bs{v}}}\bigg\{\left[\mathcal E^{(1)}_p(\substack{\sss [\bs{u},\bs{x}\rangle \\ \sss [\bs{v},\bs{y}\rangle})+\mathcal E'^{(1)}_p(\substack{\sss [\bs{u},\bs{x}\rangle \\ \sss [\bs{v},\bs{y}\rangle})\right]-\left[E^{(2)}_p(\substack{\sss [\bs{u},\bs{x}\rangle \\ \sss [\bs{v},\bs{y}\rangle})+E'^{(2)}_p(\substack{\sss [\bs{u},\bs{x}\rangle \\ \sss [\bs{v},\bs{y}\rangle})\right]-E''^{(1)}_p(\substack{\sss [\bs{u},\bs{x}\rangle \\ \sss [\bs{v},\bs{y}\rangle})\bigg\}\nn\\
&\le \frac{1}{2}\big[F_{1,p_c}(\bs{x},\bs{y})+F_{2,p_c}(\bs{x},\bs{y})\big]+\frac{1}{4}\big[2F_{3,p}(\bs{x},\bs{y})+F_{4,p}(\bs{x},\bs{y})\big]\nn\\
&+\frac{1}{4}\big[F_{7,p}(\bs{x},\bs{y})+F_{8,p}(\bs{x},\bs{y})+F_{9,p}(\bs{x},\bs{y})+F_{10,p}(\bs{x},\bs{y})\big].
\end{align}

By applying the BK inequality to $\mathcal P_p^{\bs{u},\bs{v}}$ and using \refeq{Puvwz3}, we can show 
\begin{align}
F_{i,p_c}(\bs{x},\bs{y})&
\begin{cases}
\le CsD^{*3}(\bs{x})D^{*3}(\bs{y})&(i=1),\\
\le C's^2D^{*3}(\bs{x})D^{*3}(\bs{y})&(i=2,3,4,7,8,9,10),
\end{cases}
\end{align}
which yields \refeq{mathcalHxy3bound}. \QED

Now we prove Proposition~\ref{prp: Pi0xxt3}.
\Proof{\it Proof of Proposition~\ref{prp: Pi0xxt3}.}
We have
\begin{align}\lbeq{sumofpix1x3decom}
\sum_{\bs{x}:\bs{x}(t)=3}\Pi_{p_c}^{\sss(0)}(\bs{x})&=2d\Pi_{p_c}^{\sss(0)}(e_1,3)+2d(d-1)\Pi_{p_c}^{\sss(0)}(2e_1+e_2,3)\nn\\
&+\frac{4}{3}d(d-1)(d-2)\Pi_{p_c}^{\sss(0)}(e_1+e_2+e_3,3).
\end{align}
Since we have \refeq{mathcalHxx3bound}, it suffices to estimate the first term of \refeq{pi0x3} for the case where $\bs{x}(s)=3$. By \refeq{Pxyx2pic}, it can be decomposed as 
\begin{align}\lbeq{mainpi3xx}
&p_c^2\sum_{\substack{\bs{u},\bs{v}\\ \bs{x}-\bs{u}\prec \bs{x}-\bs{v}}}D(\bs{x}-\bs{u})D(\bs{x}-\bs{v})\mathcal P_{p_c}^{\bs{u},\bs{v}}\nn\\
&=\frac{1}{2}\Big[\raisebox{-11pt}{\includegraphics[scale=0.75]{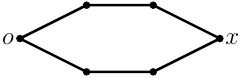}}-\raisebox{-11pt}{\includegraphics[scale=0.75]{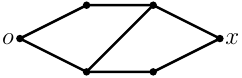}}-\raisebox{-11pt}{\includegraphics[scale=0.75]{M23}}\Big]+\text{\bf Er}(\bs{x}(s),3).
\end{align}
where
\begin{align}
\raisebox{-11pt}{\includegraphics[scale=0.75]{M1}}&=p_c^6D^{*3}(\bs{x})^2-2\underbrace{\raisebox{-17pt}{\includegraphics[scale=0.78]{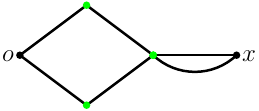}}}_{=p_c^6\eta_1[1,2,1,1](\bs{x})}+\underbrace{\raisebox{-4pt}{\includegraphics[scale=0.78]{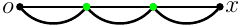}}}_{=p_c^6\eta_3[2](\bs{x})},\lbeq{M1}\\
\raisebox{-11pt}{\includegraphics[scale=0.75]{M2}}&=\underbrace{\raisebox{-11pt}{\includegraphics[scale=0.75]{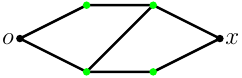}}}_{=p_c^7\eta_2[1](\bs{x})}-2\underbrace{\raisebox{-17pt}{\includegraphics[scale=0.78]{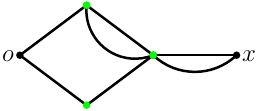}}}_{=p_c^7\eta_1[1,2,1,2](\bs{x})}+\underbrace{\raisebox{-4pt}{\includegraphics[scale=0.78]{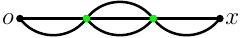}}}_{=p_c^7\eta_3[3](\bs{x})}, \lbeq{ee}\\
\raisebox{-11pt}{\includegraphics[scale=0.75]{M23}}&=p_c^3D^{*3}(\bs{x})\underbrace{\raisebox{-17pt}{\includegraphics[scale=0.78]{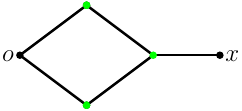}}}_{=p_c^5\eta_1[1,1,1,1](\bs{x})}-\underbrace{\raisebox{-17pt}{\includegraphics[scale=0.75]{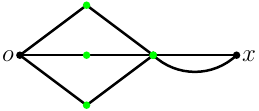}}}_{=p_c^8\eta_1[2,2,1,1](\bs{x})}+\text{\bf Er}(\bs{x}(s),3). \lbeq{M2error}
\end{align}
Here, for $l,m,n,r\ge 1$, $\eta_1[l,m,n,r](\bs{x})$, $\eta_2[l](\bs{x})$ and $\eta_3[l](\bs{x})$ are defined as follow:
\begin{align}
&\eta_1[l,m,n,r]:=\sum_{\substack{\bs{u} \\ \bs{u}(t)=2}}D^{*2}(\bs{u})^lD(\bs{x}-\bs{u})^m\sum_{\substack{\bs{s} \\ \bs{s}(t)=1}}D(\bs{s})^nD(\bs{u}-\bs{s})^r\nn\\
&\hskip20mm=
\begin{cases}
s^{l+m+n+r-1}[1+2^{l+1}s^l-(2^{l+2}-1)s^{l+1}]&[\bs{x}(s)=e_1],\\
s^{2l+m+n+r}[1+2^{l+1}]&[\bs{x}(s)=2e_1+e_2],\\
3\cdot 2^{l+1}s^{2l+m+n+r}&[\bs{x}(s)=e_1+e_2+e_3].
\end{cases}\lbeq{eta1}\\
&\eta_2[l]:=\sum_{\substack{\bs{u} \\ \bs{u}(t)=2}}D^{*2}(\bs{u})D(\bs{x}-\bs{u})\sum_{\substack{\bs{s} \\ \bs{s}(t)=1}}D(\bs{s})^lD(\bs{u}-\bs{s})D^{*2}(\bs{x}-\bs{s})\nn\\
&\hskip8mm=
\begin{cases}
s^{l+4}[5-2s-8s^2]&[\bs{x}(s)=e_1],\\
8s^{l+6}&[\bs{x}(s)=2e_1+e_2],\\
24s^{l+6}&[\bs{x}(s)=e_1+e_2+e_3],
\end{cases}\lbeq{eta3}\\
&\eta_3[l]:=\sum_{\substack{\bs{s}, \bs{u} \\  \bs{s}(t)=1\\ \bs{u}(t)=2}}D(\bs{s})^2D(\bs{u}-\bs{s})^lD(\bs{x}-\bs{u})^2
=
\begin{cases}
3s^{l+3}[1-s]&[\bs{x}(s)=e_1],\\
3s^{l+4}&[\bs{x}(s)=2e_1+e_2],\\
6s^{l+4}&[\bs{x}(s)=e_1+e_2+e_3].
\end{cases}\lbeq{eta3def}
\end{align}

Therefore, combining \refeq{D3x3}, \refeq{eta1}\refeq{eta3def} and \refeq{eta3} yields
\begin{align}\lbeq{mainpi3final}
p_c^2\sum_{\substack{\bs{u},\bs{v}\\ \bs{x}-\bs{u}\prec \bs{x}-\bs{v}}}D(\bs{x}-\bs{u})D(\bs{x}-\bs{v})\mathcal P_p^{\bs{u},\bs{v}}&=
\begin{cases}
\frac{7p_c^6}{2}s^4-14s^5+O(d^{-6})&[\bs{x}(s)=e_1],\\
s^6+O(d^{-7})&[\bs{x}(s)=2e_1+e_2],\\
9p_c^6s^6-3s^7+O(d^{-8})&[\bs{x}(s)=e_1+e_2+e_3].
\end{cases}
\end{align}

Thus, by \refeq{mathcalHxx3bound}, \refeq{pi0x3} and \refeq{mainpi3final}, we conclude \refeq{prpPi0xxt3}, and together with \refeq{sumofpix1x3decom}, we obtain \refeq{sumPi0x3}.

For \refeq{LLL}, we just apply \refeq{mathcalHxy3bound} and \refeq{Pxyx2pic} to \refeq{lmm:oxoy}. 
\QED

\subsection{The case where the time component equals 4}\label{case4}
In this section, we provides an estimate of $\Pi_p^{(0)}(\bs{x})$ for $\bs{x}(t)=4$, proving the following proposition. 
\begin{prp}\label{prp: Pi0xxt4}For $\bs{x}$ such that $\bs{x}(t)=4$, 
\begin{align}\lbeq{prpPi0xxt4}
\Pi_{p_c}^{\sss(0)}(\bs{x})&=
\begin{cases}
4s^4+O(d^{-5})&[\bs{x}(s)=o],\\
53s^6+O(d^{-7})&[\bs{x}(s)=e_1+e_2],\\
156s^8+O(d^{-9})&[\bs{x}(s)=e_1+e_2+e_3+e_4].
\end{cases}
\end{align}
Then, we have
\begin{align}\lbeq{sumofpix4final}
\sum_{\bs{x}:\bs{x}(t)=4}\Pi_{p_c}^{\sss(0)}(\bs{x})=37s^4+O(d^{-5}). 
\end{align}
\end{prp}
\Proof{\it Proof of Proposition~\ref{prp: Pi0xxt4}.}
First we note that
\begin{align}\lbeq{D4x}
D^{*4}(\bs{x})=
\begin{cases}
3s^2(1-s)&[\bs{x}(s)=o],\\
12s^3(1-2s)&[\bs{x}(s)=e_1+e_2],\\
24s^4&[\bs{x}(s)=e_1+e_2+e_3+e_4].
\end{cases}
\end{align}
and $D^{*4}(\bs{x})=O(d^{-4})$ for $\bs{x}=(2e_1,4), (e_1+e_2+2e_3,4), (2e_1+2e_2,4)$, and $(3e_1+e_2,4)$. Therefore, by applying the BK inequality, we have
\begin{align}\lbeq{sumofpix4}
\sum_{\substack{\bs{x}\\ \bs{x}(t)=4}}\Pi_{p_c}^{\sss(0)}(\bs{x})&=\Pi_{p_c}^{\sss(0)}(o,4)+2d(d-1)\Pi_{p_c}^{\sss(0)}(e_1+e_2,4)\nn\\
&+\frac{2}{3}d(d-1)(d-2)(d-3)\Pi_{p_c}^{\sss(0)}(e_1+e_2+e_3+e_4,4)+O(d^{-5}). 
\end{align}

By \refeq{Hxxbdx4} and \refeq{pi0x3}, \refeq{D4x}, we obtain
\begin{align}\lbeq{Pix0x4}
\Pi_{p_c}^{\sss(0)}(\bs{x})&=\frac{p_c^2}{2}\sum_{\substack{\bs{u}, \bs{v}\\ \bs{u}\neq\bs{v}}}D(\bs{x}-\bs{u})D(\bs{x}-\bs{v})\mathcal P_p^{\bs{u},\bs{v}}+\text{\bf Er}(\bs{x}(s),4)
\end{align}
where
\begin{align}
\text{\bf Er}(\bs{x}(s),4)\le 
\begin{cases}
O(d^{-5})&[\bs{x}(s)=o],\\
O(d^{-7})&[\bs{x}(s)=e_1+e_2],\\
O(d^{-9})&[\bs{x}(s)=e_1+e_2+e_3+e_4].
\end{cases}
\end{align}
By \refeq{LLL} and \refeq{D4x}, we obtain
\begin{align}\lbeq{pix4final}
\Pi_{p_c}^{\sss(0)}(\bs{x})&=\frac{1}{2}\raisebox{-11pt}{\includegraphics[scale=0.75]{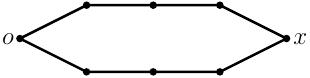}}+\text{\bf Er}(\bs{x}(s),4)
\nn\\
&=\frac{p_c^8}{2}D^{*4}(\bs{x})^2-\underbrace{\raisebox{-15pt}{\includegraphics[scale=0.7]{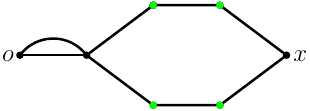}}}_{=p_c^8\xi_1(\bs{x})}+\frac{3}{2}\underbrace{\raisebox{-15pt}{\includegraphics[scale=0.7]{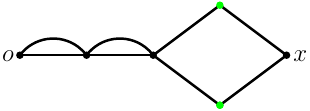}}}_{=p_c^8\xi_2(\bs{x})}\nn\\
&-\frac{1}{2}\underbrace{\raisebox{-15pt}{\includegraphics[scale=0.7]{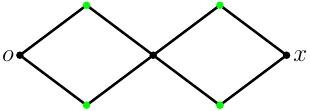}}}_{=p_c^8\xi_3(\bs{x})}-\frac{1}{2}\underbrace{\raisebox{-0pt}{\includegraphics[scale=0.7]{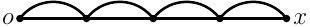}}}_{=p_c^8\xi_4(\bs{x})}+\text{\bf Er}(\bs{x}(s),4)
\end{align}
where
\begin{align}
&\xi_1(\bs{x})
=
\begin{cases}
O(d^{-5}),\\
18s^6+O(d^{-7}),\\
144s^8,
\end{cases}
\xi_2(\bs{x})
=
\begin{cases}
O(d^{-5}),\\
2s^6+O(d^{-7}),\\
48s^8,
\end{cases}
\lbeq{xi12}\\
&
\xi_3(\bs{x})
=
\begin{cases}
s^4+O(d^{-6}),\\
8s^6+O(d^{-7}),\\
96s^8,
\end{cases}
\xi_4(\bs{x})
=
\begin{cases}
O(d^{-6}),\\
O(d^{-7}),\\
24s^8.
\end{cases}
\lbeq{xi34}
\end{align}
Here, the values in each brace in \refeq{xi12}--\refeq{xi34} from top to bottom, correspond to $\bs{x}(s)=o$, for $\bs{x}(s)=e_1+e_2$ and for $\bs{x}(s)=e_1+e_2+e_3+e_4$ respectively. 

Combining \refeq{D4x} and \refeq{pix4final}--\refeq{xi34}, we conclude \refeq{prpPi0xxt4}. Together with \refeq{sumofpix4}, this yields \refeq{sumofpix4final}. \QED

\section{Analysis of $\Pi_p^{(1)}(\bs{x})$}\label{Sectionpi1}
In this section, we provides estimates for $\Pi_{p}^{\sss(1)}(\bs{x})$, recalling the definition: 
\begin{align}
\Pi_{p}^{\sss(1)}(\bs{x})=\sum_{b_1 \in \mathbb B}\mathbb P_p\big(\{\bs{o}\Rightarrow \bs{\underline{b}_1}\}\cap E(b_1,\bs{x}: \tilde C^{b_1}(\bs{o})\big).
\end{align}
We consider the case where $\bs{x}(t)=2$, $\bs{x}(t)=3$ and $\bs{x}(t)=4$ separately, in Section~\ref{Pi12case}, \ref{Pi13case} and \ref{Pi14case}, respectively. For the proof, we use the following lemma which is needed to estimate of the leading term of $\Pi_{p}^{\sss(1)}(\bs{x})$. 
\begin{lmm}\label{lmm: sumPxxosPi0}
\begin{align}\lbeq{sumPxxosPi0}
\sum_{\bs{s}}\mathcal P_{p,[\bs{o},\bs{s}\rangle}^{\bs{x},\bs{x}}=2\Pi_{p}^{\sss(0)}(\bs{x})+
\begin{cases}
\mathcal H_{p}(\bs{x},\bs{x})+\frac{1}{8}F_{9,p}(\bs{x},\bs{x})+\text{\bf Er}(\bs{x}(s),2)&[\bs{x}(t)=2],\\
\mathcal H_{p}(\bs{x},\bs{x})+\text{\bf Er}(\bs{x}(s),3)&[\bs{x}(t)=3],\\
\text{\bf Er}(\bs{x}(s),4)&[\bs{x}(t)=4].
\end{cases}
\end{align}
\end{lmm}
\Proof{\it Proof of Lemma~\ref{lmm: sumPxxosPi0}.}
Since the measure is invariant under the time reversal and spatial translation, we obtain
\begin{align}\lbeq{i}
\mathcal P_{p,[\bs{o},\bs{s}\rangle}^{\bs{x},\bs{x}}
&=\mathbb P_{p}\big(\{ \bs{o}\rightarrow[\bs{x}-\bs{s}, \bs{x} \rangle\} \circ \{\bs{o}\rightarrow \bs{x}\}\big)\nn\\
&=\sum_{\substack{\omega_1: \bs{o}\rightarrow [\bs{x}-\bs{s}, \bs{x} \rangle \\ \omega_2: \bs{o}\rightarrow \bs{x} \\ \omega_1\cap\omega_2=\varnothing }}\mathbb P_{p}\big(\{\omega_1,\omega_2~\text{occupied}\} \cap E_{\prec_{\bs{s}}}(\omega_1) \cap E_{\succ_{\bs{s}}}(\omega_2:\omega_1 )\big)
\end{align}
where we write $E_{\prec_{\bs{s}}}(\omega_1)$ and $E_{\succ_{\bs{s}}}(\omega_2:\omega_1)$ to denote $E_{\prec_{\bs{s}}}(\omega_1)$ and $E_{\succ_{\bs{s}}}(\omega_2:\omega_1)$, respectively in the ordering where $\bs{s}$ is minimal of support of $D$. We denote this ordering by $\prec_{\bs{s}}$ or $\succ_{\bs{s}}$, and define $V_{\prec_{\bs{s}}}^{\bs{o}\rightarrow \bs{x}}(\cdot)$ and $V_{\succ_{\bs{s}}}^{\bs{o}\rightarrow \bs{x}}(\cdot:\cdot)$ in a similar manner. When we  divide $\omega_1$ and $\omega_2$ as $\omega_1=\omega_1'\cup [\bs{x}-\bs{s}, \bs{x}\rangle$ and $\omega_2=\omega'_2\cup [\bs{t},\bs{x}\rangle$ with $\bs{t}\neq \bs{x}-\bs{s}$, similarly as in \refeq{Esucc0}--\refeq{Esucc}, we obtain
\begin{align}
E_{\prec_{\bs{s}}}(\omega'_1\cup[\bs{x}-\bs{s},\bs{x}\rangle)&=E_{\prec_{\bs{s}}}(\omega'_1)\cap V_{\prec_{\bs{s}}}^{\bs{o}\rightarrow \bs{x}}(\bs{s}),\lbeq{ii}\\
E_{\succ_{\bs{s}}}(\omega'_2\cup[\bs{t},\bs{x}\rangle: \omega'_1\cup [\bs{x}-\bs{s},\bs{x}\rangle)&=E_{\succ_{\bs{s}}}(\omega'_2: \omega'_1)\cap V_{\succ_{\bs{s}}}^{\bs{o}\rightarrow \bs{x}}(\bs{x}-\bs{t}: \omega'_1).\lbeq{iii}
\end{align}
We note that when there exists an occupied path from $\bs{o}$ to $\bs{x}-\bs{s}$ and $[\bs{x}-\bs{s},\bs{x}\rangle$ is open, the event $V_{\prec_{\bs{s}}}^{\bs{o}\rightarrow \bs{x}}(\bs{s})$ follows naturally. Therefore, by applying \refeq{ii}--\refeq{iii}, we obtain
\begin{align}
\sum_{\bs{s}}\mathcal P_{p,[\bs{o},\bs{s}\rangle}^{\bs{x},\bs{x}}
&=\sum_{\substack{\bs{s}, \bs{t}\\ \bs{t}\neq\bs{x}-\bs{s}}}D(\bs{s})D(\bs{x}-\bs{t})\mathcal P_{p}^{\bs{x}-\bs{s},\bs{t}}-\mathcal H_{p}^{\bs{s}}(\bs{x},\bs{x})
\end{align}
where, writing 
\begin{align}
Y^{\bs{s}}_{\sss \omega'_1,\omega'_2}&:=\{\omega'_1,\omega'_2~\text{occupied}\}\cap E_{\prec_{\bs{s}}}(\omega'_1)\cap E_{\succ_{\bs{s}}}(\omega'_2:\omega'_1),\\
X_{\sss \omega'_1, \omega'_2}^{\bs{s},\sss [\bs{u},\bs{x}\rangle, [\bs{v},\bs{y}\rangle}&:=\{[\bs{u},\bs{x}\rangle, [\bs{v},\bs{y}\rangle~\text{open}\}\cap Y^{\bs{s}}_{\sss \omega'_1,\omega'_2},
\end{align}
we define 
\begin{align}\lbeq{mathcalHs}
\mathcal H_{p}^{\bs{s}}(\bs{x},\bs{x})=\sum_{\substack{\bs{s}, \bs{t}\\ \bs{t}\neq\bs{x}-\bs{s}}}\sum_{\substack{\omega'_1: \bs{o}\rightarrow \bs{x}-\bs{s} \\ \omega'_2: \bs{o}\rightarrow \bs{t} \\ \omega_1\cap\omega_2=\varnothing }}\mathbb P_{p}\big(X_{\sss \omega'_1, \omega'_2}^{\bs{s},\sss [\bs{x}-\bs{s},\bs{x}\rangle, [\bs{t},\bs{x}\rangle} \cap V_{\succ \bs{s}}^{\sss\bs{o}\rightarrow \bs{x}}(\bs{x}-\bs{t}: \omega_1')^c\big). 
\end{align}
Furthermore, as in \refeq{Voyc}, we obtain
\begin{align}
V_{\succ \bs{s}}^{\sss\bs{o}\rightarrow \bs{x}}(\bs{x}-\bs{t}: \omega_1')^c&=\bigcup_{\substack{\bs{w}\\ \bs{x}-\bs{w}\succ_{\bs{s}} \bs{x}-\bs{t}}}\{\bs{o}\rightarrow [\bs{w},\bs{x}\rangle\}\cap V_{\succ \bs{s}}^{\sss\bs{o}\rightarrow \bs{x}}(\bs{x}-\bs{w}: \omega_1').\lbeq{Voycs}
\end{align}
Thus, by repeatedly applying \refeq{Voycs} to \refeq{mathcalHs}, as in Section~\ref{time2and3}, $\mathcal H_{p}^{\bs{s}}(\bs{x},\bs{x})$ can be reformulated as 
\begin{align}
\mathcal H_{p}^{\bs{s}}(\bs{x},\bs{x})&=
\begin{cases}
\frac{1}{4}F_{2,p}(\bs{x},\bs{x})-\frac{1}{8}F_{4,p}(\bs{x},\bs{x})+\text{\bf Er}(\bs{x}(s),2)&[\bs{x}(t)=2]\\
\frac{1}{4}F_{2,p}(\bs{x},\bs{x})+\text{\bf Er}(\bs{x}(s),3)&[\bs{x}(t)=3]\\
\text{\bf Er}(\bs{x}(s),4)&[\bs{x}(t)=4]
\end{cases}
\nn\\
&=
\begin{cases}
\mathcal H_{p}(\bs{x},\bs{x})+\frac{1}{8}F_{9,p}(\bs{x},\bs{x})+\text{\bf Er}(\bs{x}(s),2)&[\bs{x}(t)=2],\\
\mathcal H_{p}(\bs{x},\bs{x})+\text{\bf Er}(\bs{x}(s),3)&[\bs{x}(t)=3],\\
\text{\bf Er}(\bs{x}(s),4)&[\bs{x}(t)=4].
\end{cases}
\end{align} 
For the case where $\bs{x}(t)=4$, we estimate it similary in the proof of \refeq{Hxxbdx4}. Therefore, comparing $\sum_{\bs{s}}\mathcal P_{p_c,[\bs{o},\bs{s}\rangle}^{\bs{x},\bs{x}}$ with $\Pi_{p_c}^{\sss(0)}(\bs{x})$ given in a form of \refeq{pi0x3}, we obtain \refeq{sumPxxosPi0}. \QED

\subsection{The case where the time component equals 2}\label{Pi12case}
This section provides an estimate of $\Pi_{p}^{\sss(1)}(\bs{x})$ for $\bs{x}(t)=2$, proving the following proposition.
\begin{prp}\label{prp: Pix1xt2}
For $\bs{x}$ such that $\bs{x}(t)=2$, 
\begin{align}\lbeq{prpPix1xt2}
\Pi_{p_c}^{\sss(1)}(\bs{x})=
\begin{cases}
p_c^4s^2-\frac{p_c^4}{2}[2+p_c^2]s^3+2s^4+O(d^{-5})&[\bs{x}(s)=o]\\
2p_c^4s^4+O(d^{-8})&[\bs{x}(s)=e_1+e_2].
\end{cases}
\end{align}
Then, we obtain
\begin{align}\lbeq{sumofpi1x2final}
\sum_{\bs{x}:\bs{x}(t)=2}\Pi_{p_c}^{\sss(1)}(\bs{x})=2p_c^4s^2-\frac{p_c^4}{2}[6+p_c^2]s^3+2s^4+O(d^{-5}).
\end{align}
\end{prp}

\begin{lmm}\label{lmm: Pix1xt2}
For $\bs{x}$ such that $\bs{x}(t)=2$, 
\begin{align}\lbeq{h}
\Pi_{p_c}^{\sss(1)}(\bs{x})=\sum_{\bs{s}}\mathcal P_{p_c,[\bs{o},\bs{s}\rangle}^{\bs{x},\bs{x}}. 
\end{align}
\end{lmm}
\Proof{\it Proof of Lemma~\ref{lmm: Pix1xt2}.}
When $\bs{x}(t)=2$, the only possible choice is $\bs{\underline{b}_1}=\bs{o}$, since if $\bs{\underline{b}_1}(t)=1$, $\bs{o}$ cannot be doubly connected. Moreover, the bond $[\bs{\bar{b}_1}, \bs{x}\rangle$ is the unique pivotal bond for the event $\bs{\bar{b}_1} \rightarrow \bs{x}$. Accordingly, in the case where $\bs{x}(t)=2$, we obtain \refeq{h}. 
\QED

\Proof{\it Proof of Proposition~\ref{prp: Pix1xt2}.}
Combining Proposition~\ref{prp: Pi0xxt2}, \refeq{mathcalHxx2bound}, \refeq{F9est} and Lemma~\ref{lmm: sumPxxosPi0} yields \refeq{prpPix1xt2}. Together with \refeq{sumofpi0x2} where the summand is replaced by $\Pi_{p_c}^{\sss(0)}(\bs{x})$, we obtain \refeq{sumofpi1x2final}. 
\QED

\subsection{The case where the time component equals 3}\label{Pi13case}
In this section, we estimate $\Pi_{p}^{\sss(1)}(\bs{x})$ for $\bs{x}(t)=3$, establishing the following proposition.
\begin{prp}\label{prp: Pix1xt3}
For $\bs{x}$ such that $\bs{x}(t)=3$, 
\begin{align}\lbeq{prpPix1xt3}
\Pi_{p_c}^{\sss(1)}(\bs{x})=
\begin{cases}
7p_c^6s^4-31s^5+O(d^{-6})&[\bs{x}(s)=e_1],\\
2s^6+O(d^{-7})&[\bs{x}(s)=2e_1+e_2],\\
18p_c^6s^6-6s^7+O(d^{-8})&[\bs{x}(s)=e_1+e_2+e_3].
\end{cases}
\end{align}
Then, we have
\begin{align}\lbeq{prpsumPix1xt3}
\sum_{\bs{x}:\bs{x}(t)=3}\Pi_{p_c}^{\sss(1)}(\bs{x})=10p_c^6s^3-49s^4+O(d^{-5}).
\end{align}
\end{prp}

\begin{lmm}\label{lmm: Pix1xt3}.
For $\bs{x}$ such that $\bs{x}(t)=3$, 
\begin{align}\lbeq{lmmPix1xt3}
\Pi_{p_c}^{\sss(1)}(\bs{x})&=\sum_{\bs{s}}\mathcal P_{p_c,[\bs{o},\bs{s}\rangle}^{\bs{x},\bs{x}}+\sum_{\substack{\bs{u}, \bs{v} \\ \bullet}}D(\bs{x}-\bs{u})D(\bs{x}-\bs{v})\mathcal P_p^{\bs{u},\bs{u}}(\bs{v})\nn\\
&-\sum_{\bs{s}, \bs{v}}\mathbb P_{p_c}\big(\big\{[\bs{o}, \bs{s}\rangle[\bs{s}, \bs{v}\rangle[\bs{v},\bs{x}\rangle~\text{open} \}\circ \{\bs{o}\rightarrow \bs{x}\}\big\} \cap\{\bs{v} \in \tilde C^{[\bs{o}, \bs{s}\rangle}(\bs{o})\}\big)+\text{\bf Er}(\bs{x}(s),3).
\end{align}
\end{lmm}

\Proof{\it Proof of Lemma~\ref{lmm: Pix1xt3}.}
By definition, $\Pi_{p_c}^{\sss(1)}(\bs{x})$ can be decomposed according to the time-component of $\bs{\underline{b}_1}$ as
\begin{align}\lbeq{pi1x3}
\Pi_{p_c}^{\sss(1)}(\bs{x})=\sum_{\bs{s}}\mathbb P_{p_c}\big(E([\bs{o}, \bs{s}\rangle,\bs{x}: \tilde C^{[\bs{o}, \bs{s}\rangle}(\bs{o}))\big)+\sum_{\bs{u}}\mathbb P_p\big(\{\bs{o}\Rightarrow \bs{u}\}\cap E([\bs{u},\bs{x}\rangle ,\bs{x}: \tilde C^{[\bs{u},\bs{x}\rangle}(\bs{o}))\big).
\end{align}
Since we have $E([\bs{o}, \bs{s}\rangle,\bs{x}: \tilde C^{[\bs{o}, \bs{s}\rangle}(\bs{o})) \subset  \{[\bs{o},\bs{s}\rangle \rightarrow \bs{x} \}\circ \{\bs{o}\rightarrow \bs{x}\}$, we can decompose the first term of \refeq{pi1x3} as 
\begin{align}\lbeq{PEosxCosodecom}
&\sum_{\bs{s}}\mathbb P_{p_c}\big(E([\bs{o}, \bs{s}\rangle,\bs{x}: \tilde C^{[\bs{o}, \bs{s}\rangle}(\bs{o})\big)\nn\\
&=\sum_{\bs{s}}\mathcal P_{p_c,[\bs{o},\bs{s}\rangle}^{\bs{x},\bs{x}}-\sum_{\bs{s}}\mathbb P_{p_c}\big(E([\bs{o}, \bs{s}\rangle,\bs{x}: \tilde C^{[\bs{o}, \bs{s}\rangle}(\bs{o}))^c\cap \big\{\{[\bs{o},\bs{s}\rangle \rightarrow \bs{x} \}\circ \{\bs{o}\rightarrow \bs{x}\}\big\} \big),
\end{align}
where the first term corresponds to the first term of \refeq{lmmPix1xt3}. Thus, we establsih 
\begin{align}\lbeq{fff}
&\sum_{\bs{s}}\mathbb P_{p_c}\big(\big\{\{[\bs{o},\bs{s}\rangle \rightarrow \bs{x} \} \circ \{\bs{o}\rightarrow \bs{x}\}\big\}\cap  E\big([\bs{o}, \bs{s}\rangle,\bs{x}: \tilde C^{[\bs{o}, \bs{s}\rangle}(\bs{o})\big)^c\big)\nn\\
&=\sum_{\bs{s}, \bs{v}}\mathbb P_{p_c}\big(\big\{[\bs{o}, \bs{s}\rangle[\bs{s}, \bs{v}\rangle[\bs{v},\bs{x}\rangle~\text{open} \}\circ \{\bs{o}\rightarrow \bs{x}\}\big\} \cap\{\bs{v} \in \tilde C^{[\bs{o}, \bs{s}\rangle}(\bs{o})\}\big)+\text{\bf Er}(\bs{x}(s),3).
\end{align}
We note that 
\begin{align}\lbeq{j}
&E\big([\bs{o}, \bs{s}\rangle,\bs{x}: \tilde C^{[\bs{o}, \bs{s}\rangle}(\bs{o})\big)^c\nn\\
&=\{[\bs{o},\bs{s}\rangle \rightarrow \bs{x} \}^c\cup \{\bs{x} \notin \tilde C^{[\bs{o}, \bs{s}\rangle}(\bs{o}) \}\cup \{\exists\bs{b'} \in \text{piv}(\bs{s}\rightarrow \bs{x}),  \bs{\underline{b}'} \in \tilde C^{[\bs{o}, \bs{s}\rangle}(\bs{o}) \}
\end{align}
and, the events $\{[\bs{o},\bs{s}\rangle \rightarrow \bs{x} \}^c$ and $\{\bs{x} \notin \tilde C^{[\bs{o}, \bs{s}\rangle}(\bs{o}) \}$ are disjoint from the event $\{[\bs{o},\bs{s}\rangle \rightarrow \bs{x} \} \circ \{\bs{o}\rightarrow \bs{x}\}$. Therefore, it suffices to consider the intersection of the event $\{[\bs{o},\bs{s}\rangle \rightarrow \bs{x} \} \circ \{\bs{o}\rightarrow \bs{x}\}$ with the event $\{\exists\bs{b'} \in \text{piv}(\bs{s}\rightarrow \bs{x}),  \bs{\underline{b}'} \in \tilde C^{[\bs{o}, \bs{s}\rangle}(\bs{o}) \}$.

Since $\bs{x}(t)-\bs{s}(t)=2$, we have
\begin{align} 
\{\exists\bs{b'} \in \text{piv}(\bs{s}\rightarrow \bs{x}),  \bs{\underline{b}'} \in \tilde C^{[\bs{o}, \bs{s}\rangle}(\bs{o}) \}=\bigcup_{\bs{v}} \{[\bs{s},\bs{v}\rangle \in \text{piv}(\bs{s}\rightarrow \bs{x})\} \cap \{\bs{v} \in \tilde C^{[\bs{o}, \bs{s}\rangle}(\bs{o})\}.
\end{align} 
We note that when $[\bs{s}, \bs{v}\rangle$ is pivotal for the event $\bs{s}\rightarrow \bs{x}$, $[\bs{v},\bs{x}\rangle$ is also pivotal for the event. Since by applying the BK inequality and writing $[\bs{o}, \bs{s}\rangle[\bs{s}, \bs{v}\rangle[\bs{v},\bs{x}\rangle$ by $\omega_{\bs{s},\bs{v}}$, we bound the contribution from the event $\{[\bs{s},\bs{v}\rangle \in \text{piv}(\bs{s}\rightarrow \bs{x})\}^c$ as 
\begin{align}
&\sum_{\bs{s}, \bs{v}}\mathbb P_{p_c}\big(\big\{\{\omega_{\bs{s},\bs{v}}~\text{occupied} \}\circ \{\bs{o}\rightarrow \bs{x}\}\big\} \cap  \{[\bs{s},\bs{v}\rangle \in \text{piv}(\bs{s}\rightarrow \bs{x})\}^c \cap \{\bs{v} \in \tilde C^{[\bs{o}, \bs{s}\rangle}(\bs{o})\}\big)\nn\\
&\le \sum_{\bs{s}, \bs{v}}\sum_{ \substack{\bs{r}, \bs{w} \\ \bs{r}(t)=0,1, \bs{r}\neq\bs{s}\\ \bs{w}(t)=2,3, \bs{w}\neq \bs{v}}}\mathbb P_{p_c}\big(\{\omega_{\bs{s},\bs{v}}~\text{occupied} \}\circ \{\bs{o}\rightarrow \bs{r} \rightarrow \bs{w} \rightarrow \bs{x}\} \circ \{\bs{r} \rightarrow \bs{v}\} \circ \{ \bs{s}\rightarrow \bs{w}\}\big)\nn\\
&\le Cp_c^8s^2D^{*3}(\bs{x})^2\nn\\
&\le \text{\bf Er}(\bs{x}(s),3),
\end{align}
we obtain \refeq{fff}. 

In the remainder of the proof, we consider the second term of \refeq{pi1x3}. Namely, we aim to show 
\begin{align}\lbeq{sumPodbuE}
&\sum_{\bs{u}}\mathbb P_p\big(\{\bs{o}\Rightarrow \bs{u}\}\cap E([\bs{u},\bs{x}\rangle ,\bs{x}: \tilde C^{[\bs{u},\bs{x}\rangle}(\bs{o}))\big)\nn\\
&=\sum_{\substack{\bs{u}, \bs{v} \\ \bs{u}\neq \bs{v}}}D(\bs{x}-\bs{u})D(\bs{x}-\bs{v})\mathcal P_p^{\bs{u},\bs{u}}(\bs{v})+\text{\bf Er}^{\sss(1)}(\bs{x}(s),3).
\end{align}

First we note that 
\begin{align}
E\big([\bs{u},\bs{x}\rangle ,\bs{x}: \tilde C^{[\bs{u},\bs{x}\rangle}(\bs{o})\big)&=\{[\bs{u},\bs{x}\rangle~\text{open}\}\cap \{\bs{x}\in \tilde C^{[\bs{u}, \bs{x}\rangle}(\bs{o})\}\nn\\
&=\{\bs{o}\rightarrow \bs{x}\} \circ \{[\bs{u},\bs{x}\rangle~\text{open}\}\nn\\
&=\bigcup_{\substack{\omega:\bs{o}\rightarrow \bs{x} \\ \omega \cap [\bs{u},\bs{x} \rangle =\varnothing}}\{\omega~\text{occupied}\}\cap E_{\succ}(\omega:[\bs{u}, \bs{x}\rangle)\cap \{[\bs{u},\bs{x}\rangle~\text{open}\}.
\end{align}
Therefore, the left-han side of \refeq{sumPodbuE} can be rewritten by
\begin{align}
&\sum_{\substack{\omega:\bs{o}\rightarrow \bs{x} \\ \omega \cap [\bs{u},\bs{x} \rangle =\varnothing}}\mathbb P_{p_c}\Big(\{\bs{o} \Rightarrow \bs{u} \}\cap \{\omega~\text{occupied}\}\cap E_{\succ}(\omega:[\bs{u}, \bs{x}\rangle)\cap \{[\bs{u},\bs{x}\rangle~\text{open}\}\Big)\nn\\
&=\sum_{\substack{\bs{u}, \bs{v} \\ \bs{u}\neq \bs{v}}}\sum_{\omega: \bs{o}\rightarrow \bs{v}}\mathbb P_{p_c}\Big(\{\bs{o} \Rightarrow \bs{u} \}\cap \{\omega~\text{occupied}\}\cap E_{\succ}(\omega \cup [\bs{v}, \bs{x}\rangle:[\bs{u}, \bs{x}\rangle)\cap \big\{\substack{[\bs{u},\bs{x}\rangle \\ [\bs{v},\bs{x}\rangle}~\text{open}\big\}\Big). 
\end{align}
By \refeq{Esucc}, we have $E_{\succ}(\omega \cup [\bs{v}, \bs{x}\rangle:[\bs{u}, \bs{x}\rangle)=E_{\succ}(\omega)\cap  V_{\succ}^{\bs{o}\rightarrow \bs{x}}(\bs{v}-\bs{x}: [\bs{u}, \bs{x}\rangle)$. If we disregrad $V_{\succ}^{\bs{o}\rightarrow \bs{x}}(\bs{v}-\bs{x}: [\bs{u}, \bs{x}\rangle)$, we obtain the leading term of \refeq{sumPodbuE}. By applying \refeq{Voyc}, the contribution from $V_{\succ}^{\bs{o}\rightarrow \bs{x}}(\bs{v}-\bs{x}: [\bs{u}, \bs{x}\rangle)^c$ can be estimated as
\begin{align}
&\sum_{\substack{\bs{u}, \bs{v} \\ \bs{u}\neq \bs{v}}}\sum_{\omega: \bs{o}\rightarrow \bs{v}}\mathbb P_{p_c}\Big(\{\bs{o} \Rightarrow \bs{u} \}\cap \{\omega~\text{occupied}\}\cap E_{\succ}(\omega)\cap  V_{\succ}^{\bs{o}\rightarrow \bs{x}}(\bs{v}-\bs{x}: [\bs{u}, \bs{x}\rangle)^c \cap \big\{\substack{[\bs{u},\bs{x}\rangle \\ [\bs{v},\bs{x}\rangle }~\text{open}\big\}\Big)\nn\\
&\le p_c^3\sum_{\substack{\bs{u}, \bs{v}, \bs{w} \\ \bs{u}\neq \bs{v}, \bs{w} \neq \bs{u} \\ \bs{w}-\bs{x}\succ \bs{v}-\bs{x}}}D(\bs{x}-\bs{u})D(\bs{x}-\bs{v})D(\bs{x}-\bs{w})\mathcal P_p^{\bs{u},\bs{u}}(\bs{v},\bs{w})\nn\\
&\le Cs^2D^{*3}(\bs{x})^2\nn\\
&\le \text{\bf Er}(\bs{x}(s),3)
\end{align}
which contributes to the error term of \refeq{sumPodbuE}. 
\QED

\Proof{\it Proof of Proposition~\ref{prp: Pix1xt3}.}For the first term of \refeq{lmmPix1xt3}, combining Lemma~\ref{lmm: sumPxxosPi0} with Proposition~\ref{prp: Pi0xxt3} and \refeq{mathcalHxx3bound}, we obtain
\begin{align}
\sum_{\bs{s}}\mathcal P_{p_c,[\bs{o},\bs{s}\rangle}^{\bs{x},\bs{x}}=
\begin{cases}
7p_c^6s^4-30s^5+O(d^{-6})&[\bs{x}(s)=e_1],\\
2s^6+O(d^{-7})&[\bs{x}(s)=2e_1+e_2],\\
18p_c^6s^6-6s^7+O(d^{-8})&[\bs{x}(s)=e_1+e_2+e_3].
\end{cases}
\end{align}

Next we estimate the second term of \refeq{lmmPix1xt3}, observing that combining \refeq{ouovowestimate}--\refeq{ouovowestimate1} and \refeq{osxoy2} provides
\begin{align}
\Big|\mathcal P_p^{\bs{u},\bs{u}}(\bs{w})-\raisebox{-12pt}{\includegraphics[scale=0.75]{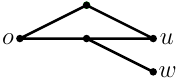}}\Big|&\le 
CD^{*2}(\bs{u})^2D^{*2}(\bs{w}),\lbeq{Puumainerror1}\\
\Big|\mathcal P_p^{\bs{u},\bs{u}}(\bs{w})-\Big[\raisebox{-12pt}{\includegraphics[scale=0.75]{zetauuwmain}}+\frac{1}{2}\raisebox{-12pt}{\includegraphics[scale=0.75]{Ppuuw}}~\Big]\Big|&\le 
CsD^{*2}(\bs{u})^2D^{*2}(\bs{w}).\lbeq{Puumainerror2}
\end{align}
We establish
\begin{align}\lbeq{dd}
\sum_{\substack{\bs{u}, \bs{v} \\ \bullet }}D(\bs{x}-\bs{u})D(\bs{x}-\bs{v})\mathcal P_p^{\bs{u},\bs{u}}(\bs{v})=
\begin{cases}
4s^5+O(d^{-6})&[\bs{x}(s)=e_1],\\
O(d^{-7})&[\bs{x}(s)=2e_1+e_2],\\
6s^7+O(d^{-9})&[\bs{x}(s)=e_1+e_2+e_3].
\end{cases}
\end{align}
In the case where $\bs{x}(s)=e_1$, by \refeq{Puumainerror2}, we obtain
\begin{align}\lbeq{gg}
&\Big|\sum_{\substack{\bs{u}, \bs{v} \\ \bullet }}D(\bs{x}-\bs{u})D(\bs{x}-\bs{v})\mathcal P_p^{\bs{u},\bs{u}}(\bs{v})-\Big[\raisebox{-11pt}{\includegraphics[scale=0.75]{M2}}+\frac{1}{2}\raisebox{-11pt}{\includegraphics[scale=0.75]{F2main}}\Big]\Big|\nn\\
&\le Cs^2D^{*2}(\bs{x})^2\nn\\
&\le \text{\bf Er}^{\sss(1)}(\bs{x}(s),3).
\end{align}
For the case $\bs{x}(s)=2e_1+e_2$ and $\bs{x}(s)=e_1+e_2+e_3$, by \refeq{Puumainerror1} and the bound $\sup_{\bs{x}\neq\bs{o}}D^{*2}(\bs{x})\le O(d^{-2})$, we obtain
\begin{align}\lbeq{ggg}
\Big|\sum_{\substack{\bs{u}, \bs{v} \\ \bullet }}D(\bs{x}-\bs{u})D(\bs{x}-\bs{v})\mathcal P_p^{\bs{u},\bs{u}}(\bs{v})-\raisebox{-11pt}{\includegraphics[scale=0.75]{M2}}\Big|&\le Cs^2D^{*2}(\bs{x})^2\nn\\
&\le \text{\bf Er}(\bs{x}(s),3).
\end{align}
The diagrams appearing in \refeq{gg}--\refeq{ggg} can be estimated as in \refeq{ee} and as in \refeq{mainF2xxx3}. Thus, we conclude \refeq{dd}. 

Finally, we consider the third term of \refeq{lmmPix1xt3}, aiming to show
\begin{align}\lbeq{ffff}
&\sum_{\bs{s}, \bs{v}}\mathbb P_{p_c}\big(\big\{[\bs{o}, \bs{s}\rangle[\bs{s}, \bs{v}\rangle[\bs{v},\bs{x}\rangle~\text{open} \}\circ \{\bs{o}\rightarrow \bs{x}\}\big\} \cap\{\bs{v} \in \tilde C^{[\bs{o}, \bs{s}\rangle}(\bs{o})\}\big)\nn\\
&=
\begin{cases}
5s^5+O(d^{-6})&[\bs{x}(s)=e_1],\\
O(d^{-7})&[\bs{x}(s)=2e_1+e_2],\\
6s^7+O(d^{-9})&[\bs{x}(s)=e_1+e_2+e_3].
\end{cases}
\end{align}
Since we obtain
\begin{align}
&\big\{\{\omega_{\bs{s},\bs{v}}~\text{occupied} \}\circ \{\bs{o}\rightarrow \bs{x}\}\big\} \cap\{\bs{v} \in \tilde C^{[\bs{o}, \bs{s}\rangle}(\bs{o})\}\nn\\
&=\bigcup_{ \substack{\omega_1: \bs{o}\rightarrow \bs{x} \\ \omega_1\cap \omega_{\bs{s},\bs{v}}=\varnothing }}\bigcup_{ \substack{\omega_2: \bs{o}\rightarrow \bs{v} \\ \omega_2\cap [\bs{o},\bs{s} \rangle=\varnothing }}\{\omega_{\bs{s},\bs{v}}, \omega_1,\omega_2~\text{occupied} \}\cap E_{\succ}(\omega_1:\omega_{\bs{s},\bs{v}})\cap E_{\succ}(\omega_2:[\bs{o},\bs{s}\rangle),
\end{align}
and the contributions from $E_{\succ}(\omega_1:\omega_{\bs{s},\bs{v}})^c$ and $E_{\succ}(\omega_2:[\bs{o},\bs{s}\rangle)^c$ are accounted for $\text{\bf Er}(\bs{x}(s),3)$, it suffices to estimate 
\begin{align}
\sum_{\bs{s}, \bs{v}}\sum_{ \substack{\omega_1: \bs{o}\rightarrow \bs{x} \\ \omega_1\cap \omega_{\bs{s},\bs{v}}=\varnothing }}\sum_{ \substack{\omega_2: \bs{o}\rightarrow \bs{v} \\ \omega_2\cap [\bs{o},\bs{s} \rangle=\varnothing }}\mathbb P_{p_c}\big(\{\omega_{\bs{s},\bs{v}}, \omega_1,\omega_2~\text{occupied} \}\big)=\raisebox{-11pt}{\includegraphics[scale=0.75]{F2main}}+\raisebox{-11pt}{\includegraphics[scale=0.75]{M2}}
\end{align}
where we distinguish the two cases: $\omega_1\cap \omega_2=\varnothing$ and $\omega_1\cap \omega_2\neq \varnothing$. The first and second diagrams can be estimated as in \refeq{mainF2xxx3} and \refeq{ee}, respectively. This yields \refeq{ffff}, which together with \refeq{dd}, gives \refeq{prpPix1xt3}. Then, \refeq{prpsumPix1xt3} follows from \refeq{prpPix1xt3} and \refeq{sumofpix1x3decom} with $\Pi_{p_c}^{\sss (0)}$ replaced by  $\Pi_{p_c}^{\sss (1)}$. 
\QED

\subsection{The case where the time component equals 4}\label{Pi14case}
In this section, we estimate $\Pi_{p}^{\sss(1)}(\bs{x})$ for $\bs{x}(t)=4$, proving the following proposition. 
\begin{prp}\label{prp: Pix1xt4}
For $\bs{x}$ such that $\bs{x}(t)=4$,
\begin{align}\lbeq{prpPix1xt4}
\Pi_{p_c}^{\sss(1)}(\bs{x})=
\begin{cases}
\frac{17}{2}s^4+O(d^{-5})&[\bs{x}(s)=o],\\
108s^6+O(d^{-7})&[\bs{x}(s)=e_1+e_2],\\
324s^8+O(d^{-9})&[\bs{x}(s)=e_1+e_2+e_3+e_4].
\end{cases}
\end{align}
Then, we have
\begin{align}\lbeq{sumofpi4x1}
\sum_{\bs{x}:\bs{x}(t)=4}\Pi_{p_c}^{\sss(1)}(\bs{x})=76s^4+O(d^{-5}).
\end{align}
\end{prp}

\begin{lmm}\label{lmm: Pix1xt4}
For $\bs{x}$ such that $\bs{x}(t)=4$,
\begin{align}\lbeq{lmmeqPix1xt4}
\Pi_{p_c}^{\sss(1)}(\bs{x})&=\sum_{\bs{s}}\mathcal P_{p_c,[\bs{o},\bs{s}\rangle}^{\bs{x},\bs{x}}+\sum_{\substack{[\bs{u}, \bs{v}\rangle \\ \bs{u}(t)=2}}\mathcal P_{p_c}^{\bs{u},\bs{u}}\mathcal P_{p_c,[\bs{u},\bs{v}\rangle}^{\bs{x},\bs{x}}+\text{\bf Er}(\bs{x}(s),4).
\end{align}
\end{lmm}

\Proof{\it Proof of Lemma~\ref{lmm: Pix1xt4}. }
Similarly, as in \refeq{pi1x3}, we have
\begin{align}\lbeq{pi1x4}
\Pi_{p_c}^{\sss(1)}(\bs{x})&=\sum_{\bs{s}}\mathbb P_{p_c}\big(E([\bs{o}, \bs{s}\rangle,\bs{x}: \tilde C^{[\bs{o}, \bs{s}\rangle}(\bs{o})\big)+\sum_{\substack{b_1 \in \mathbb B \\ \bs{\underline{b}_1}(t)=2,3}}\mathbb P_{p_c}\big(\{\bs{o}\Rightarrow \bs{\underline{b}_1}\}\cap E(b_1,\bs{x}: \tilde C^{b_1}(\bs{o})\big). 
\end{align}
We analyze the first and the second term separately. We aim to establish 
\begin{align}
\sum_{\bs{s}}\mathbb P_{p_c}\big(E([\bs{o}, \bs{s}\rangle,\bs{x}: \tilde C^{[\bs{o}, \bs{s}\rangle}(\bs{o})\big)&=\sum_{\bs{s}}\mathcal P_{p_c,[\bs{o},\bs{s}\rangle}^{\bs{x},\bs{x}}+\text{\bf Er}(\bs{x}(s),4),\lbeq{pi1x4first}\\
\sum_{\substack{b_1 \in \mathbb B \\ \bs{\underline{b}_1}(t)=2,3}}\mathbb P_{p_c}\big(\{\bs{o}\Rightarrow \bs{\underline{b}_1}\}\cap E(b_1,\bs{x}: \tilde C^{b_1}(\bs{o})\big)&=\sum_{\substack{[\bs{u}, \bs{v}\rangle \\ \bs{u}(t)=2}}\mathcal P_{p_c}^{\bs{u},\bs{u}}\mathcal P_{p_c,[\bs{u},\bs{v}\rangle}^{\bs{x},\bs{x}}+
\text{\bf Er}(\bs{x}(s),4)\lbeq{odbb1Eb1xC}
\end{align}
which together with \refeq{pi1x4}, yields \refeq{lmmeqPix1xt4}. 

We begin with \refeq{pi1x4first}. As in \refeq{PEosxCosodecom}, it suffices to show
\begin{align}
\sum_{\bs{s}}\mathbb P_{p_c}\big( \big\{\{[\bs{o},\bs{s}\rangle \rightarrow \bs{x} \}\circ \{\bs{o}\rightarrow \bs{x}\}\big\} \cap E\big([\bs{o}, \bs{s}\rangle,\bs{x}: \tilde C^{[\bs{o}, \bs{s}\rangle}(\bs{o})\big)^c\big)\le \text{\bf Er}(\bs{x}(s),4). 
\end{align}
By \refeq{j} and the reasoning given below it, we consider the intersection of the event $\{[\bs{o},\bs{s}\rangle \rightarrow \bs{x} \}\circ \{\bs{o}\rightarrow \bs{x}\}$ with the event $\{\exists\bs{b'} \in \text{piv}(\bs{s}\rightarrow \bs{x}),  \bs{\underline{b}'} \in \tilde C^{[\bs{o}, \bs{s}\rangle}(\bs{o}) \}$. By applying the BK inequality, we obtain
\begin{align}
&\sum_{\bs{s}}\mathbb P_{p_c}\big( \big\{\{[\bs{o},\bs{s}\rangle \rightarrow \bs{x} \}\circ \{\bs{o}\rightarrow \bs{x}\}\big\} \cap \{\exists\bs{b'} \in \text{piv}(\bs{s}\rightarrow \bs{x}),  \bs{\underline{b}'} \in \tilde C^{[\bs{o}, \bs{s}\rangle}(\bs{o}) \}\big)\nn\\
&\le \sum_{\substack{\bs{s},\bs{u}, \bs{r} \\ \bs{u}(t)=2,3 \\ \bs{r}\neq \bs{o} }}\sum_{\substack{\bs{w} \\ \bs{w}=\bs{r}, \bs{o}}}\mathbb P_{p_c}\big(\{[\bs{o},\bs{s}\rangle \rightarrow \bs{u} \rightarrow \bs{x} \}\circ \{\bs{o}\rightarrow \bs{r} \rightarrow \bs{u}\}\circ \{\bs{w} \rightarrow \bs{x} \}\big)\nn\\
&\le Cs(1+s)D^{*4}(\bs{x})^2\nn\\
&=\text{\bf Er}(\bs{x}(s),4),
\end{align}
from which follows \refeq{pi1x4first}.

Finally, we consider \refeq{odbb1Eb1xC}. In the case where $\bs{\underline{b}_1}(t)=2$, the bond $[\bs{\bar{b}_1}, \bs{x}\rangle$ is the unique pivotal bond for the event that $\bs{\bar{b}_1}\rightarrow \bs{x}$. Thus, by applying the BK inequality, we obtain
\begin{align}\lbeq{ob1E}
&\sum_{\substack{b_1 \in \mathbb B \\ \bs{\underline{b}_1}(t)=2}}\mathbb P_{p_c}\big(\{\bs{o}\Rightarrow \bs{\underline{b}_1}\}\cap E(b_1,\bs{x}: \tilde C^{b_1}(\bs{o})\big)-\sum_{\substack{[\bs{u}, \bs{v}\rangle \\ \bs{u}(t)=2}}\mathcal P_{p_c}^{\bs{u},\bs{u}}\mathcal P_{p_c,[\bs{u},\bs{v}\rangle}^{\bs{x},\bs{x}}\nn\\
&=\sum_{\substack{[\bs{u}, \bs{v}\rangle \\ \bs{u}(t)=2}}\mathbb P_{p_c}\big(\{\bs{o}\Rightarrow \bs{u}\}\cap E(b_1,\bs{x}: \tilde C^{b_1}(\bs{o}))\cap  \big\{\{[\bs{u}, \bs{v}\rangle \rightarrow \bs{x}\}\circ \{\bs{o}\rightarrow \bs{x}\}\big\}^c\big)\nn\\
&\le\sum_{\substack{[\bs{u}, \bs{v}\rangle \\ \bs{u}(t)=2}}\sum_{\substack{\bs{w} \\ \bs{w}\neq \bs{u}}} \mathbb P_{p_c}\big(\{\bs{o}\rightarrow \bs{u}\}\circ \{\bs{o}\rightarrow \bs{w}\rightarrow [\bs{u}, \bs{v}\rangle \rightarrow \bs{x}\}\circ \{\bs{w}\rightarrow \bs{x}\}\big)\nn\\
&\le p_c^2\sum_{\substack{[\bs{u}, \bs{v}\rangle \\ \bs{u}(t)=2}}\sum_{\substack{\bs{w} \\ \bs{w}\neq \bs{u}}} \tau_{p_c}(\bs{u})\tau_{p_c}(\bs{w})\tau_{p_c}(\bs{u}-\bs{w})D(\bs{v}-\bs{u})D(\bs{x}-\bs{v})\tau_{p_c}(\bs{x}-\bs{w})\nn\\
&\le Cp_c^9sD^{*4}(\bs{x})^2\nn\\
&=\text{\bf Er}(\bs{x}(s),4)
\end{align}

When $\bs{\underline{b}_1}(t)=3$, we have 
$E(b_1,\bs{x}: \tilde C^{b_1}(\bs{o}))=\{[\bs{\underline{b}_1},\bs{x}\rangle~\text{open}\}\cap \{\bs{x} \in \tilde C^{b_1}(\bs{o})\}$, since $\bs{\bar{b}_1}=\bs{x}$ and $\text{piv}(\bs{x}\rightarrow \bs{x})=\varnothing$. Therefore, as in \refeq{ob1E}, applying the BK inequality, we obtain 
\begin{align}
&\sum_{\substack{b_1 \in \mathbb B \\ \bs{\underline{b}_1}(t)=3}}\mathbb P_{p_c}\big(\{\bs{o}\Rightarrow \bs{\underline{b}_1}\}\cap E(b_1,\bs{x}: \tilde C^{b_1}(\bs{o}))\big)\nn\\
&\le \sum_{\substack{\bs{u}\\ \bs{u}(t)=3}}\sum_{\substack{\bs{w}\\ \bs{w}\neq\bs{u}}}\sum_{\substack{\bs{z}\\ \bs{z}(t)\le 2}}\mathbb P_{p_c}\big(\{\bs{o}\rightarrow\bs{z}\rightarrow  [\bs{u},\bs{x}\rangle \}\circ \{\bs{o}\rightarrow \bs{u}\}\circ \{\bs{z}\rightarrow [\bs{w},\bs{x}\rangle \}\big)\nn\\
&\le Cp_c^9sD^{*4}(\bs{x})^2\nn\\
&=\text{\bf Er}(\bs{x}(s),4).
\end{align}
Together with \refeq{ob1E}, this implies \refeq{odbb1Eb1xC}. Then, we complete the proof. 
\QED

\Proof{\it Proof of Proposition~\ref{prp: Pix1xt4}.}For the first term of \refeq{lmmeqPix1xt4}, we combine Lemma~\ref{lmm: sumPxxosPi0} with Proposition~\ref{prp: Pi0xxt4} to obtain 
\begin{align}\lbeq{k}
\sum_{\bs{s}}\mathcal P_{p_c,[\bs{o},\bs{s}\rangle}^{\bs{x},\bs{x}}=
\begin{cases}
8s^4+O(d^{-5})&[\bs{x}(s)=o],\\
106s^6+O(d^{-7})&[\bs{x}(s)=e_1+e_2],\\
312s^8&[\bs{x}(s)=e_1+e_2+e_3+e_4].
\end{cases}
\end{align}

Next, we consider the second term in \refeq{lmmeqPix1xt4}. By Lemma~\ref{lmm:osxoy} and the translation invariance, we have
\begin{align}
\sum_{\bs{v}}\mathcal P_{p_c,[\bs{u},\bs{v}\rangle}^{\bs{x},\bs{x}}-\raisebox{-13pt}{\includegraphics[scale=0.85]{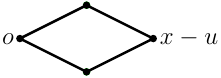}}\le \frac{p_c^6}{2}D^{*2}(\bs{x}-\bs{u})^3. 
\end{align}
Since 
\begin{align}
\frac{p_c^6}{2}\sum_{\bs{u}}\Pi_{p_c}^{\sss(0)}(\bs{u})D^{*2}(\bs{x}-\bs{u})^3\le \text{\bf Er}(\bs{x}(s),4)
\end{align}
which follows from Proposition~\ref{prp: Pi0xxt2} and \refeq{DmDn}, we obtain
\begin{align}\lbeq{a}
\sum_{\substack{[\bs{u}, \bs{v}\rangle \\ \bs{u}(t)=2}}\mathcal P_{p_c}^{\bs{u},\bs{u}}\mathcal P_{p_c,[\bs{u},\bs{v}\rangle}^{\bs{x},\bs{x}}&=\sum_{\bs{u}}\Pi_{p_c}^{\sss(0)}(\bs{u})\sum_{\bs{v}}\mathcal P_{p_c,[\bs{u},\bs{v}\rangle}^{\bs{x},\bs{x}}\nn\\
&=\sum_{\bs{u}}\Pi_{p_c}^{\sss(0)}(\bs{u})\raisebox{-13pt}{\includegraphics[scale=0.85]{2squre}}+\text{\bf Er}(\bs{x}(s),4).
\end{align}
Noting that $\raisebox{-11pt}{\includegraphics[scale=0.8]{2squre}}=0$ when $\bs{x}-\bs{u}=\pm 2e_i$ for $i=1,\cdots, d$, we deduce from  Proposition~\ref{prp: Pi0xxt2} and \refeq{Pi02xxmain} that 
\begin{align}\lbeq{b}
\sum_{\bs{u}}\Pi_{p_c}^{\sss(0)}(\bs{u})\raisebox{-13pt}{\includegraphics[scale=0.85]{2squre}}=
\begin{cases}
\frac{1}{2}s^4+O(d^{-5})&[\bs{x}(s)=o],\\
2s^6+O(d^{-7})&[\bs{x}(s)=e_1+e_2],\\
12s^8+O(d^{-10})&[\bs{x}(s)=e_1+e_2+e_3+e_4].
\end{cases}
\end{align}
which together with \refeq{lmmeqPix1xt4}, \refeq{a} and \refeq{k} provides \refeq{prpPix1xt4}. Furthermore, combining these with \refeq{sumofpix4} with $\Pi_{p_c}^{\sss(0)}(\bs{x})$ replaced by $\Pi_{p_c}^{\sss(1)}(\bs{x})$, we obtain \refeq{sumofpi4x1}. 
\QED

\section{Analysis of $\Pi_p^{(2)}(\bs{x})$}\label{Sectionpi2}
In this section, we estimate $\Pi_{p}^{\sss(2)}(\bs{x})$, as defined by 
\begin{align}
\Pi_{p}^{\sss(2)}(\bs{x})=\sum_{b_1, b_2 \in \mathbb B}\mathbb P_p\big(\{\bs{o}\Rightarrow \bs{\underline{b}_1}\}\cap E(b_1,\bs{\underline{b}_2}: \tilde C^{b_1}(\bs{o})) \cap E(b_2,\bs{x}: \tilde C^{b_2}(\bs{\bar{b}_1}))\big).
\end{align}
We note that $\bs{\underline{b}_2}$ cannot coincides with $\bs{\bar{b}_1}$, since if this were the case, $\bs{\underline{b}_2}$ cannot be reached from $\bs{o}$ without using $b_1$. 

\subsection{The case where the time component equals 3}
In this section, we estimate $\Pi_{p}^{\sss(2)}(\bs{x})$ for $\bs{x}(t)=3$, establishing the following proposition. 
\begin{prp}\label{prp: Pix2xt3}
For $\bs{x}$ such that $\bs{x}(t)=3$,
\begin{align}\lbeq{prpPix2xt3}
\Pi_{p}^{\sss(2)}(\bs{x})=
\begin{cases}
3s^5+O(d^{-6})&[\bs{x}(s)=e_1],\\
O(d^{-7})&[\bs{x}(s)=2e_1+e_2],\\
6s^7+O(d^{-9})&[\bs{x}(s)=e_1+e_2+e_3].
\end{cases}
\end{align}
Then, we have 
\begin{align}\lbeq{sumofpix2x3final}
\sum_{\bs{x}:\bs{x}(t)=3}\Pi_{p_c}^{\sss(2)}(\bs{x})=4s^4+O(d^{-5}). 
\end{align}
\end{prp}

\begin{lmm}\label{lmm:Pix2decom}
For $\bs{x}$ such that $\bs{x}(t)=3$, 
\begin{align}\lbeq{lmmeqPix2decom}
\Pi_{p}^{\sss(2)}(\bs{x})=\sum_{\substack{\bs{s},\bs{u}, \bs{v}\\\bs{v}\neq\bs{u}}}\mathcal P_{p_c,[\bs{o},\bs{s}\rangle}^{\bs{u},\bs{u}}D(\bs{v}-\bs{s})D(\bs{x}-\bs{u})D(\bs{x}-\bs{v})+\text{\bf Er}(\bs{x}(s),3).
\end{align}
\end{lmm}

\Proof{\it Proof of Lemma~\ref{lmm:Pix2decom}.}
In the case where $\bs{x}(t)=3$, there is only one possible configuration for $b_1$ and $b_2$, namely, $\bs{\underline{b}_1}=\bs{o}$ and $\bs{\bar{b}_2}=\bs{x}$. Therefore, we obtain
\begin{align}
\Pi_{p}^{\sss(2)}(\bs{x})&=\sum_{\bs{s},\bs{u}}\mathbb P_p\Big(E\big([\bs{o}, \bs{s}\rangle,\bs{u}: \tilde C^{[\bs{o}, \bs{s}\rangle}(\bs{o})\big) \cap E\big([\bs{u},\bs{x}\rangle ,\bs{x}: \tilde C^{[\bs{u},\bs{x}\rangle}(\bs{s})\big)\Big)\nn\\
&=\sum_{\bs{s},\bs{u}}\mathbb P_p\Big(\big\{\{[\bs{o}, \bs{s}\rangle \rightarrow \bs{u}\} \circ \{\bs{o} \rightarrow \bs{u}\} \big\} \cap \big\{\{ [\bs{u}, \bs{x}\rangle~\text{open}\}\circ \{\bs{s} \rightarrow \bs{x} \}\big\} \Big). 
\end{align}
We note that
\begin{align}
\{ [\bs{u}, \bs{x}\rangle~\text{open}\}\circ \{\bs{s} \rightarrow \bs{x} \}=\bigcup_{\substack{\bs{v}\\\bs{v}\neq\bs{u}}}\big\{[\bs{u}, \bs{x}\rangle, [\bs{s},\bs{v}\rangle, [\bs{v},\bs{x}\rangle~\text{open}\big\}\cap  E_{\succ}\big([\bs{s},\bs{v}\rangle\cup[\bs{v},\bs{x}\rangle; [\bs{u},\bs{x}\rangle \big).
\end{align}
Since by using that $E_{\succ}\big([\bs{s},\bs{v}\rangle\cup[\bs{v},\bs{x}\rangle; [\bs{u},\bs{x}\rangle \big)^c \subset \bigcup_{\substack{\bs{w} \\ \bs{w}\neq \bs{u},\bs{v}}}\{[\bs{s},\bs{w}\rangle [\bs{w},\bs{x}\rangle~\text{open} \}$ and by applying the BK inequality, we obtain
\begin{align}
&\sum_{\substack{\bs{s},\bs{u}, \bs{v}\\\bs{v}\neq\bs{u}}}\mathbb P_{p_c}\Big(\big\{\{[\bs{o}, \bs{s}\rangle \rightarrow \bs{u}\} \circ \{\bs{o} \rightarrow \bs{u}\} \big\} \cap \big\{[\bs{u}, \bs{x}\rangle, [\bs{s},\bs{v}\rangle, [\bs{v},\bs{x}\rangle~\text{open}\big\}\cap  E_{\succ}\big([\bs{s},\bs{v}\rangle\cup[\bs{v},\bs{x}\rangle; [\bs{u},\bs{x}\rangle \big)^c\Big)\nn\\
&\le p_c^5\sum_{\substack{\bs{s},\bs{u}, \bs{v},\bs{w} \\\bs{v}\neq\bs{u} \\ \bs{w}\neq\bs{u},\bs{v}}}\mathbb P_{p_c}\Big(\{[\bs{o}, \bs{s}\rangle \rightarrow \bs{u}\} \circ \{\bs{o} \rightarrow \bs{u}\}\Big)D(\bs{x}-\bs{v})D(\bs{x}-\bs{w})D(\bs{x}-\bs{u})D(\bs{w}-\bs{s})D(\bs{v}-\bs{s})\nn\\
&\le p_c^7\sum_{\bs{s},\bs{u}}D(\bs{s})D(\bs{u}-\bs{s})\tau_{p_c}(\bs{u})D^{*2}(\bs{x}-\bs{s})^2D(\bs{x}-\bs{u})\nn\\
&\le p_c^9s^2D^{*3}(\bs{x})^2\nn\\
&=\text{\bf Er}^{\sss (2)}(\bs{x}(s),3).
\end{align}
Thus, we obtain
\begin{align}
\Pi_{p}^{\sss(2)}(\bs{x})
&=\sum_{\substack{\bs{s},\bs{u}, \bs{v}\\\bs{v}\neq\bs{u}}}\mathbb P_{p_c}\Big(\big\{\{[\bs{o}, \bs{s}\rangle \rightarrow \bs{u}\} \circ \{\bs{o} \rightarrow \bs{u}\} \big\} \cap \big\{[\bs{u}, \bs{x}\rangle, [\bs{s},\bs{v}\rangle, [\bs{v},\bs{x}\rangle~\text{open}\big\}\Big)+\text{\bf Er}^{\sss (2)}(\bs{x}(s),3)\nn\\
&=\sum_{\substack{\bs{s},\bs{u}, \bs{v}\\\bs{v}\neq\bs{u}}}\mathcal P_{p_c,[\bs{o},\bs{s}\rangle}^{\bs{u},\bs{u}}D(\bs{v}-\bs{s})D(\bs{x}-\bs{u})D(\bs{x}-\bs{v})+\text{\bf Er}(\bs{x}(s),3)
\end{align}
which complete the proof. \QED

\Proof{\it Proof of Proposition~\ref{prp: Pix2xt3}.}
By Lemma~\ref{lmm:osxoy} and Lemma~\ref{lmm:Pix2decom}, we obtain
\begin{align}
\Pi_{p}^{\sss(2)}(\bs{x})-\raisebox{-11pt}{\includegraphics[scale=0.75]{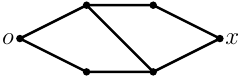}}
&\le \frac{p_c^6s^2}{2}D^{*3}(\bs{x})^2+\text{\bf Er}^{\sss (2)}(\bs{x}(s),3)=\text{\bf Er}^{\sss (2)}(\bs{x}(s),3). 
\end{align}
Thus, it suffices to estimate the leading term, namely $\raisebox{-10pt}{\includegraphics[scale=0.7]{M2a}}$, which can be estimated as in \refeq{ee}. Finally, \refeq{sumofpix2x3final} follows from \refeq{prpPix2xt3} and \refeq{sumofpix1x3decom} by replacing $\Pi_{p_c}^{\sss(0)}$ with $\Pi_{p_c}^{\sss(2)}$. 
\QED
\subsection{The case where the time component equals 4}
In this section, we show that $\Pi_{p}^{\sss(2)}(\bs{x})$ for $\bs{x}(t)=4$ contributes to the error term.  
\begin{prp}
For $\bs{x}$ such that $\bs{x}(t)=4$,
\begin{align}\lbeq{Pix2xt4final}
\Pi_{p}^{\sss(2)}(\bs{x})\le \text{\bf Er}(\bs{x}(s),4)
\end{align}
Then, we have 
\begin{align}
\sum_{\bs{x}:\bs{x}(t)=4}\Pi_{p_c}^{\sss(2)}(\bs{x})\le O(d^{-5}). 
\end{align}
\end{prp}
\Proof{\it Proof of Proposition~\ref{prp: Pix1xt4}.}
Similarly, in the case where $\bs{x}(t)=4$, the only possible configuration for $b_1$ is $\bs{\underline{b}_1}=\bs{o}$. For $b_2$, we consider the two cases: $\bs{\underline{b}_2}(t)=2$ or $\bs{\underline{b}_2}(t)=3$. Consequently, writing 
$\bs{b_2}=[\bs{u},\bs{v}\rangle$, we can rewrite $\Pi_{p}^{\sss(2)}(\bs{x})$ as
\begin{align}\lbeq{Pi2x4bd}
\Pi_{p}^{\sss(2)}(\bs{x})&=\sum_{\bs{s}}\sum_{\substack{[\bs{u},\bs{v}\rangle \\ \bs{u}(t)=2,3} }\mathbb P_p\big(E([\bs{o},\bs{s}\rangle,\bs{u}: \tilde C^{[\bs{o},\bs{s}\rangle}(\bs{o})) \cap E([\bs{u},\bs{v}\rangle,\bs{x}: \tilde C^{[\bs{u},\bs{v}\rangle}(\bs{s}))\big).
\end{align}

When $\bs{u}(t)=2$, the only pivotal bond for the event that $\bs{v}\rightarrow \bs{x}$ is $[\bs{v},\bs{x}\rangle$. Hence, 
by the BK inequality, the contribution from the case where $\bs{u}(t)=2$ can be bounded by 
\begin{align}
\sum_{\bs{s}}\sum_{\substack{[\bs{u},\bs{v}\rangle \\ \bs{u}(t)=2} }\mathbb P_p\big(\{\bs{o} \rightarrow \bs{u} \}\circ \{[\bs{o},\bs{s}\rangle \rightarrow [\bs{u},\bs{v}\rangle \rightarrow \bs{x} \} \circ \{\bs{s} \rightarrow \bs{x} \}\big)&\le p_c^9sD^{*4}(\bs{x})^2\nn\\
&=\text{\bf Er}(\bs{x}(s),4). 
\end{align}

For the case where $\bs{u}(t)=3$, if the event $\{\bs{s} \Rightarrow \bs{u}\}$ occurs, then by the BK inequality, it can be bounded as 
\begin{align}
&\sum_{\bs{s},\bs{u}}\mathbb P_p\big(E([\bs{o},\bs{s}\rangle,\bs{u}: \tilde C^{[\bs{o},\bs{s}\rangle}(\bs{o})) \cap E([\bs{u},\bs{x}\rangle,\bs{x}: \tilde C^{[\bs{u},\bs{x}\rangle}(\bs{s}))\cap \{\bs{s} \Rightarrow \bs{u}\} \big)\nn\\
&\le \sum_{\bs{s},\bs{u}}\sum_{\substack{\bs{v},\bs{w} \\ \bs{v}\neq \bs{w} }}\mathbb P_p\big(\{[\bs{o},\bs{s}\rangle \rightarrow \bs{w} \rightarrow [\bs{u},\bs{x}\rangle\} \circ \{\bs{s} \rightarrow \bs{v}  \rightarrow \bs{u}\}  \circ \{\bs{v} \rightarrow \bs{x}\} \circ \{\bs{o} \rightarrow \bs{w}\} \big)\nn\\
&+\sum_{\bs{s},\bs{u}}\sum_{\substack{\bs{w} \\ \bs{w} \neq \bs{u} }}\mathbb P_p\big(\{[\bs{o},\bs{s}\rangle \rightarrow \bs{w} \rightarrow [\bs{u},\bs{x}\rangle\} \circ \{\bs{s} \rightarrow \bs{u}\}  \circ \{\bs{w} \rightarrow \bs{x}\} \circ \{\bs{o} \rightarrow \bs{w}\} \big)\nn\\
&\le Cp_c^{10}s^2D^{*4}(\bs{x})^2\nn\\
&=\text{\bf Er}(\bs{x}(s),4). 
\end{align}
If the event $\{\bs{s} \Rightarrow \bs{u}\}$ does not occur, then we have 
\begin{align}
&\sum_{\bs{s},\bs{u}}\mathbb P_p\big(E([\bs{o},\bs{s}\rangle,\bs{u}: \tilde C^{[\bs{o},\bs{s}\rangle}(\bs{o})) \cap E([\bs{u},\bs{x}\rangle,\bs{x}: \tilde C^{[\bs{u},\bs{x}\rangle}(\bs{s}))\cap \{\bs{s} \Rightarrow \bs{u}\}^c\big)\nn\\
&\le  \sum_{\bs{s},\bs{u}}\sum_{\substack{\bs{w} \\ \bs{w}\neq\bs{u} }}\mathbb P_p\big(\{[\bs{o},\bs{s}\rangle \rightarrow \bs{w} \rightarrow [\bs{u},\bs{x}\rangle\} \circ \{\bs{w}  \rightarrow \bs{x}\}  \circ \{\bs{o} \rightarrow \bs{w}\} \big)\nn\\
&\le Cp_c^9sD^{*4}(\bs{x})^2\nn\\
&=\text{\bf Er}(\bs{x}(s),4).
\end{align}
Therefore, we conclude \refeq{Pix2xt4final}.\QED

\appendix 
\section{Proofs of Lemmas~\ref{lmm:osxoy} and \ref{lmm:ouovowestimate}}\label{Appendix}
In this appendix, we provide the proofs of Lemmas~\ref{lmm:ouovowestimate} and \ref{lmm:osxoy}, in Section~\ref{Lmm1} and \ref{Lmm2}, respectively. 

\subsection{Proof of Lemma~\ref{lmm:ouovowestimate}}\label{Lmm1}
We here restate Lemma~\ref{lmm:ouovowestimate}. 
\begin{lmm}\label{lmm:Aouovowestimate}
Let $\bs{u},\bs{v}$ be points such that $\bs{u}(t)=\bs{v}(t)=2$, and  let $\bs{w}$ and $\bs{z}$ be distincts points with $\bs{w}(t)=\bs{z}(t)=2$, where $\bs{w}\neq \bs{u}, \bs{v}$ and $\bs{z}\neq \bs{u}, \bs{v}$. Then, for $p\le p_c$, 
\begin{align}\lbeq{Aouovowestimate}
\Big|\mathcal P_p^{\bs{u},\bs{v}}(\bs{w})-\zeta_p(\bs{u},\bs{v},\bs{w})\Big|\le 
p^6D^{*2}(\bs{u})D^{*2}(\bs{v})D^{*2}(\bs{w})
\end{align}
where
\begin{align}
\mathcal P_{p,[\bs{o},\bs{s}\rangle}^{\bs{x},\bs{y}}&:=\mathbb P_p\big(\{[\bs{o},\bs{s}\rangle \rightarrow \bs{x} \}\circ \{\bs{o}\rightarrow \bs{y}\}\big), \\
\zeta_p(\bs{u},\bs{v},\bs{w})&:=\frac{2-\delta_{\bs{u},\bs{v}}}{2}p\sum_{\substack{\bs{s} \\ \bs{s}(t)=1}}D(\bs{w}-\bs{s})\big[\mathcal P_{p,[\bs{o},\bs{s}\rangle}^{\bs{u},\bs{v}}+\mathcal P_{p,[\bs{o},\bs{s}\rangle}^{\bs{v},\bs{u}}\big].\lbeq{wsosuovouosvA}
\end{align}
Furthermore, in the case $\bs{u}=\bs{v}$, 
\begin{align}\lbeq{Aouovowestimate1}
\Big|\mathcal P_p^{\bs{u},\bs{u}}(\bs{w})-\Big[\zeta_p(\bs{u},\bs{u},\bs{w})+\frac{1}{2}\raisebox{-12pt}{\includegraphics[scale=0.75]{Ppuuw}}~\Big]\Big|\le 
p^8sD^{*2}(\bs{u})^2D^{*2}(\bs{w})
\end{align}
Additionally, 
\begin{align}\lbeq{APuvwz3}
\mathcal P_p^{\bs{u},\bs{v}}(\bs{w},\bs{z})\le 
\begin{cases}
Cs^2D^{*2}(\bs{u})D^{*2}(\bs{v})&[\bs{u} \neq\bs{v}],\\
Cs^2D^{*2}(\bs{u})D^{*2}(\bs{w})&[\bs{u} =\bs{v}].
\end{cases}
\end{align}
\end{lmm}
\Proof{\it Proof of Lemma~\ref{lmm:Aouovowestimate}.}We first consider \refeq{Aouovowestimate} and \refeq{Aouovowestimate1}, noting that
\begin{align}\lbeq{Xowwx} 
\mathcal P_p^{\bs{u},\bs{v}}(\bs{w})=\sum_{\substack{\omega_1:\bs{o} \rightarrow \bs{u}\\ \omega_2:\bs{o}\rightarrow \bs{v}\\ \omega_1\cap \omega_2=\varnothing}}\mathbb P_p\big(Y_{\sss \omega_1,\omega_2}\cap \{\bs{o}\rightarrow \bs{w}\}\big),
\end{align}
where $Y_{\sss \omega_1,\omega_2}$ is defined in \refeq{Ydef}. In the case where $\bs{u}=\bs{v}$, we impose the  additional condition $\omega_1\prec \omega_2$. Since the event $\{\bs{o} \rightarrow \bs{w}\}$ can be rewritten by a union of disjoint events
\begin{align}
\{\bs{o} \rightarrow \bs{w}\}=\bigcup_{\omega_3: \bs{o}\rightarrow \bs{w}}\{\omega_3~\text{occupied}\}\cap E_{\prec}(\omega_3), 
\end{align}
the left-hand side of \refeq{Xowwx} can be rewritten as
\begin{align}\lbeq{Xo3E3}
\mathcal P_p^{\bs{u},\bs{v}}(\bs{w})=\sum_{\substack{\omega_1:\bs{o} \rightarrow \bs{u}\\ \omega_2:\bs{o}\rightarrow \bs{v}\\\omega_3:\bs{o}\rightarrow \bs{w} \\ \omega_1\cap \omega_2=\varnothing}}\mathbb P_p\big(Y_{\sss \omega_1,\omega_2}\cap \{\omega_3~\text{occupied}\}\cap E_{\prec}(\omega_3)\big). 
\end{align}
For the additional conditions on the relation between $\omega_1,\omega_2$ and $\omega_3$, we consider three cases: $\omega_1\cap \omega_3\neq \varnothing$, $\omega_2\cap \omega_3\neq \varnothing$ and $\omega_1\cap \omega_2\cap \omega_3=\varnothing$. 

We begin with the case where $\omega_1\cap \omega_2\cap \omega_3=\varnothing$. In the case where $\bs{u}\neq\bs{v}$, since $\omega_1\cap \omega_2 \cap \omega_3=\varnothing$, the events $\{\omega_i ~\text{occupied}\}$ with $i=1,2,3$ are mutually independent, dropping the condition $\omega_1\cap \omega_2 \cap \omega_3=\varnothing$ in the final step, we obtain
\begin{align}
\sum_{\substack{\omega_1:\bs{o} \rightarrow \bs{u}\\ \omega_2:\bs{o}\rightarrow \bs{v}\\\omega_3:\bs{o}\rightarrow \bs{w} \\ \omega_1\cap \omega_2=\varnothing \\ \omega_1\cap \omega_2 \cap \omega_3=\varnothing}}\mathbb P_p\big(Y_{\sss \omega_1,\omega_2}\cap \{\omega_3~\text{occupied}\}\cap E_{\prec}(\omega_3)\big) &\le \sum_{\substack{\omega_1:\bs{o} \rightarrow \bs{u}\\ \omega_2:\bs{o}\rightarrow \bs{v}\\\omega_3:\bs{o}\rightarrow \bs{w} \\ \omega_1\cap \omega_2 \cap \omega_3=\varnothing}}\mathbb P_p\big(\{\omega_1,\omega_2, \omega_3 ~\text{occupied}\}\big)\nn\\
&\le p^6D^{*2}(\bs{u})D^{*2}(\bs{v})D^{*2}(\bs{w}).
\end{align}
In the case where $\bs{u}=\bs{v}$, observing that the events $E_{\prec}(\omega_1)$ and $E_{\succ}(\omega_2:\omega_1)$ are included in the event $Y_{\sss \omega_1,\omega_2}$, we obtain
\begin{align}\lbeq{omegaempty}
&\sum_{\substack{\omega_1,\omega_2:\bs{o} \rightarrow \bs{u}\\ \omega_3:\bs{o}\rightarrow \bs{w} \\ \omega_1\cap \omega_2=\varnothing \\ \omega_1\cap \omega_2 \cap \omega_3=\varnothing \\\omega_1\prec \omega_2}}\mathbb P_p\big(Y_{\sss \omega_1,\omega_2}\cap \{\omega_3~\text{occupied}\}\cap E_{\prec}(\omega_3)\big)-\sum_{\substack{\omega_1,\omega_2:\bs{o} \rightarrow \bs{u} \\\omega_3:\bs{o}\rightarrow \bs{w} \\ \omega_1\cap \omega_2=\varnothing\\ \omega_1\cap \omega_2 \cap \omega_3=\varnothing \\\omega_1\prec \omega_2 }}\mathbb P_p\big(\{\omega_1,\omega_2, \omega_3 ~\text{occupied}\}\big)\nn\\
&\le \sum_{\substack{\omega_1,\omega_2:\bs{o} \rightarrow \bs{u}\\ \omega_3:\bs{o}\rightarrow \bs{w} \\ \omega_1\cap \omega_2 \cap \omega_3=\varnothing}}\mathbb P_p\Big(\{\omega_1,\omega_2, \omega_3 ~\text{occupied}\}\cap \big(E_{\prec}(\omega_1)\cap E_{\succ}(\omega_2:\omega_1)\cap E_{\prec}(\omega_3)\big)^c\Big)\nn\\
&\le p^6D^{*2}(\bs{u})^2D^{*2}(\bs{w})\big[2D^{*2}(\bs{u})+D^{*2}(\bs{w})\big]\nn\\
&\le CsD^{*2}(\bs{u})^2D^{*2}(\bs{w}).
\end{align}
Thus, it suffices to estimate the second term in the left-hand side of \refeq{omegaempty}, and we obtain 
\begin{align}
\sum_{\substack{\omega_1,\omega_2:\bs{o} \rightarrow \bs{u}\\ \omega_3:\bs{o}\rightarrow \bs{w} \\ \omega_1\cap \omega_2=\varnothing\\ \omega_1\cap \omega_2 \cap \omega_3=\varnothing \\\omega_1\prec \omega_2}}\mathbb P_p\big(\{\omega_1,\omega_2, \omega_3 ~\text{occupied}\}\big)&=\sum_{\substack{\bs{s},\bs{t},\bs{z} \\ \bs{u}-\bs{s}\prec \bs{u}-\bs{t} \\ \bullet}}\mathbb P_p\big(\{[\bs{o},\bs{s}\rangle[\bs{s},\bs{u}\rangle,[\bs{o},\bs{t}\rangle[\bs{t},\bs{u}\rangle,[\bs{o},\bs{z}\rangle[\bs{z},\bs{w}\rangle~\text{open}\}\big)\nn\\
&=\frac{1}{2}\raisebox{-14pt}{\includegraphics[scale=0.85]{Ppuuw}}.
\end{align}

Next we consider the remaining cases where  $\omega_1\cap \omega_3\neq \varnothing$ and $\omega_2\cap \omega_3\neq \varnothing$, noting that these events can not occur at the same time, and that when these events occur, $\omega_3$ shares only one bond incident to $\bs{o}$ with either $\omega_1$ or $\omega_2$. Therefore, for the case where $\omega_1\cap \omega_3\neq \varnothing$, we obtain
\begin{align}\lbeq{o1o2o3}
&\sum_{\substack{\omega_1:\bs{o} \rightarrow \bs{u}\\ \omega_2:\bs{o}\rightarrow \bs{v}\\\omega_3:\bs{o}\rightarrow \bs{w} \\ \omega_1\cap\omega_2=\varnothing \\\omega_1\cap\omega_3\neq \varnothing}}\mathbb P_p\big(Y_{\sss \omega_1,\omega_2}\cap \{\omega_3~\text{occupied}\}\cap E_{\prec}(\omega_3)\big)\nn\\
&=\sum_{\bs{s}}\sum_{\substack{\omega_1:[\bs{o},\bs{s}\rangle \rightarrow \bs{u}\\ \omega_2:\bs{o}\rightarrow \bs{v}\\ \omega_1\cap \omega_2=\varnothing}}\mathbb P_p\big(Y_{\sss \omega_1,\omega_2}\cap \{[\bs{s},\bs{w}\rangle~\text{open} \}\cap E_{\prec}([\bs{o},\bs{s}\rangle \cup[\bs{s},\bs{w}\rangle)\big),
\end{align}
where the condition $\omega_1 \prec \omega_2$ is required when $\bs{u}=\bs{v}$. The right-hand side of \refeq{o1o2o3} can be decomposed into two parts: one corresponding to the case where we disregard $ E_{\prec}([\bs{o},\bs{s}\rangle \cup[\bs{s},\bs{w}\rangle)$ and the other where we replace $ E_{\prec}([\bs{o},\bs{s}\rangle \cup[\bs{s},\bs{w}\rangle)$ with its complement $ E_{\prec}([\bs{o},\bs{s}\rangle \cup[\bs{s},\bs{w}\rangle)^c$. 

For the second case, noting that 
\begin{align}
E_{\prec}([\bs{o},\bs{s}\rangle \cup[\bs{s},\bs{w}\rangle)^c=\bigcup_{\substack{\omega_3': \bs{o}\rightarrow \bs{w}\\\omega'_3\prec [\bs{o},\bs{s}\rangle \cup[\bs{s},\bs{w}\rangle}}\{\omega_3~\text{occupied}\}\cap E_{\prec}(\omega'_3),
\end{align}
we obtain 
\begin{align}\lbeq{YswEc}
&\sum_{\bs{s}}\sum_{\substack{\omega_1:[\bs{o},\bs{s}\rangle \rightarrow \bs{u}\\ \omega_2:\bs{o}\rightarrow \bs{v}\\ \omega_1\cap \omega_2=\varnothing}}\mathbb P_p\big(Y_{\sss \omega_1,\omega_2}\cap \{[\bs{s},\bs{w}\rangle~\text{open} \}\cap E_{\prec}([\bs{o},\bs{s}\rangle \cup[\bs{s},\bs{w}\rangle)^c\big)\nn\\
&\le p\sum_{\bs{s}}D(\bs{w}-\bs{s})\sum_{\substack{\omega_1:[\bs{o},\bs{s}\rangle \rightarrow \bs{u}\\ \omega_2:\bs{o}\rightarrow \bs{v}\\ \omega_1\cap \omega_2=\varnothing}}\sum_{\substack{\omega_3': \bs{o}\rightarrow \bs{w}\\\omega'_3\prec [\bs{o},\bs{s}\rangle \cup[\bs{s},\bs{w}\rangle}}\mathbb P_p\big(\{\omega_1,\omega_2, \omega'_3 ~\text{occupied}\}\big)\nn\\
&\le \frac{p^7}{2}D^{*2}(\bs{v})D^{*2}(\bs{w})\sum_{\substack{\bs{s} \\ \bs{s}(t)=1}}D(\bs{w}-\bs{s})D(\bs{s})D(\bs{u}-\bs{s})\nn\\
&\le \frac{p^7s}{2}D^{*2}(\bs{u})D^{*2}(\bs{v})D^{*2}(\bs{w}). 
\end{align}

In the first case, since $\bs{w}$ is distinct from both $\bs{u}$ and $\bs{v}$, the events $Y_{\sss \omega_1,\omega_2}$ and $\{[\bs{s},\bs{w}\rangle~\text{open} \}$ are independent. Thus, we obtain
\begin{align}\lbeq{Yswa}
\sum_{\bs{s}}\sum_{\substack{\omega_1:[\bs{o},\bs{s}\rangle \rightarrow \bs{u}\\ \omega_2:\bs{o}\rightarrow \bs{v}\\ \omega_1\cap \omega_2=\varnothing}}\mathbb P_p\big(Y_{\sss \omega_1,\omega_2}\cap \{[\bs{s},\bs{w}\rangle~\text{open} \}\big)=p\sum_{\substack{\bs{s} \\ \bs{s}(t)=1}}D(\bs{w}-\bs{s})\sum_{\substack{\omega_1:[\bs{o},\bs{s}\rangle \rightarrow \bs{u}\\ \omega_2:\bs{o}\rightarrow \bs{v}\\ \omega_1\cap \omega_2=\varnothing}}\mathbb P_p\big(Y_{\sss \omega_1,\omega_2}\big).
\end{align}
Since \refeq{i} also holds when replacing one instance of $\bs{x}$ with $\bs{y}$, we obtain for $\bs{x}$ and $\bs{y}$ satisfying $\bs{x}(t)=\bs{y}(t)=2$ 
\begin{align}\lbeq{osuov1}
\mathcal P_{p,[\bs{o},\bs{s}\rangle}^{\bs{x},\bs{y}}
&=\sum_{\substack{\omega_1: \bs{o}\rightarrow [\bs{x}-\bs{s}, \bs{x} \rangle \\ \omega_2: \bs{o}\rightarrow \bs{y} \\ \omega_1\cap\omega_2=\varnothing }}\mathbb P_{p}\big(\{\omega_1,\omega_2~\text{occupied}\} \cap E_{\succ_{\bs{s}}}(\omega_2:\omega_1 )\big),
\end{align}
the latter summation in \refeq{Yswa} can be decomposed as 
\begin{align}
\sum_{\substack{\omega_1:[\bs{o},\bs{s}\rangle \rightarrow \bs{u}\\ \omega_2:\bs{o}\rightarrow \bs{v}\\ \omega_1\cap \omega_2=\varnothing}}\mathbb P_p\big(Y_{\sss \omega_1,\omega_2}\big)&=\frac{\delta_{\bs{u},\bs{v}}}{2}\mathcal P_{p,[\bs{o},\bs{s}\rangle}^{\bs{u},\bs{v}}\nn\\
&-\sum_{\substack{\omega_1:[\bs{o},\bs{s}\rangle \rightarrow \bs{u}\\ \omega_2:\bs{o}\rightarrow \bs{v}\\ \omega_1\cap \omega_2=\varnothing}}\mathbb P_p\big(\{\omega_1,\omega_2~\text{occupied}\} \cap E_{\prec}(\omega_1)^c\cap E_{\succ}(\omega_2:\omega_1)\big). 
\end{align}
Here, the factor $\frac{1}{2}$ comes from the condition $\omega_1\prec \omega_2$. 
In a similar manner as in the derivation \refeq{YswEc}, we establish
\begin{align}\lbeq{DYEcE}
&\sum_{\bs{s}}D(\bs{w}-\bs{s})\sum_{\substack{\omega_1:[\bs{o},\bs{s}\rangle \rightarrow \bs{u}\\ \omega_2:\bs{o}\rightarrow \bs{v}\\ \omega_1\cap \omega_2=\varnothing}}\mathbb P_p\big(\{\omega_1,\omega_2~\text{occupied}\} \cap E_{\prec}(\omega_1)^c\cap E_{\succ}(\omega_2:\omega_1)\big)\nn\\
&\le  \frac{p^7s}{2}D^{*2}(\bs{u})D^{*2}(\bs{v})D^{*2}(\bs{w}). 
\end{align}
The case $\omega_2\cap \omega_3\neq \varnothing$ can be handled in a similar manner. Combining \refeq{o1o2o3}--\refeq{DYEcE}, we then obtain \refeq{Aouovowestimate} and \refeq{Aouovowestimate1}.

Finally, we prove \refeq{APuvwz3}. We partition the event under consideration into two cases: the first is 
$\big\{\left\{\{\bs{o}\rightarrow \bs{u}\}\circ\{\bs{o}\rightarrow \bs{v}\}\right\}\cap \{\bs{o}\rightarrow \bs{w}\}\big\}\circ \{\bs{o}\rightarrow \bs{z}\}$ and the second is its complement. By applying the BK inequality, the contribution from the first mentioned event can be bounded above by $\mathcal P_p^{\bs{u},\bs{v}}(\bs{w})\tau_p(\bs{z})$. 

For the contribution from the complementary event, we can analyze the cases based on where the path from $\bs{o}$ to $\bs{w}$ and the path from $\bs{o}$ to $\bs{z}$ diverge from the paths originating from $\bs{o}$ to $\bs{u}$ or $\bs{v}$.
Since we are considering the complementary event of $\big\{\left\{\{\bs{o}\rightarrow \bs{u}\}\circ\{\bs{o}\rightarrow \bs{v}\}\right\}\cap \{\bs{o}\rightarrow \bs{w}\}\big\}\circ \{\bs{o}\rightarrow \bs{z}\}$, the path from $\bs{o}$ to $\bs{z}$ must intersect with at least one of the paths from $\bs{o}$ to $\bs{u}$, $\bs{o}$ to $\bs{v}$ or $\bs{o}$ to $\bs{w}$. Moreover, we note that since the time-components of the endpoints of these paths are 2, the path from $\bs{o}$ to $\bs{z}$ cannot have intersections with more than three of these paths. When $\{\bs{o}\rightarrow \bs{u}\}\circ\{\bs{o}\rightarrow \bs{v}\}$ and $\{\bs{o}\rightarrow \bs{w}\}$ occur bond-disjointly, there are three possible patterns for how the path from $\bs{o}$ to $\bs{z}$ diverges: it can diverges from the path from $\bs{o}$ to $\bs{u}$, from $\bs{o}$ to $\bs{v}$ or from $\bs{o}$ to $\bs{w}$. However, when these events do not occur bond-disjointly, there are two possible patterns: either the paths from $\bs{o}$ to $\bs{z}$, and from $\bs{o}$ to $\bs{w}$ intersect, or they do not. Consequently, by applying the BK inequality, we obtain 
\begin{align}\lbeq{complementevent1}
&\mathcal P_p^{\bs{u},\bs{v}}(\bs{w},\bs{z})-\mathbb P_p\big(\big\{\left\{\{\bs{o}\rightarrow \bs{u}\}\circ\{\bs{o}\rightarrow \bs{v}\}\right\}\cap \{\bs{o}\rightarrow \bs{w}\}\big\}\circ \{\bs{o}\rightarrow \bs{z}\}\big)\nn\\
&\le p\sum_{(\bs{a},\bs{b}) \in \{\bs{u},\bs{v}, \bs{w}\}}\tau_p(\bs{a})\tau_p(\bs{b})\sum_{\substack{\bs{s}\\ \bs{s}(t)=1}}D(\bs{z}-\bs{s})\mathbb P_p\big([\bs{o},\bs{s}\rangle \rightarrow \bs{c}\big)\nn\\
&+p^2\sum_{\bs{a}\in \{\bs{u},\bs{v}\}}\sum_{\substack{\bs{s},\bs{t}\\ \bs{s}(t)=\bs{t}(t)=1 \\ \bs{s}\neq\bs{t}}}D(\bs{w}-\bs{t})D(\bs{z}-\bs{s})\mathbb P_p\big([\bs{o},\bs{s}\rangle \rightarrow \bs{a}\big)\mathbb P_p\big([\bs{o},\bs{t}\rangle \rightarrow \bs{\dot{a}}\big)\nn\\
&+p^2\sum_{\bs{a}\in \{\bs{u},\bs{v}\}}\tau_p(\bs{a})\sum_{\substack{\bs{s}\\ \bs{s}(t)=1}}D(\bs{w}-\bs{s})D(\bs{z}-\bs{s})\mathbb P_p\big([\bs{o},\bs{s}\rangle \rightarrow \bs{\dot{a}}\big)
\end{align}
where $\bs{c}$ denotes the remaining point contained in $\{\bs{u},\bs{v},\bs{w}\}$ other than $\bs{a}$ and $\bs{b}$, and $\bs{\dot{a}}$ represents the other point in $\{\bs{u}, \bs{v}\}$ distinct from $\bs{a}$. For the case where $\bs{u}=\bs{v}$, the factor $\frac{1}{2}$ on the right-hand side is required due to symmetry. The first term in the right-hand side can be bounded as 
\begin{align}
p\sum_{(\bs{a},\bs{b}) \in \{\bs{u},\bs{v}, \bs{w}\}}\tau_p(\bs{a})\tau_p(\bs{b})\sum_{\substack{\bs{s}\\ \bs{s}(t)=1}}D(\bs{z}-\bs{s})\mathbb P_p\big([\bs{o},\bs{s}\rangle \rightarrow \bs{c}\big)&\le p^3s\sum_{(\bs{a},\bs{b}) \in \{\bs{u},\bs{v}, \bs{w}\}}\tau_p(\bs{a})\tau_p(\bs{b})D^{*2}(\bs{c})\nn\\
&\le 3p^7sD^{*2}(\bs{u})D^{*2}(\bs{v})D^{*2}(\bs{w}).
\end{align}

For the second and third terms, in the case where $\bs{u}\neq\bs{v}$,we bound the factors $D(\bs{w}-\bs{t})D(\bs{z}-\bs{s})$ and $D(\bs{w}-\bs{s})D(\bs{z}-\bs{s})$ in the summands  by $s^2$ respectively, and thus we bound the sum of the second and third terms by 
\begin{align}
p^4s^2\sum_{\bs{a}\in \{\bs{u},\bs{v}\}}\{p^2D^{*2}(\bs{a})+\tau_p(\bs{a})\}D^{*2}(\bs{\dot{a}})\le 4p^6s^2D^{*2}(\bs{u})D^{*2}(\bs{v}).
\end{align}
In the case where $\bs{u}=\bs{v}$, we bound it as 
\begin{align}
&p^2\sum_{\bs{s}}D(\bs{z}-\bs{s})\mathbb P_p\big([\bs{o},\bs{s}\rangle \rightarrow \bs{u}\big)\Big[\sum_{\substack{\bs{t}\\ \bs{t}\neq\bs{s}}}D(\bs{w}-\bs{t})\mathbb P_p\big([\bs{o},\bs{t}\rangle \rightarrow \bs{u}\big)+\tau_p(\bs{u})D(\bs{w}-\bs{s})\Big]\nn\\
&\le p^6s\sum_{\bs{s}}D(\bs{s})D(\bs{u}-\bs{s})\Big[\sum_{\substack{\bs{t}\\ \bs{t}\neq\bs{s}}}D(\bs{w}-\bs{t})D(\bs{t})D(\bs{u}-\bs{t})+D^{*2}(\bs{u})D(\bs{w}-\bs{s})\Big]\nn\\
&\le p^6s^2D^{*2}(\bs{u})D^{*2}(\bs{w})+p^6sD^{*2}(\bs{u})\sum_{\bs{s}}D(\bs{s})D(\bs{u}-\bs{s})D(\bs{w}-\bs{s})\nn\\
&\le2p^6s^2D^{*2}(\bs{u})D^{*2}(\bs{w}). 
\end{align}

Therefore, we obtain
\begin{align}\lbeq{ouovowozbound} 
&\mathcal P_p^{\bs{u},\bs{v}}(\bs{w},\bs{z})\le \mathcal P_p^{\bs{u},\bs{v}}(\bs{w})\tau_p(\bs{z})+3p^7sD^{*2}(\bs{u})D^{*2}(\bs{v})D^{*2}(\bs{w})+
\begin{cases}
4p^6s^2D^{*2}(\bs{u})D^{*2}(\bs{v})&[\bs{u}\neq \bs{v}],\\
2p^6s^2D^{*2}(\bs{u})D^{*2}(\bs{w})&[\bs{u}=\bs{v}]. 
\end{cases}
\end{align}
By applying the BK inequality, we obtain
\begin{align}\lbeq{bdofzeta}
\zeta_{p_c}(\bs{u},\bs{v},\bs{w})\le p_c^3s\big\{\tau_{p_c}(\bs{v})D^{*2}(\bs{u})+\tau_{p_c}(\bs{u})D^{*2}(\bs{v})\big\}\le 2p_c^5sD^{*2}(\bs{u})D^{*2}(\bs{v}),
\end{align}
which, along with \refeq{Aouovowestimate} and the bound $\sup_{\bs{x}}D^{*2}(\bs{x})\le s$, yields
\begin{align}
\mathcal P_p^{\bs{u},\bs{v}}(\bs{w}) &\le CsD^{*2}(\bs{u})D^{*2}(\bs{v}). \lbeq{Puvw3}
\end{align}
By applying \refeq{Puvw3} to \refeq{ouovowozbound} together with $\sup_{\bs{x}}D^{*2}(\bs{x})\le s$, we obtain \refeq{APuvwz3}, and thus complete the proof of Lemma~\ref{lmm:Aouovowestimate}.\QED

\subsection{Proof of Lemma~\ref{lmm:osxoy}}\label{Lmm2}
Lemma~\ref{lmm:osxoy} is restated below. 
\begin{lmm}\label{lmm:Aosxoy}
Let $\bs{x}$ and $\bs{y}$ be points such that $\bs{x}(t)=\bs{y}(t)=2$. For the both case where $\bs{x}\neq \bs{y}$ and $\bs{x}=\bs{y}$, 
\begin{align}\lbeq{osxoy2A}
\Big|\mathcal P_{p,[\bs{o},\bs{s}\rangle}^{\bs{x},\bs{y}}-\raisebox{-12pt}{\includegraphics[scale=0.75]{Pposxymain}}\Big|\le \frac{1}{2}\raisebox{-12pt}{\includegraphics[scale=0.75]{Pposxyerror}}
\end{align}
\end{lmm}
\Proof{\it Proof of Lemma~\ref{lmm:Aosxoy}.}
By \refeq{osuov1}, we have
\begin{align}\lbeq{osuov1AA}
\mathcal P_{p,[\bs{o},\bs{s}\rangle}^{\bs{x},\bs{y}}
&=\sum_{\substack{\omega:\bs{o}\rightarrow \bs{y}\\ \omega \cap [\bs{o},\bs{s}\rangle [\bs{s},\bs{x} \rangle  =\varnothing}}\mathbb P_p\big(\{[\bs{o},\bs{s}\rangle [\bs{s},\bs{x} \rangle ,\omega~\text{occupied}\}\cap E_{\succ_{\bs{s}}}(\omega:[\bs{o},\bs{s}\rangle [\bs{s},\bs{x} \rangle)\big).
\end{align}
Since we obtain 
\begin{align}
&\sum_{\substack{\omega:\bs{o}\rightarrow \bs{y}\\ \omega \cap [\bs{o},\bs{s}\rangle [\bs{s},\bs{x} \rangle  =\varnothing}}\mathbb P_p\big(\{[\bs{o},\bs{s}\rangle [\bs{s},\bs{x} \rangle ,\omega~\text{occupied}\}\cap E_{\succ_{\bs{s}}}(\omega:[\bs{o},\bs{s}\rangle [\bs{s},\bs{x} \rangle)^c\big).\nn\\
&\le \sum_{\substack{\bs{t}\\ \bs{t}\neq \bs{s} }}\mathbb P_p\big(\{[\bs{o},\bs{s}\rangle[\bs{s},\bs{x}\rangle,[\bs{o},\bs{t}\rangle[\bs{t},\bs{y}\rangle~\text{occupied}\} \cap E_{\succ_{\bs{s}}}([\bs{o},\bs{t}\rangle[\bs{t},\bs{y}\rangle:[\bs{o},\bs{s}\rangle[\bs{s},\bs{x}\rangle)^c\big)\nn\\
&\le \sum_{\substack{\bs{t}, \bs{w}\\ \bs{t},\bs{w}\neq \bs{s} \\ \bs{w}-\bs{y}\succ \bs{t}-\bs{y}}}\mathbb P_p\big(\{[\bs{o},\bs{s}\rangle[\bs{s},\bs{x}\rangle,[\bs{o},\bs{t}\rangle[\bs{t},\bs{y}\rangle, [\bs{o},\bs{w}\rangle[\bs{w},\bs{y}\rangle~\text{occupied}\} \big)\nn\\
&\le \frac{p^6}{2}D(\bs{s})D(\bs{x}-\bs{s})D^{*2}(\bs{y})^2, 
\end{align}
together with \refeq{osuov1AA}, we establish \refeq{osxoy2A} for the both cases $\bs{x}\neq \bs{y}$ and $\bs{x}= \bs{y}$.
\QED

%\begin{figure}[t]
%\begin{align}
%&F_{1,p}(\bs{x},\bs{y})\le\raisebox{-18pt}{\includegraphics[scale=0.7]{F1x3}}+\raisebox{-18pt}{\includegraphics[scale=0.7]{F1x30}},\nn\\
%&F_{2,p}(\bs{x},\bs{y})\le Cs\Bigg[\raisebox{-18pt}{\includegraphics[scale=0.7]{F2x3}}+\raisebox{-18pt}{\includegraphics[scale=0.7]{F2x30}}\Bigg],\nn\\
%&F_{3,p}(\bs{x},\bs{y})\le Cs\Bigg[\raisebox{-18pt}{\includegraphics[scale=0.7]{F3x3}}+\raisebox{-18pt}{\includegraphics[scale=0.7]{F3x30}}~\Bigg],\nn\\
%&F_{4,p}(\bs{x},\bs{y})=\le Cs^2\Bigg[\raisebox{-18pt}{\includegraphics[scale=0.7]{F4x3}}+\raisebox{-18pt}{\includegraphics[scale=0.7]{F4x30}}~\Bigg],\nn\\
%&F_{7,p}(\bs{x},\bs{y})\le \raisebox{-18pt}{\includegraphics[scale=0.7]{F7x3}},~~
%F_{8,p}(\bs{x},\bs{y})\le Cs\raisebox{-18pt}{\includegraphics[scale=0.7]{F8x3}},\nn\\
%&F_{9,p}(\bs{x},\bs{y})\le Cs^2\raisebox{-18pt}{\includegraphics[scale=0.7]{F9x3}},~~F_{10,p}(\bs{x},\bs{y})\le Cs\Bigg[\raisebox{-18pt}{\includegraphics[scale=0.7]{F10x3}}+\raisebox{-18pt}{\includegraphics[scale=0.7]{F10x30}}\Bigg].\nn
%\end{align}
%\caption{Diagrammatic representations of $F_{1,p}(\bs{x},\bs{y})$ through $F_{10,p}(\bs{x},\bs{y})$ excluding $F_{5,p}(\bs{x},\bs{y})$ and $F_{6,p}(\bs{x},\bs{y})$, for $\bs{x}, \bs{y}$ with $\bs{x}(t)=\bs{y}(t)=3$.}
%\label{fig:B}
%\end{figure}

\section*{Acknowledgements}
The work of the author was supported by the National Center for Theoretical Sciences, Mathematics Division, Taiwan. 
The author thanks the Issac Newton Institute for Mathematical Sciences, Cambridge, for its hospitality during the programme {\it Stochastic systems for anomalous diffusion}, where the part of work on this paper was carried out.


\begin{thebibliography}{99}
\bibitem{B77}J. Blease
\newblock Directed bond percolation on hypercubical lattices. 
\newblock \emph{J. Phys. C}, \textbf{10} (1977): 925--936.


\bibitem{B772}J. Blease.
\newblock Pair connectedness for directed bond percolation on some $2d$ lattices by series methods. 
\newblock \emph{J. Phys. C}, \textbf{10} (1977): 3461--3476.


\bibitem{bhh21}D.C. Brydges, T. Helmuth and M. Holmes.
\newblock The continuous-time lace expansion.
\newblock \emph{Comm. Pure Appl. Math.}, \textbf{74} (2021): 2251--2309.

\bibitem{BK52}T.H. Berlin and M. Kac. 
\newblock The spherical model of a ferromagnet. 
\newblock \emph{Phys Rev}, \textbf{86} (1952): 821--835.


\bibitem{CD83}J. T. Cox and R. Durrett.
\newblock Oriented percolation in dimensions $d\ge 4$: bounds and asymptotic formulas.
\newblock \emph{Math. Proc. Camb. Phi. Soc.}, \textbf{93} (1983): 151--162.



\bibitem{GF74}P.R. Gerber and M.E. Fisher. 
\newblock Critical temperatures of classical $n$-vector models on hypercubic lattices. 
\newblock \emph{Phys Rev B}, \textbf{10} (1974), 4697--4703.

\bibitem{G10}B.T. Graham. 
\newblock Borel-type bounds for the self-avoiding walk connective constant. 
\newblock \emph{J.Phys. A, Math. Theor.}, \textbf{43} (2010), 13, Id/No235001.


\bibitem{HS95}T. Hara.
\newblock Mean field critical behavior for correlation length for percolation in high dimensions. 
\newblock \emph{Probab. Theory Relat. Fields}, \textbf{86} (1990): 337--385.

\bibitem{HS93}T. Hara and G. Slade.
\newblock Unpublished note
\newblock (1993).

\bibitem{HS95}T. Hara and G. Slade.
\newblock The Self-Avoiding Walk and Percolation Critical Points in High Dimensions.
\newblock \emph{Comb. Probab. Commun.}, \textbf{4} (1995): 197--215.



\bibitem{HS04}R. van der Hofstad and A. Sakai.
\newblock Gaussian scaling for the critical spread-out contact process above the upper critical dimension. 
\newblock \emph{Electron. J. Probab.}, \textbf{9}, no.9 (2004): 710--769.


\bibitem{HS05}R. van der Hofstad and A. Sakai.
\newblock Critical points for spread-out self-avoiding walk, percolation and the contact process above the upper critical dimensions.
\newblock \emph{Probab. Theory Relat. Fields}, \textbf{132} (2005): 438--470.



\bibitem{HGS05}R. van der Hofstad and G. Slade.
\newblock Asymptotic expansion in $n^{-1}$ for percolation critical values on the $n$-cube and $\mathbb Z^n$.
\newblock \emph{Random Struct. Algorithm}, \textbf{27} (2005): 331--357.

\bibitem{HM22}M. Hydenreich and K. Matzke.
\newblock Expansion for the critical point of site percolation: the first three terms.
\newblock \emph{Comb. Probab. Commun.}, \textbf{31} (2022): 430--454.

\bibitem{L97}T.M. Ligget.
\newblock \emph{THE 1996 WALD MEMORIAL LECTURES.}
\newblock \emph{Ann. Probab.}, \textbf{25}, No.1 (1997): 1--29.

\bibitem{ms93}N. Madras and G. Slade.
\newblock \emph{The Self-Avoiding Walk}, Birkh\"auser, Boston, MA,  Boston, Inc., Boston, (1993). xiv+425pp.


\bibitem{NY93}B.G. Nguyen and W.-S. Yang. 
\newblock Gaussian limit for critical oriented percolation in high dimensions. 
\newblock \emph{J. Stat. Phys.}, \textbf{78} (1995): 841--876.

\bibitem{NY932}B.G. Nguyen and W.-S. Yang. 
\newblock Triangle condition for oriented percolation in high dimensions.
\newblock \emph{Ann. Probab.}, \textbf{21}, No.4 (1993): 1809--1844.


\bibitem{S01}A. Sakai. 
\newblock Mean-field critical behavior for the contact process. 
\newblock \emph{J. Stat. Phys.}, \textbf{104} (2001): 111--143.

\bibitem{S018}A. Sakai. 
\newblock Hyperscaling for oriented percolation in $1+1$ space-time dimensions. 
\newblock \emph{J. Stat. Phys.}, \textbf{171} (2018): 462--469.


\bibitem{g004}G. Slade.
\newblock \emph{The Lace Expansion and its Applications}, Springer, (2006). Lecture Notes in
 Mathematics Vol. 1879. Ecole d’Eté de Probabilités de Saint–Flour XXXIV–2004. 

\bibitem{S79}A. D. Sokal. 
\newblock An improvement of Watson's theorem on Borel summability
\newblock \emph{J. Math. Phys.}, \textbf{21} (1979): 261--263.







\end{thebibliography}
\end{document}